\documentclass[12pt]{article}

\usepackage[T1]{fontenc}
\usepackage{lmodern}
\usepackage{amsthm,amsmath,amssymb, latexsym, mathtools,secdot}
\usepackage{bm}
\usepackage{bbm}
\usepackage{mathrsfs}

\usepackage[dvipsnames]{xcolor}
\usepackage[shortlabels]{enumitem}
\usepackage{microtype}
\usepackage[colorlinks=false]{hyperref}
\usepackage{bookmark}
\hypersetup{
    breaklinks=true,
    colorlinks,
		linktocpage,
    linkcolor={blue!80!black},
    citecolor={blue!80!black},
    urlcolor={blue!80!black}
}
\usepackage{stmaryrd}
\usepackage[nameinlink,capitalize]{cleveref}

\newif\ifdraft
\draftfalse
\ifdraft
\usepackage[notcite,notref]{showkeys}

\fi

\numberwithin{equation}{section}

\ifdraft
\usepackage[textsize=tiny]{todonotes}
\else
\usepackage[disable]{todonotes}
\fi

\newcommand{\vertiii}[1]{{\left\vert\kern-0.25ex\left\vert\kern-0.25ex\left\vert #1\right\vert\kern-0.25ex\right\vert\kern-0.25ex\right\vert}}

\newcommand*{\Gref}[1]{\hyperref[fSDE]{Eq($#1$)}}
\newcommand*{\Gtag}[1]{\tag{Eq($#1$)}}

\newcommand{\chain}{{\mathscr{C}}}
\newcommand{\W}{\mathring{W}}

\newcommand{\cC}{\mathcal{C}}

\newcommand{\A}{\mathcal{A}}
\newcommand{\B}{\mathcal{B}}
\newcommand{\F}{\mathcal{F}}
\newcommand{\G}{\mathcal{G}}
\renewcommand{\H}{\mathcal{H}}

\def \N{\mathbb{N}}

\def \R{\mathbb{R}}
\def \Z{\mathbb{Z}}

\DeclareMathOperator{\I}{\mathbbm{1}}
\newcommand{\indic}{\I}

\DeclareMathOperator{\Law}{Law}%

\newcommand{\nn}{\nonumber}

\def\E{\hskip.15ex\mathsf{E}\hskip.10ex}
\def\P{\mathsf{P}}
\DeclareMathOperator{\Var}{Var}
\DeclareMathOperator{\var}{-var}

\newcommand{\bra}[1]{{\llbracket #1 \rrbracket}}

\usepackage{aliascnt}
\theoremstyle{plain}
\newtheorem{theorem}{Theorem}[section]
\newcommand{\sharedthm}[2]{%
  \newaliascnt{#1}{theorem}%
  \newtheorem{#1}[#1]{#2}%
  \aliascntresetthe{#1}%
  \crefname{#1}{#2}{#2s}%
}
\sharedthm{lemma}{Lemma}
\sharedthm{proposition}{Proposition}
\sharedthm{corollary}{Corollary}
\theoremstyle{remark}
\sharedthm{definition}{Definition}
\sharedthm{notation}{Notation}
\sharedthm{remark}{Remark}
\sharedthm{example}{Example}
\sharedthm{condition}{Condition}
\sharedthm{assumption}{Assumption}
\newtheorem*{assumption*}{Assumption}
\def \defby{\vcentcolon =}

\newcommand{\dd}{\, d}
\newcommand{\from}{\colon}

\newcommand{\1}{\indic}

\newcommand{\Rd}{\R^d}
\newcommand{\CC}{\cC}
\newcommand\tand{\quad\text{and}\quad}
\newcommand\twhen{\quad\text{when}\quad}
\newcommand\tas{\quad\text{a.s.}}
\newcommand{\vr}{\var}
\newcommand{\les}{\lesssim}
\DeclarePairedDelimiter{\norm}{\lVert}{\rVert}
\DeclarePairedDelimiter{\abs}{\lvert}{\rvert}
\DeclarePairedDelimiterX{\inp}[2]{\langle}{\rangle}{#1, #2}
\DeclarePairedDelimiterX\set[1]\lbrace\rbrace{ #1 }
\renewcommand{\epsilon}{\varepsilon}
\newcommand{\eps}{\epsilon}

\newcommand{\wt}{\widetilde}

\begin{document}

\title{Strong Well-Posedness and Path-by-Path Uniqueness for Fractional SDEs with Integrable Drifts}

	\author{
	Oleg Butkovsky%
	\thanks{Weierstrass Institute, Mohrenstrasse 39, 10117 Berlin, FRG. Email: \texttt{oleg.butkovskiy@gmail.com}}
	\and
	\setcounter{footnote}{3}
	Khoa  L\^e
	\thanks{School of Mathematics, University of Leeds, U.K. Email: \texttt{k.le@leeds.ac.uk}}
	\and
	Toyomu Matsuda%
	\thanks{formerly
    EPFL, Bâtiment MA, 1015 Lausanne, Switzerland,
    Email: \texttt{toyomu.matsuda@epfl.ch}
	}
}

\maketitle

\begin{abstract}
  In dimension $d$, we consider stochastic differential equations driven by additive fractional Brownian motion with Hurst parameter $H<1/2$ and with a time-dependent drift that is $p$-integrable in space and $q$-integrable in time.
  Under the subcriticality condition $\frac{Hd}{p}+\frac1q<1-H$, complemented by two integrability conditions needed for Girsanov's theorem, we prove strong existence and path-by-path uniqueness of solutions.
  We also establish stability with respect to the initial data and the drift, both pathwise and in moment norms, and construct a H\"older continuous semiflow.
  When $H$ approaches $1/2$, our conditions reduce, up to an endpoint case, to the condition of Krylov and R\"ockner (2005) for Brownian motion, and in dimension one they complement the results of Nualart and Ouknine (2003), which rely on comparison arguments.
  Our methods combine the stochastic sewing lemma, Lyons' rough path theory, and the quantitative John--Nirenberg inequality for stochastic processes of vanishing mean oscillation.
\bigskip



\end{abstract}

\tableofcontents

\section{Introduction}
Given an integer dimension $d\ge1$, a finite time horizon $T>0$, a time-dependent vector field $b: [0, T] \times \mathbb{R}^d \to \mathbb{R}^d$ and a driving signal $\gamma: [0, T] \to \mathbb{R}^d$, we consider the following integral equation:
\begin{align}\label{eq:intro_de}
  x_t = x_0 + \int_0^t b(r, x_r) \dd r + \gamma_t, \quad t \in[0,T].
\end{align}
When $b$ is Lipschitz, well-posedness of \eqref{eq:intro_de} can be easily established using Picard's fixed-point argument. Interestingly, in the case where $\gamma$ is sampled from a stochastic process, well-posedness of \eqref{eq:intro_de} can be restored even for irregular $b$, a phenomenon known as \emph{regularization by noise}.

If $\gamma$ is sampled from Brownian motion $W$, Equation~\eqref{eq:intro_de} can be interpreted as the stochastic differential equation (SDE)
\begin{align}\label{eq:intro_sde_bm}
  d X_t = b(t, X_t) \, d t + d W_t, \quad X_0 = x_0.
\end{align}
The study of \eqref{eq:intro_sde_bm} for irregular drift $b$ has a long history and can be traced back to the classical works of Zvonkin~\cite{zvonkin74} and Veterennikov~\cite{veretennikov81}. In particular, the problem of strong well-posedness (existence of strong solutions and pathwise uniqueness) has been extensively investigated. 
The seminal work by Krylov and Röckner in \cite{krylov05} obtained  strong well-posedness for \eqref{eq:intro_sde_bm} with 
$b$ merely belonging to $L_q([0, T]; L_p(\mathbb{R}^d))$ for $p,q,d$ satisfying the condition
\begin{align}\label{eq:parameter_for_BM}
  \frac{d}{2 p} + \frac{1}{q} < \frac{1}{2}, \quad p \geq 2.
\end{align}
The case when $p=q$ was considered earlier by Gy\"{o}ngy and Mart\'{\i}nez in \cite{MR1864041}.

The groundbreaking results of \cite{MR1864041,krylov05} have inspired further research in this field. Zhang extended the analysis to include multiplicative noise in \cite{zhang10}. More recently, the critical case has been addressed by Röckner and Zhao in \cite{Rockner:2023aa, rocknerzhao21} and independently by Krylov in \cite{Krylov:2021aa}. Other directions of study include stability and numerical approximations with sharp rates, as discussed in \cite{leling21, galeatiling22, le22,dareiotis2021quantifying}.

Davie \cite{davie07} brought a new insight into the field, 
by considering uniqueness for each trajectory of the driving signal.
Specifically, it is shown that given a bounded measurable function $b$ on $\mathbb{R}^d$ 
the integral equation   
\begin{align*}
  x_t = x_0 + \int_0^t b(x_r) \dd r + \gamma_t, \quad t \geq 0,
\end{align*}
has a unique solution
for almost all realization $\gamma$ of Brownian motion.
Uniqueness proven in Davie's result is stronger than pathwise uniqueness, in that uniqueness is understood 
in the sense of ordinary integral equation for each realization; in particular, solutions are not assumed to be adapted stochastic processes.
This stronger notion of uniqueness is nowadays called \emph{path-by-path uniqueness}.
Recently, path-by-path uniqueness under Krylov--Röckner condition \eqref{eq:parameter_for_BM} has been  
obtained in \cite{anzeletti2023pathbypath}.

Regularization by noise is not solely specific to Brownian noise. 
For SDEs, two major classes of noise have been considered. One is Lévy noise, as studied in \cite{Kim:2014aa,Athreya:2020aa,Kremp:2022aa,Raynal:2022aa} to name a few.
The other is \emph{fractional Brownian} noise, which is the focus of this article.
When the driving signal in \eqref{eq:intro_de} is sampled from fractional Brownian motion $B^H$ with Hurst parameter $H \in(0,1)$, 
one leads to the SDE
\begin{equation}
  \label{eq:main-sde}
  d X_t = b(t, X_t) \, d t + d B^H_t, \quad X_0 = x\in\Rd.
\end{equation}
Such SDE has a unique strong solution whenever the drift $b$ is Lipschitz continuous, as one can easily see from standard fixed-point arguments.
The case when $b$ is irregular was considered first by Nualart and Ouknine in\cite{Nualart:2003tz}, surprisingly prior to \cite{krylov05}. Building upon their earlier work \cite{Nualart:2002uv}, Nualart and Ouknine showed regularization by fractional noise for SDEs with integrable drifts although only in dimension one.
Their arguments are based on Girsanov's theorem and the comparison principle, the 
latter being specific to dimension one. 
It is worth noting that the arguments from \cite{Nualart:2003tz} also works for the case of driving Brownian noise.
Yet, due to the use of the comparison principle, their arguments cannot be applied to higher dimensional settings.
As one will see in the current article, proving regularization by fractional noise in multi-dimensions 
requires more robust methods, which forego known techniques in both one-dimensional and Markovian settings.

Different approaches have been developed for regularization by fractional Brownian noise in multi-dimensions.  
Catellier and Gubinelli in \cite{catelliergubinelli}, building upon \cite{davie07}, introduce a path-by-path approach lying on the framework of nonlinear Young differential equations.  
Baños, Nilssen and Proske in \cite{Banos:2020vv} employ tools from Malliavin calculus to construct strong solutions. In \cite{khoa20}, the second named author  initiates the \emph{stochastic sewing} techniques and proves the strong well-posedness of {\eqref{eq:main-sde}} under the condition
\begin{align}\label{con.kle}
  \frac{Hd}{p} + \frac{1}{q} < \frac{1}{2} - H, \quad p \geq 2, \quad q>2.
\end{align}

These encouraging works have facilitated recent dramatic developments in the field: further results on singular fractional SDEs \cite{Harang:2021aa, Gerencser:2023aa,Galeati:2022aa,Galeati:2022ab,Banos:2022ut,anzeletti22,galeatigerencser22,gassiatmardy22,li2023nonstationary,BLMytnik}, 
stochastic partial differential equations with singular drift \cite{Butkovsky:2019aa,athreya2021wellposedness,Bechtold:2023aa,HRZ}, approximation of singular SDEs 
\cite{Butkovsky:2021aa,dareiotis2021quantifying}, 
regularization by multiplicative noise \cite{dareiotisgerencser22,catellierduboscq22} 
and extensions of the stochastic sewing lemma \cite{Li:2022aa,Gerencser:2023aa,Le:2023aa,kern2023stochastic,matsuda22,BLMytnik}.

In the related work \cite{galeatigerencser22} by Galeati and Gerencsér, using stochastic sewing techniques, the authors show strong well-posedness and path-by-path uniqueness for Equation \cref{eq:main-sde} where the drift belongs to $L_q([0, T]; \CC^\alpha_x(\mathbb{R}^d))$  
for $\alpha\in\R$ and $q\in(1,\infty]$ satisfying
\begin{align}\label{con.gg}
  -H \alpha+\frac1{\min(q,2)} <1-H.
\end{align}
In view of the embedding $L_p(\Rd) \hookrightarrow \CC^{-d/p}(\Rd)$, their results allow drifts in $L_q([0, T]; L_p(\mathbb{R}^d))$ with $p\in[1,\infty]$, $q\in(1,\infty]$ satisfying 
\begin{align}\label{con.ggpq}
  \frac{Hd}{p}+\frac1{\min(q,2)} <1-H.
\end{align}
We make an important observation that neither \eqref{con.kle} nor \eqref{con.ggpq} is optimal because they fail to recover \eqref{eq:parameter_for_BM} for $H = 1/2$.
In the companion work \cite{BLMytnik}, weak existence for the SDE \eqref{eq:main-sde} with time-homogeneous drift 
$b = b(x)$ is obtained under the optimal condition $H \in (0, 1)$ and $\frac{Hd}{p} < 1 - H$. 
Weak existence for time-dependent drifts, even beyond the subcritical regime, was subsequently obtained by Butkovsky and Gallay \cite{BGal23}.
However, the problem of uniqueness was left for future works.
The main objective of the present work is to address the question of uniqueness for \eqref{eq:main-sde} when $H < 1/2$.
Precisely, we prove strong existence and path-by-path uniqueness for \eqref{eq:main-sde} under  conditions which are
compatible with the Brownian case treated by Krylov and Röckner (\cref{thm:main_uniqueness}), together with stability estimates (\cref{thm:main_stability}) and the existence of a H\"older continuous semiflow (\cref{thm:main_semiflow}).
Our key estimates on additive functionals (\cref{subsec:add_fc_fbm}) are proved for general Gaussian processes satisfying a local nondeterminism condition; the exponential integrability and linearization results (\cref{subsec:girsanov,subsec:linearization}) are stated for the fractional Brownian motion and the Riemann--Liouville process, the latter being needed in our conditioning arguments.
For path-by-path uniqueness, we follow the strategy of \cite{anzeletti2023pathbypath}, where the Brownian case under the Krylov--Röckner condition is treated.

The rest of the paper is organized as follows. In \cref{sec:main_results} we state our main results, \cref{thm:main_uniqueness,thm:main_stability,thm:main_semiflow}, discuss the conditions on the parameters and outline the strategies of the proofs. \cref{s:3} collects preliminary material: auxiliary estimates for stochastic processes, Girsanov's theorem for fractional Brownian motion, and linear Young differential equations. In \cref{sec:int_bd} we establish the key bounds on additive functionals of Gaussian processes, including their variational regularity and exponential integrability, together with a linearization identity. \cref{sec:strong.existence} is devoted to pathwise stability and the existence of strong solutions. In \cref{sec:stability} we derive stability in moment norms by means of VMO estimates and a rough path lift of the linearized equation. \cref{sec:flow_and_path_by_path_uniqueness} constructs a semiflow and deduces path-by-path uniqueness, and \cref{sec:proof_main_results} completes the proofs of the main results. The approximation of the drift by smooth functions is discussed in \cref{app:aux}.

\medskip 
\textbf{Acknowledgements}. The main results of this paper were announced by KL at the Oberwolfach workshop \emph{Directions in Rough Analysis} in November 2024, see \cite{OWR2024LBM}. The authors are very grateful Lucio Galeati, M\'at\'e Gerencs\'er and Hannes Kern for many interesting discussions and valuable comments on an earlier draft of the paper. OB is  funded by the Deutsche Forschungsgemeinschaft (DFG, German Research Foundation) under Germany's Excellence Strategy --- The Berlin Mathematics Research Center MATH+ (EXC-2046/1, project ID: 390685689, sub-project EF1-22). KL is funded by the Engineering \& Physical Sciences Research Council (EPSRC), grant
number EP/Y016955/1. The main part of the research was done while TM was a PhD student at Freie Universität Berlin with financial support by the DFG through the IRTG 2544.

\medskip
\textbf{AI declaration.} 
Near the completion of this work, AI model Sonnet 5.5. was used for proofreading and checking consistency. 
The final version of the paper was carefully reviewed and edited by the authors, and all mathematical content was written by the authors.

\section{Main results}\label{sec:main_results}

We begin by introducing  the basic notation. Fix a dimension $d\in\N \defby \{1, 2, \ldots\}$ and a time interval $T>0$. Let $B^H=\{B_t^H, t\in[0,T]\}$ be a $d$-dimensional fractional Brownian motion with the Hurst parameter $H\in(0,1)$. That is, $B^H$ is a $d$-dimensional vector, whose components are independent and each of them is a mean zero Gaussian process with covariance
$$
\E B_s^{H,i} B_t^{H,i}=\frac12(t^{2H}+s^{2H}-|t-s|^{2H}),\quad s,t\in[0,T],\,\, i=1,\hdots, d.
$$

We recall (see, e.g., \cite[Section~5.1.3]{nualart06}) that given a fractional Brownian motion $B^H$ one can construct on the same probability space a standard $d$-dimensional Brownian motion $W$ such that
\begin{equation}\label{eq:fbm_kernel} 
B^H_t = \int_0^t K_H (t, s) \, d W_s, 
\end{equation}
where the kernel $K_H$ is given by 
\begin{align}
		 K_H (t, s) &:= c_H \Big[ \Big( \frac{s}{t} \Big)^{\frac{1}{2}-H} 
		(t-s)^{H-\frac{1}{2}}\notag \\
		&\quad+ 
		\Big(\frac{1}{2} -H\Big) s^{\frac{1}{2} - H} 
		\int_s^t u^{H-\frac{3}{2}} (u-s)^{H-\frac{1}{2}} \,d u
		\Big],\qquad \text{for }H<\frac12; \label{eq:kernel_def}\\
		 K_H (t, s) &:= c_H s^{\frac12-H} 
		\int_s^t u^{H-\frac{1}{2}} (u-s)^{H-\frac{3}{2}} \,d u,\qquad \text{for }H>\frac12. \notag
\end{align}
The constant $c_H$ is chosen so that $\E[|B_t^{H, i}-B_s^{H, i}|^2] = \abs{t-s}^{2H}$. 

We say that  an adapted process $(B^H_t)_{t \geq 0}$ is  an $(\mathcal{F}_t)$-\emph{fractional Brownian motion} (starting at $0$) if there exists a standard $d$-dimensional $(\mathcal{F}_t)_{t \geq 0 }$-Brownian motion $W$ such that \eqref{eq:fbm_kernel} holds. Since the natural filtrations generated by $W$ and $B^H$ coincide, 
for which we write $(\F_t^{B^H})$,
we see that $B^H$ is always $(\F_t^{B^H})$-fractional Brownian motion. We refer the reader to \cite[Section~5]{nualart06}, \cite{picard11,nourdin12} for more details on 	fractional Brownian motion.

As usual, by a \textit{weak solution} to  \eqref{eq:main-sde}, we mean a couple of  continuous processes $(X,B^H)$ on a complete filtered probability space $(\Omega, \F, (\mathcal{F}_t)_{t \in [0, T]}, \P)$ such that $B^H$ is an $(\mathcal{F}_t)$-fractional Brownian motion, $X$ is adapted to $(\mathcal{F}_t)$, $\int_0^T \abs{b(r, X_r)} \dd r < \infty$ and 
\begin{equation}\label{fSDE}	\Gtag{x;b}
X_t = x + \int_0^t b(r,X_r)\,dr + B^H_t, \quad t \in [0, T].
\end{equation}
We say that a weak solution $(X,B^H)$ is a \textit{strong solution}  to  \eqref{eq:main-sde}, if it is adapted to the standard augmentation of filtration generated by $B^H$.

We say that \textit{weak uniqueness} holds for \eqref{eq:main-sde} if the law of $X$ is unique. 
We say that \textit{pathwise uniqueness} for \eqref{eq:main-sde} holds if for any two weak solutions of  \eqref{eq:main-sde} $(X^1,B^H)$, $(X^2,B^H)$  defined on the same filtered probability space with the same $B^H$  and the same initial condition one has $\P(X^1_t=X^2_t\, \text{for all }t\in[0,T])=1$. Finally, we say that  \textit{path-by-path uniqueness} holds if for any probability space $(\Omega, \F,  \P)$ on which $B^H$ is defined, there exists a set of full measure $\Omega^*\subset\Omega$ such that 
for any $x\in\R^d$, $\omega\in\Omega^*$ there exists a unique solution to \eqref{fSDE}.

Clearly, path-by-path uniqueness is the strongest form of uniqueness for an equation. We refer to \cite[Appendix~B]{anzeletti2023pathbypath} and \cite{ShaWr} for further discussions of the relations between the three notions of uniqueness.

Our first main result is strong well-posedness and path-by-path uniqueness of SDE \eqref{eq:main-sde} 
under the conditions compatible with the Brownian case of Krylov--Röckner \cite{krylov05}. 

Let $p,q\in[1,\infty]$. Let $L_qL_p([0,T]):= L_q([0,T]; L_p(\R^d))$ be  the space of measurable functions 
$b \from [0, T] \times \mathbb{R}^d \to \mathbb{R}^d$ such that the following mixed norm
\begin{equation*}
\|b\|_{L_q([0,T]; L_p(\R^d))} := \Bigl( \int_0^T \bigl( \int_{\mathbb{R}^d} |b(t,x)|^p \,dx \bigr)^{\frac{q}{p}} \,dt \Bigr)^{\frac{1}{q}}
\end{equation*}
is finite. Here we make the usual modification if $p = \infty$ or $q = \infty$.
When the temporal domain is clear from the context, we simply write $L_qL_p$ for $L_qL_p([0,T])$.

We will assume the following conditions for our main results.
\begin{subequations}\label{eq:parameters}
\begin{equation}
    p, q\in[1,\infty], \quad H \in (0,1/2),
    \label{eq:parameters_range}
\end{equation}
\begin{equation}\label{eq:parameters_subcritical}
    \frac{Hd}{p} + \frac{1}{q} < 1 - H,
\end{equation}
\begin{equation}\label{eq:parameters_girsanov}  
  p \geq 2 Hd,
\end{equation}
\begin{equation}\label{eq:parameters_cameron_martin}
    \frac{Hd}{p} < \frac{p^{-1} - q^{-1}}{p^{-1} - 2^{-1}} \Big( \frac{1}{2} - H \Big) 
    \quad \text{if } p \le (1-H)^{-1} \text{ and } q > 2. 
\end{equation}
\end{subequations}

We will clarify the purpose of these conditions in \cref{subsec:parameters}. 
Note that, by substituting $H = 1/2$ to \cref{eq:parameters_subcritical,eq:parameters_girsanov,eq:parameters_cameron_martin}, we recover the Krylov--Röckner condition \eqref{eq:parameter_for_BM}, up to the endpoint case $p=2$ (see \cref{subsec:parameters}).

\begin{theorem}[Strong existence and path-by-path uniqueness]
  \label{thm:main_uniqueness}
Assume that the parameters $H$, $d$, $p$, $q$ satisfy \eqref{eq:parameters}. Let $b\in L_q([0,T]; L_p(\R^d))$. Then the following holds:
  \begin{enumerate}[label={(\roman*)}, font=\normalfont, leftmargin=*]
    \item \emph{(Strong existence)} For each $x \in \mathbb{R}^d$, there exists a stochastic process 
            $X = (X_t)_{t \in [0, T]}$ adapted to the standard augmentation of the filtration generated by $B^H$ such that 
      almost surely  
      $\int_0^T |b(r,X_r)| \, d r < \infty$ and 
      \begin{equation*}
        X_t = x + \int_0^t b(r,X_r) \,d r + B^H_t, \quad \forall t \in [0, T].
      \end{equation*}
    \item \emph{(Path-by-path uniqueness)} There exists a set $\Omega^* \subset \Omega$ of full probability such that,
      for every $\omega \in \Omega^*$ and every $x \in \mathbb{R}^d$, the integral equation
      \begin{equation*}
        y_t = x + \int_0^t b(r,y_r) \,d r + B^H_t(\omega), \quad \forall t \in [0, T],
      \end{equation*}
      has exactly one continuous solution $y \colon [0, T] \to \mathbb{R}^d$ satisfying $\int_0^T |b(r,y_r)| \, d r < \infty$.
  \end{enumerate}
  \end{theorem}

In order to state our stability result, we need the following additional notation. 
For a function $f\in L_qL_p$, and $\eta>0$  consider the following quantity 
 \begin{equation}\label{eq:frak_c}
\Lambda_f(\eta ):=\inf\{ u>0: \|f \I_{\{|f| > u\}}\|_{L_qL_p} \le \eta \}. 
  \end{equation}
By the monotone convergence theorem, the map $u\mapsto\|f \I_{\{|f| > u\}}\|_{L_qL_p}$ is non-increasing and right-continuous for all $p,q\in[1,\infty]$. Being its generalized inverse, $\Lambda_f$ has the property that
$\|f \I_{\{|f| > \Lambda_f(\eta )\}}\|_{L_qL_p}\le \eta $ whenever $\Lambda_f(\eta)<\infty$.
If $q<\infty$, the dominated convergence theorem shows that $\|f \I_{\{|f| > u\}}\|_{L_qL_p}\to0$ as $u\to\infty$, and hence $\Lambda_f(\eta )$ is finite for any $\eta >0$.
\begin{remark}\label{rem:Lambda_q_infty}
  If $q=\infty$, then $\Lambda_f(\eta)$ may be infinite: for instance, if $p<\infty$ and $f(t,x)=t^{-d/p}\I_{\{|x|\le t\}}$, then $\|f \I_{\{|f| > u\}}\|_{L_\infty L_p}=\|f\|_{L_\infty L_p}$ for every $u>0$.
  However, the conditions \eqref{eq:parameters} are open with respect to $1/q$: if they hold with $q=\infty$, then they also hold with $q$ replaced by some $\bar q\in(2,\infty)$, and $L_\infty L_p\hookrightarrow L_{\bar q}L_p$ on $[0,T]$.
    Hence, whenever $\Lambda_f$ appears below with $q=\infty$, it is understood that $L_qL_p$ in \eqref{eq:frak_c} is replaced by $L_{\bar q}L_p$ for such a fixed $\bar q$. In particular, $\Lambda_f(\eta)$ is always finite.

  The same reduction is used throughout \cref{subsec:girsanov,subsec:refined_pathwise_estimates,sub:existence_strong_sol,sec:flow_and_path_by_path_uniqueness}: the statements there are either qualitative (exponential moments, almost sure bounds, identification of limits) or involve a H\"older exponent $\gamma<1-H-\frac{Hd}p-\frac1q$, and for $q=\infty$ they follow from the corresponding statements with a finite $\bar q$, chosen so large that \eqref{eq:parameters} holds with $\bar q$ and $\gamma<1-H-\frac{Hd}p-\frac1{\bar q}$, since $\|f\|_{L_{\bar q}L_p}\le T^{1/\bar q}\|f\|_{L_\infty L_p}$. Accordingly, in the proofs of these parts we assume $q<\infty$ without further mention. The results of \cref{subsec:add_fc_fbm,sec:stability}, which involve the sharp exponent $1-H-\frac{Hd}p-\frac1q$, are proved directly for $q=\infty$ as well.
\end{remark}
  
In the theorem below and throughout the article, $\mathcal{B}^{-1}_p := \mathcal{B}_{p, \infty}^{-1}$ is the Besov space on $\mathbb{R}^d$ of regularity $-1$ and integrability $(p,\infty)$.
\begin{theorem}[Stability]\label{thm:main_stability}
  Suppose that the condition \eqref{eq:parameters} is satisfied.
  For each $i\in\{1,2\}$, let $x^i\in \Rd$ and $b^i\in L_q([0, T]; L_p(\mathbb{R}^d))$, and 
  let $X^i$ be the solution to \hyperref[fSDE]{\emph{Eq(}$x^i; b^i$\emph{)}}.
  Then we have the following stability estimates.
  \begin{enumerate}[label={(\roman*)}, font=\normalfont, leftmargin=*]
    \item \emph{(Pathwise stability)} 
    Suppose that $p > 2 H d$. There exist some non-negative random variables $A$ and $D$ such that for any $m \in (0, \infty)$ with some 
    \begin{equation}\label{eq:C_dep}
    C = C(d, H, p, q, T, m, \norm{b^1}_{L_qL_p}, \norm{b^2}_{L_qL_p})    
    \end{equation}  
    we have 
    \begin{equation}\label{eq:A_and_D_bd}
            \norm{A}_{L_m(\Omega)} \leq C (1 + \norm{b^1}_{L_qL_p})^2,  
      \quad 
      \norm{D}_{L_m(\Omega)} \leq C  \norm{b^1 - b^2}_{L_q \mathcal{B}_p^{-1}},
    \end{equation}
    and 
    \begin{equation}\label{eq:X_1_X_2_pw_bd}
      \sup_{t \in [0, T]} \abs{X^1_t - X^2_t} 
      \leq e^{A} (\abs{x^1 - x^2} + D) \quad\text{almost surely}.
    \end{equation}
    If $p = 2 H d$, we similarly have non-negative random variables $A$ and $D$ such that for any $m \in (0, \infty)$ with some 
    $\eta  = \eta (d, H, p, q, T, m)$ and 
    \begin{equation}\label{eq:C_dep_2}
      C = C(d, H, p, q, T, m, \Lambda_{b^1}(\eta ), \Lambda_{b^2}(\eta ))
    \end{equation}
    we have \eqref{eq:A_and_D_bd} and \eqref{eq:X_1_X_2_pw_bd}.    
   
    \item \emph{(Stability in moment norms)} 
      If $p > 2 H d$, for every $m \in (0, \infty)$ we have a constant $C$ whose dependence is as in \eqref{eq:C_dep} such that
      for every $0 \leq s < t \leq T$ we have
      \begin{equation}\label{eq:moment_stab}
        \norm{X^1_t - X^1_s - (X^2_t - X^2_s)}_{L_m(\Omega)} 
        \leq C (\abs{x^1 - x^2} + \norm{b^1-b^2}_{L_q \mathcal{B}_p^{-1}})
        (t-s)^{1 - H - \frac{Hd}{p} - \frac{1}{q}}.
      \end{equation}
            Similarly, if $p = 2 H d$, for every $m \in (0, \infty)$ we have an $\eta  = \eta (d, H, p, q, T, m)$ and a $C$ whose dependence is as in \eqref{eq:C_dep_2} such that \eqref{eq:moment_stab} holds with $C$ replaced by $Ce^{C(1\vee\|b^1\|_{L_qL_p})^{1/\gamma^2}}$, where $\gamma\defby1 - H - \frac{Hd}{p} - \frac{1}{q}$.
  \end{enumerate} 
\end{theorem}
\begin{remark}
  When $p=2Hd$, it is not expected to have uniform stability estimates with respect to the norms of $b^1,b^2$. This phenomenon also happens for SDEs with Brownian motion and drift in $L_d(\Rd)$, see \cite[Theorem 3.7]{Krylov:2021aa}.
\end{remark}

\begin{definition}\label{def:semiflow}
    A measurable map  
    \begin{align*}
      \Phi \from \set{(s, t, x, \omega) \in [0, T]^2 \times \mathbb{R}^d \times \Omega : s \leq t} 
      \ni (s, t, x, \omega) \mapsto \Phi^{s, x}_t(\omega)
    \end{align*}
    is called a \emph{semiflow} to \eqref{eq:main-sde} if for almost all $\omega$ we have the following property. 
    \begin{enumerate}[label={(\roman*)}, font=\normalfont, leftmargin=*]
      \item For each $s \leq t$ and $x$, the map $\omega \mapsto \Phi^{s, x}_t(\omega)$ is measurable with respect to the completion of $\sigma(B^H_r:r\le t)$. 
      \item For $s < u < t$ and $x \in \mathbb{R}^d$, we have $\Phi^{s, x}_t(\omega) = \Phi^{u, y}_t(\omega)$ 
      with $y \defby \Phi^{s, x}_u(\omega)$.
      \item For each $s \geq 0$ and $x \in \mathbb{R}^d$, the path $\phi^{s, x}_t \defby \Phi^{s, x}_t(\omega)$ satisfies 
      $\int_s^T \abs{b(r, \phi^{s, x}_r)} \dd r < \infty$ and
      \begin{align*}
      \phi^{s, x}_t = x + \int_s^t b(r, \phi^{s, x}_r) \dd r + B^H_t(\omega) - B^H_s(\omega), \quad \forall t \in [s, T].
      \end{align*}
    \end{enumerate}

    A semiflow is called locally $\kappa$-Hölder continuous if $\Phi$ is continuous with respect to $(s, t, x)$  
    and if for every $K \in (1, \infty)$ and for almost all $\omega$ there exists 
    a finite constant $C = C(\kappa, K, \omega)$ such that  
    \begin{align*}
      \abs{\Phi^{s, x}_t(\omega) - \Phi^{s, y}_t(\omega)} \leq C \abs{x-y}^{\kappa}
    \end{align*}
    for every $s \leq t$ and $(x, y)$ with $\abs{x} \vee \abs{y} \leq K$.
  \end{definition}
That the null set in \cref{def:semiflow} does \emph{not} depend on $(s, x)$ is crucial. In the terminology of \cite[Section~2]{anzeletti2023pathbypath}, for almost every $\omega$, $(s,t,x)\mapsto\Phi^{s,x}_t(\omega)$ is then a semiflow generated by the integral equation \eqref{fSDE} driven by the path $B^H(\omega)$.
\begin{theorem}[Semiflow]\label{thm:main_semiflow}
  Under \eqref{eq:parameters}, there exists a semiflow $\Phi$ to \eqref{eq:main-sde} in the sense of \cref{def:semiflow} such that, for each $s < t$ and $x \in \mathbb{R}^d$, the random variable $\Phi^{s,x}_t$ is measurable with respect to the completion of $\sigma(B^H_r - B^H_s : r \in [s, t])$, and such that the following H\"older estimate holds for almost all $\omega$: for every small $\epsilon>0$ and every $K \in [1, \infty)$, there exists a $C = C(\omega, \epsilon, K)$ such that, for every $(s, t, x)$ and $(s', t', x')$ with $\abs{x} \vee \abs{x'} \leq K$ (where we set $\Phi^{s,x}_t\defby x$ for $t<s$),
  \begin{equation}\label{eq:semiflow_holder}
    \abs{\Phi_t^{s, x}(\omega) - \Phi^{s', x'}_{t'}(\omega)} \leq C (\abs{s-s'}^{H - \epsilon} + \abs{t-t'}^{H - \epsilon} + \abs{x-x'}^{1-\epsilon}).
  \end{equation}
  In particular, $\Phi$ is locally $\kappa$-H\"older continuous for every $\kappa\in(0,1)$.
\end{theorem}

\subsection{On parameters}\label{subsec:parameters}
We explain our conditions \eqref{eq:parameters}. 
Most important is the \emph{subcritical} condition \eqref{eq:parameters_subcritical} under which the fractional noise dominates at small scales. 
To reveal the subcritical condition, we follow an argument of Galeati and Gerencsér in \cite{galeatigerencser22}. 
If $X$ solves the SDE \eqref{eq:main-sde}, then the scaled process $X^{(\lambda)}_t \defby \lambda^{-H} X_{\lambda t}$
solves the SDE 
\begin{align}\label{eq:scaled_SDE}
  \mathrm{d} X^{(\lambda)}_t = \lambda^{1 - H - \frac{Hd}{p} - \frac{1}{q}} b^{(\lambda)}(t, X_t^{(\lambda)}) \,\mathrm{d} t 
  + \mathrm{d} B^{(\lambda)}_t,
\end{align}
where 
\begin{align*}
  b^{(\lambda)}(t, x) \defby \lambda^{\frac{Hd}{p} + \frac{1}{q}} b(\lambda t, \lambda^H x), 
  \quad B^{(\lambda)}_t \defby \lambda^{-H} B_{\lambda t}^H.
\end{align*}
We note that $\norm{b^{(\lambda)}}_{L_qL_p} = \norm{b}_{L_qL_p}$ and that $B^{(\lambda)}$ is equal to $B^H$ in law.
Domination of the noise at small scales entails that the order of the drift term 
is smaller than that of the driving noise as $\lambda$ vanishes. This enforces that
\begin{align*}
  1 - H - \frac{Hd}{p} - \frac{1}{q}>0,
\end{align*} 
which is equivalent to \eqref{eq:parameters_subcritical}.

Let us explain technical conditions \cref{eq:parameters_girsanov,eq:parameters_cameron_martin}. 
Both conditions are due to Girsanov's arguments. The condition \eqref{eq:parameters_cameron_martin} ensures that the drift term
\begin{align*}
  t \mapsto \int_0^t b(r, X_r) \dd r
\end{align*}
belongs to the Cameron--Martin space $\mathcal{H}$ of the fractional Brownian motion, and the condition \eqref{eq:parameters_girsanov} ensures that 
the Cameron--Martin norm of the drift satisfies the strong Novikov's condition, that is 
\begin{align*}
  \mathbb{E}\Big[ \exp \Big( \lambda \norm[\Big]{\int_0^{\cdot} b(r, X_r) \dd r}_{\mathcal{H}}^2 \Big)  \Big] < \infty , \quad \forall \lambda \geq 0.
\end{align*}

If $q = \infty$, which includes the case where $b$ is time-independent, our condition \eqref{eq:parameters} simplifies as follows. 
The condition \eqref{eq:parameters_cameron_martin} becomes
\begin{align}\label{eq:parameters_q_inf_cameron_martin}
  H + \frac{H d}{p} < \frac{1}{2} + \frac{H d}{2} \quad\text{whenever } p \le (1-H)^{-1}.
\end{align}
Hence, the condition \eqref{eq:parameters} reduces to $H \in (0, 1/2)$ and 
\begin{align*}
  p \in 
  \begin{cases}
    [1, \infty] \cap (2H/(1-H), \infty], & \text{if } d = 1, \\
    [\max\{1, 4H\}, \infty] \setminus \{4H\}, & \text{if } d =2, \\
        [\max\{1,2Hd\}, \infty], &\text{if } d \geq 3.
  \end{cases}
\end{align*}
In particular, when $d = 1$ and  $H < 1/3$, it allows $b \in L_\infty L_1$.

When $H = 1/2$, the conditions \eqref{eq:parameters_subcritical}--\eqref{eq:parameters_cameron_martin} become \eqref{eq:parameter_for_BM}, except for the endpoint $p=2$: indeed, \eqref{eq:parameters_subcritical} becomes $\frac{d}{2p}+\frac1q<\frac12$, which forces $q>2$ and $p>d$, and then \eqref{eq:parameters_cameron_martin}, written in the form $\frac{Hd}{p}(\frac1p-\frac12)<(\frac1p-\frac1q)(\frac12-H)$ (which is equivalent for $H<\frac12$, since then $p\le(1-H)^{-1}<2$), excludes $p\le2$ when $H=\frac12$. The endpoint $p=2$, which is only relevant for $d=1$, is thus not recovered.
In this sense, our results can be viewed as an extension of Krylov--Röckner \cite{krylov05} when the driving process is a fractional Brownian motion with $H<1/2$. 
In the case when $H>1/2$, we expect that an amount of positive regularity, say in Sobolev scales, is necessary for uniqueness of strong solutions. However, our arguments do not extend trivially to the cases $H\ge1/2$ and we leave this problem open for future investigations.

\subsection{Strategies of proofs}\label{subsec:strategy}
\textbf{Strong existence.} We approximate $b$ by smooth functions $b^n$, let $X^n$ be the solution to \Gref{x;b^n}, and show that $(X^n)$ converges in probability. As in \cite[Section~6]{khoa20} (motivated by \cite{MR1302388}), the differences solve the linearized equation
\begin{align}\label{intro.X_m_n}
	X^m_t - X^n_t = \int_0^t (X^m_r - X^n_r) \dd V^{m, n}_r + R^{m, n}_t,
\end{align}
where $V_t^{m, n} =\int_0^t\int_0^1\nabla b^m(r,\theta X^m_r+(1- \theta)X^n_r)\dd \theta \dd r$ and $R_t^{m, n} =\int_0^t(b^m-b^n)(r,X^n_r)\dd r$.
A Gr\"onwall-type argument for the Young equation \eqref{intro.X_m_n} applies as soon as $V^{m,n}$ and $R^{m,n}$ have finite $\ell$-variation for some $\ell<2$ (\cref{lem:linearYDE}). In \cite{khoa20}, $V^{m,n}$ was measured in H\"older norms, which leads to the non-optimal condition \eqref{con.kle}; here we work with $\ell$-variation instead.
By Girsanov's theorem, it suffices to consider the additive functional $U=\int_0^\cdot\nabla b^m(r,B^H_r)\dd r$. The key new estimate is \cref{L:firstb}: when bounding the germ $A_{s,t}=\E^s\int_s^tb^m(r,B^H_r)\dd r$ in the stochastic sewing lemma, we use not only the local nondeterminism of $B^H_r-\E^sB^H_r$, as in \cite{khoa20}, but also the randomness of the $\mathcal{F}_s$-measurable part $\E^sB^H_r$. Combined with the stochastic sewing lemma with controls \cite{Le:2023aa} and a Kolmogorov criterion in variation (\cref{lem:variation-kolmogorov}), this shows that $U$ has finite $\ell$-variation for every $\ell>(1-H-Hd(\frac1p-\frac1{p\vee2}))^{-1}$ (\cref{lem:integral_along_fbm}), and this lower bound is less than $2$ under \eqref{eq:parameters} (\cref{lem.minusoneOK}).
For $p<2$, the latter only requires $\frac{Hd}{p}+H<\frac12+\frac{Hd}{2}$, that is, \eqref{eq:parameters_cameron_martin} with $q=\infty$; it is an interesting question whether \eqref{eq:parameters_cameron_martin} can be replaced by this weaker condition in our main results.
Since $b$ is discontinuous, identifying the limit of $(X^n)$ as a solution requires the path-by-path regularization estimates of \cref{subsec:refined_pathwise_estimates}, inspired by \cite{anzeletti2023pathbypath}.

\textbf{Girsanov's theorem.} The above and the arguments below rely on Girsanov's theorem, for which we show the Novikov-type bound $\E\exp(\lambda\|\psi\|^2_{\mathcal{H}})<\infty$ for all $\lambda>0$, where $\psi=\int_0^\cdot b(r,X_r)\dd r$ and $\mathcal{H}$ is the Cameron--Martin space of $B^H$. Since $\mathcal{H}$ contains the Sobolev space $W^{H+\frac{1}{2} + \varepsilon}_2$, we bound exponential moments of Sobolev norms of $\psi$, combining the quantitative John--Nirenberg inequality for processes of vanishing mean oscillation \cite{le22} with Sobolev and Gagliardo--Nirenberg embeddings; this is where the conditions \eqref{eq:parameters_girsanov} and \eqref{eq:parameters_cameron_martin} enter (\cref{subsec:girsanov}).

\textbf{Stability and path-by-path uniqueness.} Following \cite{davie07,Shaposhnikov:2016aa,galeatigerencser22,anzeletti2023pathbypath}, path-by-path uniqueness follows from (i) pathwise regularization estimates for $(t,x)\mapsto\int_0^tb(r,B^H_r+x)\dd r$ (\cref{lem:fbm_regularizing}), which we deduce from moment bounds obtained by stochastic sewing and the John--Nirenberg inequality, and (ii) an almost Lipschitz semiflow, which we construct from stability estimates in $L_m(\Omega)$ (\cref{sec:flow_and_path_by_path_uniqueness}).
For two solutions $X^1,X^2$, the difference $X^\Delta=X^1-X^2$ solves an equation of the form \eqref{intro.X_m_n}, and the Gr\"onwall argument yields the pathwise bound (\cref{l:pathwise.stabitliy})
\begin{align*}
  \sup_t|X^\Delta_t|\le Ce^{C\|V\|^\ell_{\ell\var}}(|x^1-x^2|+\|R\|_{\ell\var}).
\end{align*}
This does not give moment bounds, since $\|V\|^\ell_{\ell\var}$ need not be exponentially integrable under \eqref{eq:parameters} (see the beginning of \cref{sec:stability}), and the techniques of \cite{cass2013,galeatigerencser22} do not seem to apply.
Inspired by \cite{FHL21,athreya2021wellposedness}, we instead estimate moments of $X^\Delta$ directly, using that $V$ has vanishing mean oscillation: $\|\E^s|V_t-V_s|\|_{L_\infty(\Omega)}\lesssim(t-s)^\gamma$ with $\gamma=1-H-\frac{Hd}p-\frac1q$. Since $\gamma$ can be arbitrarily close to $0$, we lift $V$ to the collection $\mathbb{V}=(\mathbb{V}^w)_{w\in\Theta}$ of its iterated integrals (\cref{subsec:lift}), treat the linearized equation as a rough differential equation, and run a Gr\"onwall argument in $L_m(\Omega)$ based on the John--Nirenberg inequality and the sewing lemma (\cref{lem:Z_Gr\"onwall}). This yields the moment stability estimate \eqref{eq:moment_stab}, which is then used (\cref{sec:flow_and_path_by_path_uniqueness}) to construct the semiflow and prove path-by-path uniqueness through an argument of Shaposhnikov \cite{Shaposhnikov:2016aa} and its refinement in \cite{anzeletti2023pathbypath}.

\subsection{Notation}\label{subsec.notation}
Throughout the paper we fix the dimension $d$, the Hurst parameter $H \in (0, 1)$, the integrability parameters $p, q \in [1, \infty]$, and the final time $T \in (0, \infty)$. 
Although our main results assume that $H < 1/2$, some results (especially those in \cref{subsec:add_fc_fbm}) hold even for $H \in [1/2, 1)$.
We write $\Delta_{[s, t]} \defby \set{(u, v) : s \leq u < v \leq t}$.
For $(s, t) \in \Delta_{[0, T]}$, $\|b\|_{L_q([s,t]; L_p(\R^d))}$ denotes the mixed norm of \cref{sec:main_results} on the time interval $[s,t]$.

Given 
two-parameter maps $(A_{s, t})_{(s, t) \in \Delta_{[0, T]}}$, we write
\begin{align}\label{eq:delta_A_def}
\delta A_{s, u, t} \defby A_{s, t} - A_{s, u} - A_{u, t}, \quad s < u < t.
\end{align}
We denote by $\norm{\cdot}_{\rho\vr; [s, t]}$ the $\rho$-variation norm \cite[Definition~5.1]{friz_victoir_2010}. That is, 
\begin{align*}
  \norm{f}_{\rho\vr; [s, t]} 
  \defby \Big( \sup_{\pi} \sum_{[u, v] \in \pi} \abs{f_v - f_u}^{\rho} \Big)^{\frac{1}{\rho}},
\end{align*}
where the supremum is over all partitions $\pi$ of the interval $[s, t]$. 
We simply write $\norm{\cdot}_{\rho\vr} \defby \norm{\cdot}_{\rho\vr; [0, T]}$. 
The following inequalities hold: for $s < u < t$
\begin{equation}\label{eq:variation_ineq}
  (\norm{f}_{\rho\vr; [s, u]}^{\rho} + \norm{f}_{\rho\vr; [u, t]}^{\rho})^{\frac{1}{\rho}}
  \leq \norm{f}_{\rho\vr; [s, t]} 
  \leq \norm{f}_{\rho\vr; [s, u]} + \norm{f}_{\rho\vr; [u, t]}.
\end{equation}

For a filtered probability space $(\Omega, \mathcal{F}, (\mathcal{F}_t)_{t \in [0, T]}, \P)$, we write 
\begin{align*}
  \norm{G}_{L_m(\Omega)} \defby \Big( \int_{\Omega} \abs{G(\omega)}^m \dd \P(\omega) \Big)^{\frac{1}{m}},
\end{align*}
with usual convention for $m = \infty$. 
We write $\E$ or $\E^{\P}$ for the expectation with respect to $\P$.
When there is no confusion on the filtration $(\mathcal{F}_t)$, we denote by $\E^t$ 
the conditional expectation with respect to $\mathcal{F}_t$. Furthermore, given a $\sigma$-algebra $\mathcal{G} \subseteq \mathcal{F}$, 
we write 
\begin{align*}
  \norm{G}_{L_m(\Omega) \vert \mathcal{G}} \defby \E[ \abs{G}^m \vert \mathcal{G}]^{\frac{1}{m}}.
\end{align*}

\textbf{Convention on constants}. Throughout the paper, $C$ denotes a positive constant whose value may change from line to line, with its dependence always specified in the corresponding statement. However, for the sake of brevity, we will not explicitly state the dependencies of the constants on the parameters $H$, $d$, $p$, and $q$, which are considered to be fixed.

\vspace{1ex}
\textbf{Function spaces}. We write $\mathcal{B}_{p, q}^{\alpha}$ for the Besov space on $\mathbb{R}^d$, whose precise definition can be found in \cite[Section~2.7]{BCD11}. In our paper, we always set $q = \infty$, and simply write $\mathcal{B}_{p}^{\alpha} \defby \mathcal{B}_{p, \infty}^{\alpha}$.
We write $W^{\alpha}_m([0, T])$ for the Sobolev space on $[0, T]$ equipped with norm 
\begin{align*}
  \norm{\psi}_{W^{\alpha}_m([0, T])} \defby \norm{\psi}_{L_m([0, T])} + 
  \Big( \int_0^T \int_0^T \frac{\abs{\psi(s) -\psi(t)}^m}{\abs{s-t}^{1 + \alpha m}} \dd s \dd t \Big)^{\frac{1}{m}}
\end{align*}
when $\alpha \in(0,1)$ and $m \in[1,\infty)$, and 
\begin{align*}
  \norm{\psi}_{W^{1}_m([0, T])} \defby \norm{\psi}_{L_m([0, T])} + \norm{\psi'}_{L_m([0, T])}
\end{align*}
when $m \in[1,\infty)$.

\section{Preliminaries}\label{s:3}

We collect a number of auxiliary results and estimates that will be useful for later considerations.
\subsection{Auxiliary estimates for stochastic processes}
\begin{definition}[{\cite[Definition~1.6]{friz_victoir_2010}}]\label{def:control}
Let $0\le S<T$. We say that a map $w:\Delta_{[S,T]}\to\R_+$ is a \emph{control} if it 
  is continuous, $w(s, s) = 0$ for any $s\in[S,T]$ and
 \begin{equation*}
  w(s, u) + w(u, t) \leq w(s, t) \quad \text{for all } (s,t)\in\Delta_{[S,T]},\,u\in[s,t].
 \end{equation*}
\end{definition}
We will be frequently using the fact that the maps
\begin{equation}\label{controlmaps}
(s,t)\mapsto t-s\,\text{ and }\,(s, t) \mapsto \norm{b}_{L_q([s, t]; L_p(\R^d))}^q\, \text{ are controls},
\end{equation}
where $p\in[1,\infty]$ and $q\in[1,\infty)$. 
Furthermore, given two controls $w_1$ and $w_2$ and two nonnegative numbers $\nu_1, \nu_2$ with $\nu_1 + \nu_2 \geq 1$, 
the map (\cite[Exercise~1.10]{friz_victoir_2010})
\begin{equation}\label{frizcontrol}
  (s, t) \mapsto w_1(s,t)^{\nu_1} w_2(s, t)^{\nu_2}\text{ is a control}.
\end{equation}

Next, we recall the stochastic sewing lemma with conditioning. 
In contrast to \cite[Theorem~3.1]{Le:2023aa}, where moments are taken conditionally on $\F_s$, we need moment bounds conditionally on a $\sigma$-algebra $\mathcal{G}\subset\F_S$ which is fixed throughout; see the proof of \cref{L:firstb}.
\begin{proposition}[Stochastic sewing lemma with conditioning]\label{l:SSL}
Let $0\le S <T$. Fix $m\in[2,\infty)$, $n\in[m,\infty]$. Let   $(A_{s,t})_{(s,t)\in\Delta_{[S,T]}}$ be a family of random variables in $L_m(\Omega,\R^d)$  such that $A_{s,t}$ is $\F_t$-measurable for any $(s,t)\in\Delta_{[S,T]}$. Let $\mathcal{G}\subset\F_S$ be a $\sigma$-algebra. Assume that there exist constants $\Gamma_1,\Gamma_2,\Gamma_3\ge0$, exponents $\alpha_1,\alpha_2,\alpha_3\ge0$ and $\beta_1,\beta_2,\beta_3>0$ with
\begin{equation}\label{SSLexponents}
\alpha_1+\beta_1>\frac12,\quad \alpha_2+\beta_2>\frac12,\quad \alpha_3+\beta_3>1,
\end{equation}
and controls $w_1, w_2, w_3$ such that for any $(s,t)\in\Delta_{[S,T]}$ and $u:=(s+t)/2$ we have
\begin{align}\label{SSLcond1}
&\bigl\|\,\|A_{s,t}\|_{L_m(\Omega)|\mathcal{G}}\bigr\|_{L_n(\Omega)}\le \Gamma_1 w_1(s,t)^{\alpha_1}(t-s)^{\beta_1}+\Gamma_2 w_2(s,t)^{\alpha_2}(t-s)^{\beta_2},\\
&\bigl\|\,\|\E^s [\delta A_{s,u,t}]\|_{L_m(\Omega)|\mathcal{G}}\bigr\|_{L_n(\Omega)}\le \Gamma_3 w_3(s,t)^{\alpha_3}(t-s)^{\beta_3}.
\label{SSLcond2}
\end{align}
Further, suppose that there exists a process $\A=\{\A_t: t\in[S,T]\}$, such that for any $(s,t)\in\Delta_{[S,T]}$ one has
\begin{equation}\label{SSLcond3}
\A_t-\A_s=\lim_{N\to\infty}\sum_{i=0}^{N-1} A_{s+i\frac{t-s}N,s+(i+1)\frac{t-s}N}\,\,\text{in probability}.
\end{equation}

Then there exists a constant $C=C(\alpha_1,\alpha_2,\alpha_3,\beta_1,\beta_2,\beta_3,m)$ independent of $S,T,n, \Gamma_i, w_i$ such that for every $(s,t)\in\Delta_{[S,T]}$ we have
\begin{equation}\label{SSLres}
\bigl\|\,\|\A_t-\A_s\|_{L_m(\Omega)|\mathcal{G}}\bigr\|_{L_n(\Omega)}\le C\sum_{j=1}^3\Gamma_j w_j(s,t)^{\alpha_j}(t-s)^{\beta_j}.
\end{equation}
\end{proposition}
\begin{remark}\label{rem:SSL}
  Since $n\ge m$ and $\mathcal{G}\subset\F_s$, the conditional Jensen inequality yields
  \begin{equation*}
    \|\,\|X\|_{L_m(\Omega)|\mathcal{G}}\|_{L_n(\Omega)}\le\|X\|_{L_n(\Omega)}.
  \end{equation*}
  Hence \eqref{SSLcond2} holds whenever $\|\E^s\delta A_{s,u,t}\|_{L_n(\Omega)}\le \Gamma_3 w_3(s,t)^{\alpha_3}(t-s)^{\beta_3}$; in particular, it holds with $\Gamma_3=0$ if $\E^s\delta A_{s,u,t}=0$, which is the only case used in this article.
\end{remark}
\begin{proof}[Proof of \cref{l:SSL}]
  The argument follows the dyadic scheme of \cite[Theorem~2.1]{khoa20} and \cite[Theorem~3.1]{Le:2023aa}; we only indicate the necessary modifications.
  For a random variable $X$ we write $\vertiii{X}\defby\|\,\|X\|_{L_m(\Omega)|\mathcal{G}}\|_{L_n(\Omega)}$; by the conditional Minkowski inequality, $\vertiii{\cdot}$ satisfies the triangle inequality.
  For $j=1,2,3$ we put $\theta_j\defby\alpha_j+\beta_j$ and $\hat w_j(s,t)\defby w_j(s,t)^{\alpha_j/\theta_j}(t-s)^{\beta_j/\theta_j}$, which is a control by \eqref{frizcontrol}, so that $w_j(s,t)^{\alpha_j}(t-s)^{\beta_j}=\hat w_j(s,t)^{\theta_j}$.

  Fix $(s,t)\in\Delta_{[S,T]}$. For $k\ge0$ and $i=0,\ldots,2^k$ let $t^k_i\defby s+i2^{-k}(t-s)$, $I^k_i\defby[t^k_i,t^k_{i+1}]$ and $A^k\defby\sum_{i=0}^{2^k-1}A_{t^k_i,t^k_{i+1}}$. Then
  \begin{equation*}
    A^k-A^{k+1}=\sum_{i=0}^{2^k-1}\E^{t^k_i}\delta_i+\sum_{i=0}^{2^k-1}(\delta_i-\E^{t^k_i}\delta_i)=:I^k+J^k,\qquad \delta_i\defby\delta A_{t^k_i,t^{k+1}_{2i+1},t^k_{i+1}}.
  \end{equation*}
  By \eqref{SSLcond2}, $\vertiii{I^k}\le\Gamma_3\sum_i\hat w_3(I^k_i)^{\theta_3}$.
  The summands of $J^k$ are martingale differences with respect to $(\F_{t^k_{i+1}})_i$, and $\mathcal{G}\subset\F_{t^k_0}$.
  The conditional Burkholder--Davis--Gundy inequality \cite[Lemma~2.4]{Le:2023aa}, followed by the Minkowski inequality in $L_{m/2}$ and in $L_{n/2}$ (recall $n\ge m\ge2$), yields
  \begin{equation*}
    \vertiii{J^k}\le C_m\Big(\sum_i\vertiii{\delta_i-\E^{t^k_i}\delta_i}^2\Big)^{\frac12}\le 2C_m\Big(\sum_i\vertiii{\delta_i}^2\Big)^{\frac12},
  \end{equation*}
  where we used that $\|\E^{u}X\|_{L_m(\Omega)|\mathcal{G}}\le\|X\|_{L_m(\Omega)|\mathcal{G}}$ for $u\ge S$. Since controls are monotone with respect to inclusion of intervals, \eqref{SSLcond1} yields $\vertiii{\delta_i}\le 3\sum_{j=1,2}\Gamma_j\hat w_j(I^k_i)^{\theta_j}$.
  Furthermore, $\hat w_j(I^k_i)\le w_j(s,t)^{\alpha_j/\theta_j}(2^{-k}(t-s))^{\beta_j/\theta_j}=2^{-k\beta_j/\theta_j}\hat w_j(s,t)$, and hence by superadditivity of $\hat w_j$, for any $\kappa\ge1$,
  \begin{equation*}
    \sum_i\hat w_j(I^k_i)^{\kappa}\le\max_i\hat w_j(I^k_i)^{\kappa-1}\sum_i\hat w_j(I^k_i)\le 2^{-k\frac{\beta_j}{\theta_j}(\kappa-1)}\hat w_j(s,t)^{\kappa}.
  \end{equation*}
  Applying this with $\kappa=2\theta_j>1$ for $j=1,2$ and with $\kappa=\theta_3>1$ for $j=3$, we obtain
    \begin{gather*}
    \vertiii{A^k-A^{k+1}}\le C2^{-k\varkappa}\sum_{j=1}^3\Gamma_j\hat w_j(s,t)^{\theta_j},\quad\text{where}\\
    \varkappa\defby\min\Big\{\frac{\beta_1(2\theta_1-1)}{2\theta_1},\frac{\beta_2(2\theta_2-1)}{2\theta_2},\frac{\beta_3(\theta_3-1)}{\theta_3}\Big\}>0.
  \end{gather*}
  By \eqref{SSLcond3} with $N=2^k$, $A^k\to\A_t-\A_s$ in probability, hence almost surely along a subsequence. The conditional Fatou lemma and the Fatou lemma give
  \begin{equation*}
    \vertiii{\A_t-\A_s}\le\vertiii{A_{s,t}}+\sum_{k=0}^\infty\vertiii{A^k-A^{k+1}}\le C\sum_{j=1}^3\Gamma_j\hat w_j(s,t)^{\theta_j},
  \end{equation*}
  where we used \eqref{SSLcond1} once more for $A^0=A_{s,t}$. This is \eqref{SSLres}.
\end{proof}

The following result is called the quantitative John--Nirenberg inequality.
\begin{proposition}[{\cite[Theorem~2.3 and Corollary~3.5]{le22}}]\label{prop:vmo_john_nirenberg}
Let $0\le S\le T$ and let $(\A_t)_{t\in[S,T]}$ be a continuous process taking values in $\R^d$ adapted to the filtration $(\F_t)_{t\in[S,T]}$. 
\begin{enumerate}[\rm{(}i\rm{)}]
\item \label{jnone} 

Suppose that there exist a constant $\Gamma>0$ such that for any $(s,t)\in\Delta_{[S,T]}$
\begin{equation}\label{eq:vmo_gamma}
\E^s  |\A_t-\A_s|     \le \Gamma \quad a.s.
\end{equation}
Then for any $m\ge1$ there exists a constant $C=C(m)$ such that \begin{equation}\label{jnres}
\bigl\|\,\|\A_t-\A_s\|_{L_m(\Omega)|\F_s}\bigr\|_{L_\infty(\Omega)}\le C \Gamma,\qquad (s,t)\in\Delta_{[S,T]}.
\end{equation}

\item\label{jntwo}
Assume that there exists $r\in(0,1]$ and a control $w$ such that a more restrictive condition
\begin{equation}\label{eq:vmo_control}
\E^s  |\A_t-\A_s|     \le w(s,t)^r \quad a.s.
\end{equation}
holds for all $(s,t)\in\Delta_{[S,T]}$. Then there exist positive constants $c = c(r),C=C(r)$ such that  
  \begin{equation}\label{eq:JN_exponential}
\E [e^{c w(S,T)^{-\frac{r}{1-r}} \sup_{t\in[S,T]} |\A_t - \A_S|^{\frac1{1-r}}}] \le C \quad\text{ if }r<1,
  \end{equation}
  and 
  \begin{align}
    \sup_{t\in[S,T]} |\A_t - \A_S|\le C w(S,T)\quad\text{a.s. if}\quad r=1.\label{eq:JN_exponential2}
  \end{align}
  \end{enumerate}
\end{proposition}
\begin{remark}\label{rem:vmo}
   A continuous adapted process $\mathcal{A}$ satisfying \eqref{eq:vmo_gamma} is called a process of \emph{bounded mean oscillation} (BMO),
  and a process satisfying \eqref{eq:vmo_control} with $r\in(0,1)$ is an instance of processes of \emph{vanishing mean oscillation} (VMO).
  \Cref{prop:vmo_john_nirenberg} shows that BMO processes have conditional moments of all orders and that VMO processes satisfying \eqref{eq:vmo_control} are exponentially integrable, with the rate of exponential integrability identified by the regularity index $r$. This is how we obtain moment and exponential bounds from first-moment bounds on conditional increments throughout the article.
\end{remark}
%

The following elementary result is useful to remove singularities near terminal times. Although probably well known, we learned about these techniques from \cite{bellingeri21,leling21}.
  
\begin{lemma}[Taming singularities]\label{lem:singularity}
Let $(E,d)$ be a metric space, $0\le S\le T$, and let $Y$ be a  function $[S, T] \to E$.
\begin{enumerate}[\rm{(}i\rm{)}]
\item \label{lsone} Suppose that $S>0$ and  there exist constants $\Gamma, \tau,\eta,\lambda \ge0$  such that 
\begin{equation}\label{condtame1}
d(Y_s,Y_t)\le  \Gamma s^{-\eta} (t-s)^{\tau}(1+s^{-\lambda}t^\lambda), \quad S\le s\le t\le T. 
\end{equation}
Then 
\begin{equation}\label{restame1}
d(Y_s,Y_t)\le  C \Gamma s^{-\eta} (t-s)^{\tau}, \quad S \le  s \le  t \leq T,
\end{equation}
where $C=(1+2^\lambda)(2^\tau-1)^{-1}$.

\item \label{lstwo} Suppose that $Y$ is continuous at $S$ and
that there exists constants $\Gamma,\tau,\eta \ge0$,  
$\tau>\eta$ such that
\begin{equation}\label{condtame}
d(Y_s,Y_t)\le  \Gamma s^{-\eta} (t-s)^{\tau}, \quad S<s\le t\le T. 
\end{equation}
Then 
\begin{equation}\label{restame}
d(Y_s,Y_t)\le  C  \Gamma  (t^{\gamma}-s^{\gamma})^{\tau}, \quad  S \leq s < t \leq T,
\end{equation}
where $\gamma:=1-\frac\eta\tau\in(0,1]$ and
$C=(2^\tau-1)^{-1}2^{\tau\eta/(\tau-\eta)}\tau^\tau (\tau-\eta)^{-\tau}$.
\end{enumerate}
  \end{lemma}

\begin{proof}
(i)   We use an idea from {\cite[Lemma 3.4]{leling21}}. For  $S<s\le t\le T$ we set
\begin{equation}\label{tndef}
t_n:= s + (t - s) 2^{- n},\quad n\in\Z_+.
\end{equation}
By the triangle inequality, 
\begin{equation*}
d(Y_s, Y_t) \leq \sum_{n=1}^{\infty} d(Y_{t_n}, Y_{t_{n-1}}),
\end{equation*}
where we also used that \eqref{condtame1} implies that $Y$ is continuous at $s$.
By our choice of $t_n$, we have  $t_n\ge s$ and $t_n\ge t_{n-1}/2$. Since $t_{n-1} - t_{n}= 2^{- n} ( t-s )$,  \eqref{condtame1} yields
\begin{equation*}
d(Y_s, Y_t) \leq \Gamma \sum_{n=1}^{\infty} s^{-\eta}(t-s)^{\tau}2^{-n\tau}(1+2^\lambda)=
\Gamma(1+2^\lambda)(2^\tau-1)^{-1}s^{-\eta}(t-s)^{\tau}.
\end{equation*}

(ii) Following \cite[Proposition~2.4]{bellingeri21}, we introduce the time change  $\phi (t):=t^{\mu}$, where $\mu:=\frac{\tau}{\tau-\eta}$. Note that $\mu\in(1,\infty)$ by our assumptions.

Note that if $u,v>0$ are such that $v/2\le u \le v$, then 
\begin{equation*}
v^{\mu} - u^{\mu} \le \mu v^{\mu - 1}  (v - u)\le 2^{\mu-1} \mu  u^{\mu - 1} (v-u).
\end{equation*}	
Therefore,
\begin{equation}\label{step0tsl}
d(Y_{\phi (u)},Y_{\phi (v)})\le  C_1  \Gamma u^{\mu(\tau-\eta)-\tau} (v-u)^{\tau},\quad S<u^\mu\le v^\mu\le T,\,\,  u\ge v/2,
\end{equation}
where $C_1=2^{\tau(\mu-1)} \mu^\tau$. The choice of $\mu$ ensures that $\mu(\tau-\eta)-\tau=0$. 

To remove the restriction $u\ge v/2$ we argue as in the proof of part (i). For arbitrary $(s,t)\in\Delta_{[S^{1/\mu},T^{1/\mu}]}$, we define $t_n$ as in \eqref{tndef}.
Since $t_n\ge t_{n-1}/2$, bound \eqref{step0tsl} is applicable for $u=t_n$, $v=t_{n-1}$.  
 Then, by the triangle inequality, continuity of $Y$, \eqref{step0tsl} and \eqref{condtame} we have
  for any $(s,t)\in\Delta_{[S^{1/\mu},T^{1/\mu}]}$
\begin{align*}
d(Y_{\phi(s)}, Y_{\phi(t)}) &\le \sum_{n=1}^{\infty} d(Y_{\phi(t_n)}, Y_{\phi(t_{n-1})})\le  C_1\Gamma  
\sum_{n=1}^{\infty} (t_{n-1}-t_n)^{\tau}= C_1\Gamma (t-s)^{\tau}\sum_{n=1}^{\infty} 2^{-n\tau}
\\
&= C_1 \Gamma (2^\tau-1)^{-1} (t-s)^\tau.
\end{align*}
Inequality \eqref{restame} follows now  by change of variables $(s, t) \mapsto (s^{1/\mu}, t^{1/\mu})$ and denoting $1/\mu=1-\eta/\tau$ by $\gamma$.
  \end{proof}
  
Next, we give an extension of the classical Kolmogorov continuity theorem to the case when the increments of a process are estimated by a general control. In this case, the result provides estimates on the variation norm of the process instead of the H\"older norm in the classical case.

\begin{lemma}\label{lem:variation-kolmogorov}
Let $(V_t)_{t \in [0, T]}$ be a continuous stochastic process. Assume that for some $m\in(1,\infty)$, a control $w$, $\mu \in (\frac{1}{m}, 1]$, a $\sigma$-algebra $\G\subset\F$ we have
\begin{equation}\label{condkol}
\|\,\|V_t-V_s\|_{L_m(\Omega)|\G}\|_{L_\infty(\Omega)}\le w(s,t)^\mu,\quad (s,t)\in\Delta_{[0,T]}.
\end{equation}
Let $\mu'\in  (\frac{1}{m}, \mu)$. Then there exists a \textbf{deterministic} constant $C$, which depends only on $\mu,\mu',m$,  such that for any  $0 \le s < t \le T$
\begin{equation}\label{boundkolmi}
\bigl\|\,\| V\|_{(1/\mu')\var;[s,t]}\bigr\|_{L_m(\Omega)|\G}
\le C  w(s,t)^{\mu }.
\end{equation}
\end{lemma}
\begin{proof}
Fix $\eps>0$, $S\in[0,T]$ and consider the modified control  $\wt w(s,t):=w(s,t)+\eps(t-s)$, $S\le s \le t\le T$.
Since the function  $s \mapsto \wt w(S, s)$ is strictly increasing, we can introduce the following change of time
\begin{equation*}
\tau (t):= (\wt w(S, \cdot))^{-1}(t)
\end{equation*}
and put $\wt V_t:=V_{\tau(t)}$. Then for any $S \le s < t \leq  \tau^{-1}(T)=\wt w(S,T)$, by assumptions of the lemma, we have 
\begin{equation*}
	\|\,\|\wt V_t-\wt V_s\|_{L_m(\Omega)|\G}\|_{L_\infty(\Omega)}\le \wt w(\tau(s),\tau(t))^\mu\le 
	\bigl(\wt w(S,\tau(t))-\wt w(S,\tau(s)\bigr)^\mu=(t-s)^\mu,
\end{equation*}
where the second inequality follows from the fact that $\wt w$ is a control.

Fix now $\mu'\in  (\frac{1}{m}, \mu)$ We apply now the Besov-variation embedding \cite[Corollary A.3]{friz_victoir_2010}. We see that there exists a universal constant $C=C(\mu',m)$ such that
for any $S \le s < t \leq  \tau^{-1}(T)=\wt w(S,T)$
\begin{equation*}
\|\wt V\|_{(1/\mu')\var;[s,t]} \le C (t-s)^{\mu' - \frac{1}{m}}
	\left(\int_s^t
	\int_s^t  \frac{| \wt V_u-\wt V_v |^m}{|u - v|^{1 + m \mu' }}\,dudv\right)^{1/m}.
\end{equation*}
Combining the previous two inequalities, we have
\begin{equation*}
\E\bigl[	\|\wt V\|^m_{(1/\mu')\var;[s,t]}\bigr|\G\bigr] \le C (t-s)^{m\mu' -1}
\int_s^t
	\int_s^t  \frac{\E [| \wt V_u-\wt V_v |^m|\G]}{|u - v|^{1 + m \mu' }}\,dudv\le C (t-s)^{\mu m},
\end{equation*}
for some constant $C=C(\mu,\mu',m)$. Here we used that $1+m\mu'-m\mu<1$ which ensures that the integral is finite.
Thus, for any $S \le s < t \le T$ we finally deduce that
\begin{align*}
	\bigl\|\,\| V\|_{(1/\mu')\var;[s,t]}\bigr\|_{L_m(\Omega)|\G}
	&=\bigl\|\,\| \wt V\|_{(1/\mu')\var;[\tau^{-1}(s),\tau^{-1}(t)]}\bigr\|_{L_m(\Omega)|\G}\\
	& \le C (\tau^{-1}(t)-\tau^{-1}(s))^{\mu }\\
	&= C (\wt w(S,t)-\wt w(S,s))^{\mu }\\
	&\le C ( w(S,t)- w(S,s))^{\mu }+C\eps^\mu (t-s)^\mu,
\end{align*}
where $C=C(\mu,\mu',m)$. Since $C$ does not depend on $\eps$, by letting $\eps\to0$ we get
\begin{align*}
  \bigl\|\,\| V\|_{(1/\mu')\var;[s,t]}\bigr\|_{L_m(\Omega)|\G}
\le C ( w(S,t)- w(S,s))^{\mu }
\end{align*}
for every $S\le s\le t\le T$. Since $C$ is independent of $S$, choosing  $S=s$ yields the desired bound \eqref{boundkolmi}.
  \end{proof}

\subsection{Girsanov theorem and Sobolev spaces}\label{subsec:girsanov_prelim}
Let $W$ be a standard $d$-dimensional $(\F_t)$-Brownian motion. Besides the fractional Brownian motion $B^H=\int_0^\cdot K_H(\cdot,r)\,dW_r$ from \eqref{eq:fbm_kernel}, we shall use the \emph{Riemann--Liouville fractional Brownian motion}
\begin{equation}\label{eq:RL_def}
  \tilde B^H_t\defby \tilde c_H\int_0^t (t-r)^{H-\frac12}\,dW_r,\quad t\in[0,T],
\end{equation}
where $H\in(0,\frac12)$ and $\tilde c_H\defby\big(\frac1{2H}+\int_0^\infty((1+s)^{H-\frac12}-s^{H-\frac12})^2\,ds\big)^{-1/2}$.

The constant $\tilde c_H$ is the normalizing constant of the Mandelbrot--Van Ness representation \cite{Mandelbrot:1968aa}, which we use in \cref{sec:stability}: if $(\tilde W_t)_{t\in\R}$ is a two-sided standard $d$-dimensional Brownian motion and $\mathcal{K}(t,r)\defby\tilde c_H[(t-r)_+^{H-\frac12}-(-r)_+^{H-\frac12}]$, then
\begin{equation}\label{eq:mandelbrot_van_ness}
  B^H_t=\int_{-\infty}^t\mathcal{K}(t,r)\dd\tilde W_r,\quad t\ge0,
\end{equation}
is a fractional Brownian motion; the choice of $\tilde c_H$ ensures $\E|B^H_t-B^H_s|^2=|t-s|^{2H}$.
Let $(\tilde{\mathcal F}_t)_{t\in\R}$ be the filtration generated by $\tilde W$. Then $\F^{B^H}_t\subseteq\tilde{\mathcal F}_t$, and for every $v\ge0$ we have the decomposition
\begin{equation}\label{eq:MVN_decomp}
\begin{gathered}
  B^H_{v+t}-B^H_v=\zeta^v_t+\bar B^v_t,\quad t\ge0,\quad\text{where}\\
  \zeta^v_t\defby\int_{-\infty}^v(\mathcal{K}(v+t,r)-\mathcal{K}(v,r))\dd\tilde W_r,\quad \bar B^v_t\defby\tilde c_H\int_v^{v+t}(v+t-r)^{H-\frac12}\dd\tilde W_r,
\end{gathered}
\end{equation}
$\zeta^v=\E[B^H_{v+\cdot}-B^H_v\vert\tilde{\mathcal F}_v]$ is an $\tilde{\mathcal F}_v$-measurable continuous path with $\zeta^v_0=0$, and $\bar B^v$ is independent of $\tilde{\mathcal F}_v$ and is a Riemann--Liouville fractional Brownian motion \eqref{eq:RL_def}, driven by the Brownian motion $(\tilde W_{v+t}-\tilde W_v)_{t\ge0}$. In particular, $\bar B^0$ is the Riemann--Liouville fractional Brownian motion driven by $(\tilde W_t)_{t\ge0}$.
Conversely, given a fractional Brownian motion $B^H$ on $(\Omega,\F,\P)$, by the transfer theorem \cite[Theorem~8.17]{kallenberg} we can extend the probability space so that it carries a two-sided Brownian motion $\tilde W$ for which \eqref{eq:mandelbrot_van_ness} holds almost surely. Such an extension does not change the joint law of $B^H$ and any process which is measurable with respect to $B^H$, such as a strong solution of \eqref{eq:main-sde}, nor the conditional laws of such processes given $\F^{B^H}_t$. Hence, when proving statements that concern only such laws, we may and do assume that $\tilde W$ exists on the original probability space.

In the rest of this subsection and in \cref{subsec:girsanov,subsec:linearization}, $G$ denotes either $B^H$ or $\tilde B^H$, and we write $K^G\defby K_H$ if $G=B^H$ and $K^G(t,r)\defby\tilde c_H(t-r)^{H-\frac12}$ if $G=\tilde B^H$, so that $G_t=\int_0^tK^G(t,r)\,dW_r$.
We denote by $K^G$ also the integral operator $(K^G h)(t)\defby\int_0^t K^G(t,r)h(r)\,dr$, $t\in[0,T]$, acting on $h\in L_2([0,T];\R^d)$.
This operator is injective: for $G=B^H$ see \cite[Section~5.1.3]{nualart06}, and for $G=\tilde B^H$ it is a constant multiple of the Riemann--Liouville fractional integral of order $H+\frac12$, see \cite{Gorenflo:2015vr}.
Let $\mathcal{H}^G\defby K^G(L_2([0,T];\R^d))$ be its range, that is, the space of functions $\phi\colon[0,T]\to\R^d$ of the form $\phi=K^Gh$ for some $h\in L_2([0, T];\R^d)$, which is unique and denoted by $(K^G)^{-1}\phi$. Every $\phi\in\mathcal{H}^G$ is continuous with $\phi_0=0$. We equip $\mathcal{H}^G$ with the norm
\begin{equation}\label{normH}
\|\phi\|_{\mathcal{H}^G}^2:= \int_0^T | (K^G)^{-1} \phi(t)|^2 \,d t.
\end{equation}
The space $\mathcal{H}^G$ is identified with the Cameron--Martin space of $G$.
We simply write $\mathcal{H}\defby\mathcal{H}^{B^H}$ for the Cameron--Martin space of the fractional Brownian motion.
In both cases, $(K^G)^{-1}$ is causal, that is, $((K^G)^{-1}\phi)(t)$ depends only on $(\phi_r)_{r\in[0,t]}$; in particular, $(K^G)^{-1}\phi$ is adapted if $\phi$ is adapted.

For an adapted process $\phi$ with paths in $\mathcal{H}^G$ and $t\in[0,T]$ we set
\begin{equation}\label{eq:xi_def}
\xi_{t}(\phi)=\xi^G_t(\phi):=  \exp \Big\{  - \int_0^{t}  (K^G)^{-1} \phi(s) \,dW_s
  -\frac12\int_0^{t} | (K^G)^{-1} \phi(s)|^2 \,ds \Big\}.
\end{equation}

\begin{proposition}[Girsanov via Novikov]\label{lem:girsanov}
Let $H\in(0,1)$,  $T>0$. Let $G$ be either $B^H$ or $\tilde B^H$ (in the latter case $H<\frac12$). Let $\phi\colon[0,T]\times\Omega\to\R^d$ be an adapted process with paths in $\mathcal{H}^G$. 
For a given $\lambda\in\R$, consider the condition
  \begin{equation}\label{eq:novikov}
    \E_{\P} e^{\lambda \|\phi\|_{\mathcal{H}^G}^2}    < \infty.
  \end{equation}
 Define  the measure $\wt\P$ by $d \wt\P:= \xi_{T}(\phi) d\P$.
Then the following statements hold.
\begin{enumerate}[\rm{(}i\rm{)}]
\item Suppose that \eqref{eq:novikov} holds for $\lambda=1/2$. Then $\wt \P$ is a probability measure, the process $\wt W\defby W+\int_0^\cdot (K^G)^{-1}\phi(r)\,dr$ is an $(\F_t)$-Brownian motion under $\wt\P$, and $G_t+\phi_t=\int_0^tK^G(t,r)\,d\wt W_r$, $t\in[0,T]$. In particular, the law of $G+\phi$ under $\wt\P$ coincides with the law of $G$ under $\P$, and if $G=B^H$, then $B^H+\phi$ is an $(\F_t)$-fractional Brownian motion under $\wt\P$.

\item Suppose that \eqref{eq:novikov} holds for every  $\lambda \in \mathbb{R}$. Then, for any $\lambda \in \mathbb{R}$, we have   
  \begin{align}\label{eq:gir_density_mom}
    &\E_{\P}[\xi_{T}(\phi)^{\lambda}]
     \leq 
     \bigl(\E_{\P} e^{(2\lambda^2-\lambda)\|\phi\|_{\mathcal{H}^G}^2}\bigr)^{\frac{1}{2}} \tas
     \\\shortintertext{and}
     &\E_{\wt \P}[\xi_{T}(\phi)^{\lambda}]
     \leq 
     \bigl(\E_{\wt\P} e^{(2\lambda^2+\lambda)\|\phi\|_{\mathcal{H}^G}^2}\bigr)^{\frac{1}{2}} \tas.
     \label{eq:gir_density_mom2}
  \end{align}
\end{enumerate}
\end{proposition}
\begin{proof}
(i) For $G=B^H$, this is \cite[Theorem~2]{Nualart:2002uv}; condition \eqref{eq:novikov} ensures that condition (ii) of \cite[Theorem~2]{Nualart:2002uv} holds. In general, $h\defby(K^G)^{-1}\phi$ is adapted by causality of $(K^G)^{-1}$, and \eqref{eq:novikov} with $\lambda=\frac12$ is Novikov's condition for $h$, since $\|\phi\|_{\H^G}=\|h\|_{L_2([0,T])}$. Hence $(\xi_t(\phi))_{t\in[0,T]}$ is a $\P$-martingale, and by the classical Girsanov theorem $\wt W$ is an $(\F_t)$-Brownian motion under $\wt\P$. Finally, $G_t+\phi_t=\int_0^tK^G(t,r)\,dW_r+\int_0^tK^G(t,r)h(r)\,dr=\int_0^tK^G(t,r)\,d\wt W_r$.

(ii) Fix $\lambda\in\R$. We use the Cauchy--Schwarz inequality to get that 
\begin{equation*}
\E_{\P} [\xi_{T}(\phi)^{\lambda}]=\E_{\P} [\xi_{T}(2\lambda\phi)^{\frac12} e^{(\lambda^2-\frac\lambda2)\|\phi\|^2_{\H}}]
    \le \bigl(\E_{\P} \xi_{T}(2 \lambda \phi)\bigr)^{\frac{1}{2}}\bigl(\E_{\P} e^{(2\lambda^2-\lambda)\|\phi\|_{\H^G}^2}\bigr)^{\frac{1}{2}}.
\end{equation*}
Note that the process $t \mapsto \xi_t(2 \lambda \phi)$ is a martingale under $\P$ because  the Novikov's condition
\begin{align*}
   \E_{\P} e^{\frac12 \|2 \lambda\phi\|_{\H^G}^2}    < \infty
\end{align*}
is satisfied under our assumptions.
Hence, $\E_{\P} \xi_{T}(2 \lambda \phi)=1$ and we obtain \eqref{eq:gir_density_mom}. 

We prove \cref{eq:gir_density_mom2}  similarly. We note that the process $\wt W:=W+ \int_0^{\cdot}  (K^G)^{-1} \phi(s)\,ds$ is a standard Brownian motion under $\wt \P$. Hence, we can rewrite $\xi$ as 
\begin{equation*}
	\xi_{t}(\phi)=  \exp \Big\{  - \int_0^{t}  (K^G)^{-1} \phi(s) \,d\wt W_s 
	+\frac12\int_0^{t} | (K^G)^{-1} \phi(s)|^2 \,ds \Big\},\quad t\in[0,T].
\end{equation*}
It follows from \eqref{eq:novikov} and \eqref{eq:gir_density_mom} that for any $\lambda>0$
\begin{equation*}
\E_{\wt \P} e^{\lambda \|\phi\|_{\H^G}^2} =\E_{ \P} [e^{\lambda \|\phi\|_{\H^G}^2} \xi_{T}(\phi)]\le  \bigl(\E_{ \P} e^{2\lambda \|\phi\|_{\H^G}^2}\bigr)^{\frac12}\E_{\P} \bigl(\xi_{T}(\phi)^{2}\bigr)^{\frac12} < \infty.
\end{equation*}
Put $h\defby(K^G)^{-1}\phi$ and, for $\lambda\in\R$, let $M^\lambda_t\defby\exp\{-2\lambda\int_0^th(s)\,d\wt W_s-2\lambda^2\int_0^t|h(s)|^2\,ds\}$ be the Dol\'eans exponential of $-2\lambda\int h\,d\wt W$. Since $\E_{\wt\P}\exp\{\frac12\int_0^T|2\lambda h(s)|^2\,ds\}=\E_{\wt\P}e^{2\lambda^2\|\phi\|^2_{\H^G}}<\infty$, Novikov's criterion under $\wt\P$ shows that $M^\lambda$ is a $\wt\P$-martingale, so that $\E_{\wt\P}M^\lambda_T=1$. Moreover, $\xi_T(\phi)^\lambda=(M^\lambda_T)^{\frac12}e^{(\lambda^2+\frac\lambda2)\|\phi\|^2_{\H^G}}$. Thus, by the Cauchy--Schwarz inequality,
\begin{align*}
	\E_{\wt \P} [\xi_{T}(\phi)^{\lambda}]
	\le \bigl(\E_{\wt\P} M^\lambda_T\bigr)^{\frac{1}{2}}\bigl(\E_{\wt \P} e^{(2\lambda^2+\lambda)\|\phi\|_{\H^G}^2}\bigr)^{\frac{1}{2}}
	=	\bigl(\E_{\wt \P} e^{(2\lambda^2+\lambda)\|\phi\|_{\H^G}^2}\bigr)^{\frac{1}{2}},
\end{align*}
which is  \cref{eq:gir_density_mom2}.
\end{proof}

For $\alpha\in(1/m,1)$, we define the  Sobolev space $\W^{\alpha}_m([0, T])$ of continuous functions $\psi:[0,T]\to\Rd$ such that $\psi_0=0$ and
\begin{align*}
  \norm{\psi}_{\W^{\alpha}_m([0, T])} \defby  
  \Big( \int_0^T \int_0^T \frac{\abs{\psi(s) -\psi(t)}^m}{\abs{s-t}^{1 + \alpha m}} \dd s \dd t \Big)^{\frac{1}{m}}<\infty.
\end{align*}
We also define the Sobolev space $\W^{1}_m([0, T])$ of continuous functions $\psi:[0,T]\to\Rd$ such that $\psi_0=0$ and
\begin{align*}
  \norm{\psi}_{\W^{1}_m([0, T])} \defby  \norm{\psi'}_{L_m([0, T])}<\infty,
\end{align*}
where $\psi'$ is the weak derivative of $\psi$.

By the Garsia--Rodemich--Rumsey inequality, we have for $\psi\in\W^\alpha_m([0,T])$, 
\begin{align}
  |\psi_t- \psi_s| \leq C (t-s)^{\alpha-\frac{1}{m} } w_{\psi}(s,t)^{\frac{1}{m}}  \quad \forall [s,t]\subset[0,T], \label{eq:GRR}
\end{align}
where $w_{\psi}(s,t)=\int_{[s,t]^2}\frac{|\psi_u-\psi_v|^m}{|u-v|^{1+\alpha m}}dudv$ is a control. When $\alpha=1$, \eqref{eq:GRR} also holds with the control $w_\psi(s,t)=\int_s^t |\psi'_r|^mdr$, which can be verified directly using the relation $\psi_t-\psi_s=\int_s^t \psi'_rdr$ and H\"older inequality. 
It follows that $|\psi_t|\les w_\psi(0,T)^{1/m}$ and hence $\norm{\psi}_{L_m([0,T])}\les w(0,T)^{1/m}\les \norm{\psi}_{\W^\alpha_m([0,T])} $.
From here, we deduce that there is a constant $C$ such that 
\begin{align}
  \norm{\psi}_{\W^{\alpha}_m([0, T])} \leq \norm{\psi}_{W^{\alpha}_m([0, T])} \le C \norm{\psi}_{\W^{\alpha}_m([0, T])} \quad \forall \psi\in\W^\alpha_m([0,T]). \label{eq:equivWW}
\end{align}
Consequently,  classical embeddings on Sobolev spaces $W^\alpha_m$ extends to $\W^\alpha_m$. 
In particular, the embedding $W^\alpha_m\hookrightarrow W^{\alpha-\varepsilon}_m$ (\cite[Proposition~C.1(i)]{galeatigerencser22}) implies $\W^\alpha_m\hookrightarrow \W^{\alpha-\varepsilon}_m$ and there exists a constant $C$ such that
\begin{align}\label{eq:trivialW}
  \norm{\psi}_{\W^{\alpha-\varepsilon}_m([0, T])} \leq C \norm{\psi}_{\W^{\alpha}_m([0, T])} \quad \forall \psi\in\W^\alpha_m([0,T])
\end{align}
whenever $\alpha\in(0,1]$ such that $\alpha-\varepsilon>1/m$.
Another useful embedding is the following Gagliardo--Nirenberg inequality.
\begin{lemma}[Gagliardo--Nirenberg inequality]\label{lem:GG}
    Let $H\in(0,1/2)$,  $T>0$, $1\le p_2\le 2 \le p_1< \infty$, and $\gamma\in(0,1)$ with $1-\gamma>\frac1{p_1}$.
  Take $\theta\in[0,1]$ such that
  $$
  \frac\theta{p_1}+\frac{1-\theta}{p_2}=\frac12
  $$
  and assume that 
  \begin{equation}\label{GNcond}
  H+\gamma \theta<\frac12.
  \end{equation}
  Then there exist positive  constants $\eps=\eps(\theta,\gamma,H)$ and $C=C( p_1,p_2,\gamma,T,\eps)$  such that   any continuous function $\phi\colon[0,T]\to\R^d$ with $\phi_0=0$ satisfies 
  \begin{equation}\label{mainresGG}	
    \|\phi\|_{\W^{H+\frac12+\eps}_2([0,T])}\le C \|\phi\|_{\W^{1-\gamma}_{p_1}([0,T])}+C\|\phi\|_{\W^{1}_{p_2}([0,T])}.
  \end{equation}		
  \end{lemma}
  \begin{proof}
  Fix an $\eps>0$ small enough such that $\epsilon < \gamma$ and 
  \begin{equation}\label{genGNcond}
  H+\gamma \theta<\frac12-\eps(2-\theta).
  \end{equation}	
  This is possible thanks to condition \eqref{GNcond}. 
  By the generalized Gagliardo--Nirenberg inequality \cite[Theorem~1(A)]{BrezisM}, there exists a constant $C=C(p_1,p_2,\gamma,T,\eps)$ such that 
  \begin{align}\label{gGNineq}	
    \|\phi\|_{W^{1-\gamma \theta-\eps(1-\theta)}_2([0,T])}&\le C \|\phi\|_{W^{1-\gamma}_{p_1}([0,T])}^\theta\|\phi\|_{W^{1-\eps}_{p_2}([0,T])}^{1-\theta}\nn\\
    &\le C \|\phi\|_{W^{1-\gamma}_{p_1}([0,T])}+C\|\phi\|_{W^{1-\eps}_{p_2}([0,T])},
  \end{align}		
  where we also used the elementary inequality $a^\theta b^{1-\theta}\le a+b$ valid for all non-negative $a,b$. 
    Here \cite[Theorem~1(A)]{BrezisM} is applicable since $1-\eps$ is not an integer, so that the exceptional case \cite[(1.9)]{BrezisM} does not occur (this matters for $p_2=1$).
  It remains to compare the inhomogeneous norms in \eqref{gGNineq} with the $\W$-norms in \eqref{mainresGG}.
  Since $1-\gamma \theta-\eps(1-\theta)>H+\frac12+\eps>\frac12$ by \eqref{genGNcond}, and $1-\gamma>\frac1{p_1}$, \eqref{eq:equivWW} and \eqref{eq:trivialW} give
  \begin{align*}
    \|\phi\|_{\W^{H+\frac{1}{2}+\varepsilon}_2([0,T])}\le C \|\phi\|_{W^{1-\gamma\theta - \epsilon ( 1- \theta)}_2([0,T])}, \quad
    \|\phi\|_{W^{1-\gamma}_{p_1}([0,T])}\le C \|\phi\|_{\W^{1-\gamma}_{p_1}([0,T])}.
  \end{align*}
  For the last term in \eqref{gGNineq}, we use the elementary embedding $W^{1}_{p_2}([0,T])\hookrightarrow W^{1-\eps}_{p_2}([0,T])$ of inhomogeneous spaces, valid for every $p_2\ge1$, and the Poincar\'e inequality $\|\phi\|_{L_{p_2}([0,T])}\le T\|\phi'\|_{L_{p_2}([0,T])}$ for $\phi_0=0$, which give $\|\phi\|_{W^{1-\eps}_{p_2}([0,T])}\le C(1+T)\|\phi\|_{\W^{1}_{p_2}([0,T])}$. This avoids the space $\W^{1-\eps}_{p_2}$, which is not defined for $p_2=1$. Combining these bounds with \eqref{gGNineq} yields \eqref{mainresGG}.
  \end{proof}
  
The Sobolev spaces $\W^\alpha_m$ are relevant because they are continuously embedded into the Cameron--Martin spaces $\H^G$. In fact, for $H\in(0,\frac12)$ and any $\eps\in(0,\frac12-H]$ there exists a constant $C=C(\eps,T)$ such that
\begin{equation}\label{GGbound}
	\|\phi\|_{\H^G}\le C\|\phi\|_{\W^{H+\frac12+\eps}_2([0,T])} \quad \forall \phi\in \W^{H+\frac12+\eps}_2([0,T]).
\end{equation}
For $G=B^H$, this is \cite[Proposition~C.1(iv)]{galeatigerencser22}. For $G=\tilde B^H$, the operator $K^G$ is a constant multiple of the Riemann--Liouville integral of order $H+\frac12\in(\frac12,1)$, which by \cite[Theorem~3.1]{Gorenflo:2015vr} is an isomorphism from $L_2([0,T])$ onto the closed subspace of $W^{H+\frac12}_2([0,T])$ of functions vanishing at $0$; hence \eqref{GGbound} follows from \eqref{eq:equivWW} and \eqref{eq:trivialW}.
We rely on the following convenient embedding and the Gagliardo--Nirenberg inequality \eqref{mainresGG} to estimate Sobolev norms of stochastic processes.
\begin{lemma}
\label{lem:vmo_sobolev}
 Let $(V_t)_{t\in[0,T]}$ be a continuous adapted process, $V_0=0$; assume that there exist constants  $q\in[1,\infty)$, $\beta\in[1/q,1)$, $\Lambda>0$ and a control $w$ such that 
 \begin{align}
\label{eq:cond_sobolev_embedding}
  \E^s|V_t-V_s| \leq \Lambda  w(s,t)^{\frac{1}{q}}(t-s)^{\beta-\frac{1}{q}}, \quad\forall (s,t)\in\Delta_{[0,T]}.
 \end{align}
  Then for any $\tilde \beta\in(0,\beta)$ with $\tilde\beta>\frac1q$ there exist constants $\mu>0$ and $C<\infty$, depending only on $\beta$, $\tilde\beta$ and $q$ (in particular, not on $V$, $\Lambda$, $w$ or $T$), such that
  \begin{align}
\label{eq:vmo_sobolev_embedding}
  \E \left[\exp \mu \left(\Lambda^{-1}T^{\tilde{\beta}-\beta}w(0,T)^{-\frac{1}{q}}  \norm{V}_{\W^{\tilde \beta}_q([0,T])}\right)^{\frac{1}{1-\beta}} \right]\le C.
  \end{align}
\end{lemma}
\begin{proof}
  In view of the John--Nirenberg inequality \cite[Corollary~2.5]{le22}, \eqref{eq:cond_sobolev_embedding} implies that for any $m\ge 1$,
  \begin{align*}
    \|V_t-V_s\|_{L_m(\Omega)|\F_s}\le C_m \Lambda w(s,t)^{\frac{1}{q}}(t-s)^{\beta-\frac{1}{q}}, \quad\forall (s,t)\in\Delta_{[0,T]},
  \end{align*}
    where $C_m= C(q,\beta)\,\Gamma(m(1-\beta)+1)^{\frac1m}$, $C(q,\beta)$ is a constant and $\Gamma$ is the Gamma function. 
  We then apply \cite[Proposition 2.7]{LLAfbm} to obtain that 
  \begin{align*}
        \|\norm{V}_{\W^{\tilde \beta}_q([0,T])}\|_{L_m(\Omega)}\le C C_m \Lambda T^{\beta-\tilde{\beta}} w(0,T)^{\frac{1}{q}}
  \end{align*}
  for $m\ge q$, with $C=C(\beta,\tilde\beta,q)$.
  Writing $X$ for the random variable inside the exponential in \eqref{eq:vmo_sobolev_embedding} and taking $m=\frac{k}{1-\beta}$, we get $\E[(\mu X^{\frac1{1-\beta}})^k]\le(\mu C')^k\Gamma(k+1)$ for all large $k$, with $C'=C'(\beta,\tilde\beta,q)$. Expanding the exponential, \eqref{eq:vmo_sobolev_embedding} follows for $\mu<1/C'$, with a constant depending only on $\beta$, $\tilde\beta$ and $q$.
\end{proof}
\subsection{Linear Young differential equations}
We state a version of the rough Gr\"onwall lemma taken from \cite[Lem.~A.2]{HLN2021}. 

\begin{lemma}\label{lem:rough-Gr\"onwall}
	Let $G:[0,T]\to [0,+\infty)$ be a path such that for some constants $C$, $L>0$, $\kappa\geq 1$, some control $w$ and some increasing function $\Upsilon$ on $[0,T]^2_{\leq}$ (i.e. $\Upsilon(s,t)\le \Upsilon(u,v)$ whenever $[s,t]\subset[u,v]$), one has 
	\begin{equation}\label{eq:rough-Gr\"onwall-hypothesis}
		G_t-G_s \leq C \, \sup_{r\leq t}\, G_r \, w(s,t)^{\frac{1}{\kappa}} + \Upsilon(s,t)
	\end{equation}
	for every $s<t$ satisfying $w(s,t)\leq L$. Then, it holds that
	\begin{equation}\label{eq:rough-Gr\"onwall-conclusion}
		\sup_{t\in [0,T]} G_t \leq 2 \exp\bigg( \frac{w(0,T)}{\lambda L}\bigg) \Big( G_0 + \Upsilon(0,T) \Big)
	\end{equation}
	where $\lambda= 1 \wedge \big[L (2C e^2)^{\kappa}\big]^{-1}$.
\end{lemma}
\begin{lemma}\label{lem:linearYDE}
  Let $Z,V,R$ be continuous functions on $[0,T]$ with finite $\ell$-variation for some $\ell\in(1,2)$. Assume that $dZ_t = Z_tdV_t +dR_t$.
  Then there exist finite positive constants $c=c(\ell)$, $C=C(\ell)$ such that
  \begin{align}\label{eq:Zbound}
    &\sup_{t\in[0,T]}|Z_t| \leq 4e^{c\bra{V}^\ell_{\ell\var;[0,T]}}(|Z_0|+\bra{R}_{\ell\var;[0,T]}),
    \\&\bra{Z}_{\ell\var;[s,t]} \leq C \left( \sup_{r\in[0,T]}|Z_r| \bra{V}_{\ell\var;[s,t]}   + \bra{R}_{\ell\var;[s,t]} \right) \quad\forall s\le t\le T.\label{eq:Zboundvar}
  \end{align}
\end{lemma}
\begin{proof}
  Standard estimates for Young integrals yield that
  \begin{align*}
    |\int_s^t Z_rdV_r| \le C (|Z_s|+\bra{Z}_{\ell\var;[s,t]}) \bra{V}_{\ell\var;[s,t]}, \quad \forall s\le t\le T,
  \end{align*}
  for some universal constant $C=C(\ell)$.
  We observe that $(s,t)\mapsto \bra{V}_{\ell\var;[s,t]}^\ell$ is a control, and that if a continuous function $Y$ satisfies $|Y_t-Y_s|\le C w(s,t)^{1/\ell}$ for some control $w$, then $\bra{Y}_{\ell\var;[s,t]} \le C w(s,t)^{1/\ell}$. It follows that
  \begin{align*}
    \bra{\int_0^\cdot Z_rdV_r}_{\ell\var;[s,t]} \le C (|Z_s|+\bra{Z}_{\ell\var;[s,t]}) \bra{V}_{\ell\var;[s,t]}
  \end{align*}
  and hence, using the equation for $Z$, 
  \begin{align*}
     \bra{Z}_{\ell\var;[s,t]} \leq C \left( |Z_s| + \bra{Z}_{\ell\var;[s,t]} \right) \bra{V}_{\ell\var;[s,t]}+\bra{R}_{\ell\var;[s,t]} .
  \end{align*} 
  Now choosing $(s,t)$ such that 
  \begin{align*}
    w(s,t)^{\frac{1}{\ell} }:=\bra{V}_{\ell\var;[s,t]} \le \frac{1}{2C}, 
  \end{align*}
  we find that
  \begin{align}\label{eq:Zvarsmallst}
    \bra{Z}_{\ell\var;[s,t]} \leq 2C |Z_s|w(s,t)^{\frac{1}{\ell} } + 2\bra{R}_{\ell\var;[s,t]}  .
  \end{align}
  Set $G_t=\sup_{r\in[0,t]}|Z_r|$. We see that 
  \begin{align*}
    G_t-G_s\le \sup_{r\in[s,t]}|Z_r-Z_s|\le \bra{Z}_{\ell\var;[s,t]} \leq 2C G_s w(s,t)^{\frac{1}{\ell} } + 2\bra{R}_{\ell\var;[s,t]} .
  \end{align*}
  Applying the rough Gr\"onwall lemma (\cref{lem:rough-Gr\"onwall}), we obtain that
  \begin{align*}
    \sup_{t\in[0,T]}|Z_t|=  G_T \leq 2 e^{c w(0,T)}  (G_0+2\bra{R}_{\ell\var;[0,T]}) 
  \end{align*}
  which yields \eqref{eq:Zbound}. 

  Let $s<t\le T$ be arbitrary and let $(t_i)$ be a partition $(t_i)$ of $[s,t]$ and consider the quantity
  \begin{align*}
    \sum_{t_i}|Z_{t_{i+1}}-Z_{t_i}|^\ell. 
  \end{align*}
  The number of the intervals $[t_i,t_{i+1}]$ satisfying $w(t_i,t_{i+1})>1/(2C)^\ell$ is at most $(2C)^\ell w(s,t)$ and each of these contributes at most $2^\ell G_T^\ell$.  On the other intervals with $w(t_i,t_{i+1})\le 1/(2C)^\ell$, we use \eqref{eq:Zvarsmallst} to get that
  \begin{align*}
    \sum_{t_i}|Z_{t_{i+1}}-Z_{t_i}|^\ell\le  
    C_1( G_T^\ell w(s,t)+ \bra{R}_{\ell\var;[s,t]}^\ell). 
  \end{align*}
  for a suitable constant $C_1=C_1(\ell)$. Since $(t_i)$ is arbitrary, we can take the supremum over all partitions of $[s,t]$ and obtain \eqref{eq:Zboundvar}.
\end{proof}

\section{Key bounds on integral functionals}\label{sec:int_bd}

In this section, we obtain several bounds on integral functionals of the fractional Brownian motion and the solution to equation \ref{fSDE}, which are crucial for proving the main results.
Throughout the section, we will write $\B^{\alpha}_{\rho}:= \B^{\alpha}_{\rho, \infty}$  for 
the Besov space of regularity $\alpha \in \mathbb{R}$ and integrability $p \in [1, \infty]$, 
see \cite{BCD11} or \cite[Appendix]{athreya2021wellposedness} for further details.

\subsection{Additive functionals of Gaussian processes}\label{subsec:add_fc_fbm}
We work in a more general setting than in the main results, where we fix a $d$-dimensional centered Gaussian process $(G_t)$ with its natural filtration $(\F^G_t)$. 
We assume that $G$ satisfies the following extended  local-nondeterminism property: there exist constants $H\in(0,1)$ and $\lambda>0$ such that
\begin{equation}\label{con:local_nondeterminism}
  \Var \left(\E^u G_t  \big|  \F^G_s\right)\ge \lambda^{-1} ((t-s)^{2H}-(t-u)^{2H})\mathbb{I} ,\quad s\le u\le t,
\end{equation}
where $\Var (\E^u G_t |\F^G_s)=\E^s (\E^u G_t-\E^s G_t)^{\otimes2}$ is the conditional variance of $\E^uG_t$ given $\F^G_s$, and the inequality is understood in the sense of symmetric matrices.
More generally, all results of this subsection remain valid, with the same proofs, if $(\F^G_t)$ is replaced by a filtration $(\F_t)$ to which $G$ is adapted, provided that $G$ is Gaussian conditionally on $\F_u$ with a deterministic conditional covariance, for every $u$, and that \eqref{con:local_nondeterminism} holds with respect to $(\F_t)$.
\begin{lemma}\label{lem:LND_examples}
  Let $H\in(0,\frac12)$, let $W$ be an $(\F_t)$-Brownian motion and let $G\in\{B^H,\tilde B^H\}$ be given by $G_t=\int_0^tK^G(t,r)\,dW_r$, see \cref{subsec:girsanov_prelim}. Then $G$ satisfies the assumptions of this subsection with respect to $(\F_t)$; in particular, \eqref{con:local_nondeterminism} holds with $(\F^G_t)$ replaced by $(\F_t)$, for some $\lambda=\lambda(H)$.
\end{lemma}
\begin{proof}
  For $u\le t$, since the increments of $W$ after time $u$ are independent of $\F_u$, we have $\E[G_t|\F_u]=\int_0^uK^G(t,r)\,dW_r$, and $G_t-\E[G_t|\F_u]=\int_u^tK^G(t,r)\,dW_r$ is a centered Gaussian vector independent of $\F_u$, with deterministic covariance.
  Hence, for $s\le u\le t$, $\E[G_t|\F_u]-\E[G_t|\F_s]=\int_s^uK^G(t,r)\,dW_r$ is independent of $\F_s$, and its covariance matrix is $\int_s^uK^G(t,r)^2\,dr\,\mathbb{I}$.
  For $G=B^H$ and $H<\frac12$, we have $K_H(t,r)\ge C(t-r)^{H-\frac12}$ for $0\le r\le t$ with $C=C(H)>0$ \cite[Proposition~B.2]{BLMytnik}, so that $\int_s^uK_H(t,r)^2\,dr\ge \frac{C^2}{2H}((t-s)^{2H}-(t-u)^{2H})$. For $G=\tilde B^H$, we have $\int_s^u\tilde c_H^2(t-r)^{2H-1}\,dr=\frac{\tilde c_H^2}{2H}((t-s)^{2H}-(t-u)^{2H})$.
\end{proof}

For a distribution $f$ in $L_q\B^\alpha_\rho$ for some $\alpha< 0$, $q\in[1,\infty]$, $\rho\in[1,\infty]$, we consider the additive functional $U^f$ defined as
\begin{align}
    U^f=U^f(G)\defby\lim_n \int_0^\cdot f^n(r,G_r)dr \quad\text{in probability, uniformly on }[0,T], \label{eq:Uffn}
\end{align}
for any sequence $(f^n)$ of measurable functions in $L_1([0,T];B(\Rd))$ such that $\lim_n f^n= f$ in  $L_q([0,T];\B^{\alpha-\varepsilon}_\rho)$ for some sufficiently small  $\varepsilon\in(0,1)$. Here and below, $B(\Rd)$ denotes the space of bounded measurable functions on $\Rd$ equipped with the supremum norm. That \eqref{eq:Uffn} is well-defined for suitable parameters $(\alpha,q,\rho)$ needs justifications which will be our main focus in this subsection.

The main result herein is \cref{lem:integral_along_fbm}. It shows particularly that if
$f$ belongs to $L_q([0,T]; \B^{\alpha}_\rho(\R^d))$, then under suitable  conditions on $\alpha$, $q$, $\rho$, the additive functional $U^f$ is well-defined and 
is a.s. of finite $(\frac1{1+\alpha H-Hd(\frac1\rho-\frac12)^+} + \eps)$-variation for every $\eps>0$; in particular, for $\rho\ge2$, of finite $(\frac1{1+\alpha H} + \eps)$-variation.
Even when $f$ is a measurable function, this result significantly improves the corresponding estimates in \cite[Lemma~2.8]{HRZ}, \cite[Lemma~4.2]{BLMytnik}, and \cite[Lemma~3.8]{BGal23}, where it was shown that $U^f$ has finite $(\frac1{1+\alpha H-\frac{Hd}{\rho}}+\eps)$-variation. This gain in the regularity of $U^f$ plays a crucial role in showing our main results.

This is accomplished in several steps.
First, we obtain estimates with singularities on the moments of increments of $U^f$ conditionally on $\F_{S_0}$ for a fixed time $S_0$ (\cref{L:firstb}). 
These estimates, which are more refined than those available in the existing literature (\cite{BLMytnik,BGal23}), show that the $\rho$th conditional moment of the increment of $U^f$ is bounded by a control raised to the power $\frac{1}{1+\alpha H}$. 
Next in  \cref{lem:integral_along_fbm}, we use the Kolmogorov continuity theorem in $p$-variation form (\cref{lem:variation-kolmogorov}) to show that the $(\frac{1}{1+\alpha H}+\eps)$-variation of $U^f$ has a finite conditional first moment. Crucially, unlike the usual Kolmogorov continuity theorem, the $p$-variation version does not reduce the a.s. regularity of the function. Through the John–Nirenberg inequality,  the conditional first moment ultimately yields estimates   for all moments of the $(\frac{1}{1+\alpha H}+\eps)$-variation of $U^f$.

For a symmetric positive definite $d\times d$ matrix $\Sigma$, let $p_\Sigma$ be the density of a centered normal distribution with covariance matrix $\Sigma$ defined by 
\begin{equation}\label{eq:gaussian_density}
p_\Sigma(x)=\frac1{(2\pi)^{d/2}(\det\Sigma)^{1/2}}\exp\left(-\frac12 x^\top\Sigma^{-1}x\right),\quad x\in\R^d.
\end{equation}
For a measurable function $h$ on $\R^d$, we write $P_{\Sigma}h=p_{\Sigma}*h$. 
For each $t$, we denote $\E^t=\E(\ \cdot \ | \F^G_t)$.
\begin{lemma}\label{lem.momentfG}
  Let $\alpha<0$ and $\rho\in[1,\infty]$ be some fixed constants. Then there exists a  constant $C>0$ such that  for every $s\le u\le t$ and every bounded measurable function $f\colon\R^d\to\R$, 
  \begin{equation}\label{eq:bound_Efg}
 \bigl\|\E^u f(G_t)\bigr\|_{L_\rho(\Omega)|\F_s}\le C \|f\|_{\B^{\alpha}_\rho}(t-u)^{\alpha H}\bigl((u-s)^{-\frac{Hd}{\rho}}+\I_{\{H<1/2\}}(u-s)^{-\frac{d}{2\rho}}(t-u)^{\frac{d}{2\rho}-\frac{Hd}\rho}\bigr).
 \end{equation}
\end{lemma}
\begin{proof}
  We can assume without loss of generality that $\|f\|_{\B^\alpha_\rho}<\infty$. 
  For $s\le u \le t$, define $V_{s,u,t}$ and $\Sigma_{s,u,t}$ by  $\E^u G_t=\E^s G_t+V_{s,u,t}$ and $\Sigma_{s,u,t}=\E^s[V_{s,u,t}^{\otimes2}]$.
  Note that conditioned on $\F^G_u$,  $t\mapsto G_t-\E^uG_t=V_{u,t,t}$ is a centered Gaussian process.
  Hence, we have 
  \begin{align*}
    \E^u f(G_t)=P_{\Sigma_{u,t,t}} f(\E^u G_t)=P_{\Sigma_{u,t,t}} f(\E^s G_t+V_{s,u,t}).
  \end{align*}
 For  $\rho\neq\infty$, we have
\begin{align*}
	\bigl\|\E^u f(G_t)\bigr\|_{L_\rho(\Omega)|\F_s}^\rho
	&=\int_{\R^d} |P_{\Sigma_{u,t,t}} f(\E^s G_t+x)|^\rho p_{\Sigma_{s,u,t}}(x)\,dx \\
	&\le \|p_{\Sigma_{s,u,t}}\|_{L_\infty(\R^d)}\|
	P_{\Sigma_{u,t,t}} f\|_{L_\rho(\R^d)}^\rho.
\end{align*}	
Hence, we obtain that 
\begin{align*}
	\bigl\|\E^u f(G_t)\bigr\|_{L_\rho(\Omega)|\F_s}
	&\le \|p_{\Sigma_{s,u,t}}\|_{L_\infty(\R^d)}^{1/\rho} \|
	P_{\Sigma_{u,t,t}} f\|_{L_\rho(\R^d)},
\end{align*}	
which also holds for $\rho=\infty$.
The condition \eqref{con:local_nondeterminism} ensures that $\Sigma_{s,u,t}\ge\lambda^{-1} ((t-s)^{2H}-(t-u)^{2H})\mathbb{I}$ and $\Sigma_{u,t,t}\ge\lambda^{-1}(t-u)^{2H}\mathbb{I}$,
  which implies that
\begin{align*}
  \|p_{\Sigma_{s,u,t}}\|_{L_\infty(\R^d)}\les ((t-s)^{2H}-(t-u)^{2H})^{-\frac{d}2}
  \tand
  \|P_{\Sigma_{u,t,t}} f\|_{L_\rho(\R^d)}
	&\les(t-u)^{\alpha  H}\|f\|_{\B^\alpha_\rho}.
\end{align*}
The rest of the proof is the same as \cite[Proposition~B.2(iv)]{BLMytnik}.
\end{proof}

\begin{lemma}\label{L:firstb}
Let $f\colon[0,T]\times\R^d\to\R$ be a measurable function which belongs to $L_1([0,T];B(\Rd))$ and to $L_q ([0,T];\B^{\alpha}_\rho(\R^d))$, where $q\in[1,\infty]$, $\rho\in[2,\infty]$, $\alpha< 0$. If $\rho=\infty$, the conditional norms $\|\cdot\|_{L_\rho(\Omega)|\F_{S_0}}$ below are replaced by $\|\cdot\|_{L_m(\Omega)|\F_{S_0}}$ for an arbitrary $m\in[2,\infty)$, and the constants depend also on $m$.
\begin{enumerate}[(i)]
  \item Assume that 
  \begin{equation}\label{ahcond}
    \alpha H>\max(-\frac12,\frac1q-1).
  \end{equation}
  Then there exists a constant $C=C(\alpha,\rho)>0$ such that for any $0\le S_0<  S\le T$  one has
\begin{equation}\label{mainbound1}
\Bigl\|\,\Bigl\|\int_{S}^{T} f(r,G_r)dr\Bigr\|_{L_\rho(\Omega)|\F_{S_0}}\Bigr\|_{L_\infty(\Omega)}\le  C \|f\|_{L_q ([S,T];\B^{\alpha}_\rho)}(S-S_0)^{-\frac{Hd}\rho}(T-S)^{1+\alpha H-\frac1q}		.
\end{equation}
\item Assume  that 
\begin{equation}\label{condnew}
  \alpha H>\max\bigl(-\frac12,\frac1q-1+\frac{Hd}{\rho}\bigr).
\end{equation}
Define
\begin{equation}
\gamma:=\frac{1+\alpha H-\frac1q-\frac{Hd}{\rho}}{1+\alpha H-\frac1q}\in(0,1].\label{gammag}
\end{equation}
Then there exists a constant $C=C(\alpha,\rho)>0$ such that for any $0\le S_0\le  S\le T$  one has
\begin{equation}\label{mainbound2}
  \Bigl\|\,\Bigl\|\int_{S}^{T} f(r,G_r)dr\Bigr\|_{L_\rho(\Omega)|\F_{S_0}}\Bigr\|_{L_\infty(\Omega)}\le  C \|f\|_{L_q ([S,T];\B^{\alpha}_\rho)}((T-S_0)^{\gamma}-(S-S_0)^{\gamma})^{1+\alpha H-\frac1q}.
\end{equation}
Further, we have for any $m\ge1$
\begin{equation}\label{mainbound2ob}
\Bigl\|\int_{S}^{T} f(r,G_r)dr\Bigr\|_{L_m(\Omega)}
  \le  C \|f\|_{L_q ([S,T];\B^{\alpha}_\rho)}(T-S)^{1+\alpha H-\frac1q-\frac{Hd}{\rho}},
\end{equation}
where $C=C(\alpha,\rho,m)>0$.
\end{enumerate}
\end{lemma}
\begin{proof}
\textbf{Step 1}. Fix $0\le S_0<  S\le T$.
We are going to apply stochastic sewing lemma (\cref{l:SSL}) with $(m,n):=(\rho,\infty)$ (if $\rho=\infty$, with $(m,\infty)$ for a fixed finite $m\ge2$, using that $\|\cdot\|_{L_m(\Omega)|\F_{S_0}}\le\|\cdot\|_{L_\infty(\Omega)|\F_{S_0}}$ in \eqref{Astl41}) to the following processes
\begin{align*}
\A_t:=\int_S^{t}f(r,G_r)\,dr,\quad A_{s,t}:=\E^s\int_s^t f(r,G_r)dr,\quad (s,t)\in\Delta_{[S,T]}.
\end{align*}
Applying \cref{lem.momentfG} and the H\"older inequality, we  derive for  $(s,t)\in\Delta_{[S,T]}$
	\begin{align}
		&\|\,\|  A_{s,t}\|_{L_\rho(\Omega)|\F_{S_0}}\|_{L_\infty(\Omega)}\nn\\
		&\quad\le \int_s^t\|\, \|\E^s f(r, G_r)\|_{L_\rho(\Omega)|\F_{S_0}}\|_{L_\infty(\Omega)}\, dr\nn\\
		&\quad\le C \int_s^t \|f_r\|_{\B^{\alpha}_\rho}(r-s)^{\alpha H}(S-S_0)^{-\frac{Hd}\rho}\bigl(1 +\I_{\{H<1/2\}}
		 \frac{(r-s)^{ \frac{d}{2\rho}-\frac{Hd}{\rho}}}{(S-S_0)^{\frac{d}{2\rho}-\frac{Hd}{\rho}}}\bigr)\,dr\label{Astl41spec}\\
		&\quad\le C \|f\|_{L_q ([s,t];\B^{\alpha}_\rho)}(t-s)^{1+\alpha H-\frac1q}(S-S_0)^{-\frac{Hd}\rho}
		\bigl(1 +\I_{\{H<1/2\}}
		 \frac{(t-s)^{ \frac{d}{2\rho}-\frac{Hd}{\rho}}}{(S-S_0)^{\frac{d}{2\rho}-\frac{Hd}{\rho}}}\bigr),\label{Astl41}	\end{align}
where $C=C(\alpha,\rho)$. Here we used condition \eqref{ahcond} to ensure that ${\int_s^t (r-s)^{\alpha H q'}dr<\infty}$ for $q'\in[1,\infty]$ such that $\frac1q+\frac1{q'}=1$.	

If $q<\infty$, we put $w(s,t):=\|f\|_{L_q ([s,t];\B^{\alpha}_\rho)}^q$, which defines a control (by \eqref{controlmaps}). The factor $\|f\|_{L_q ([s,t];\B^{\alpha}_\rho)}$ in the right-hand side of \eqref{Astl41} thus can be rewritten as $w(s,t)^{1/q}$. Thanks to \eqref{ahcond}, $1+\alpha H-\frac1q>0$ and $\frac1q+(1+\alpha H-\frac1q)>\frac12$; further if $H<\frac12$, then    
	$\frac1q+(1+\alpha H-\frac1q+\frac{d}{\rho}(\frac12-H))>\frac12$. Thus, condition \eqref{SSLcond1} holds.

If $q=\infty$, then we just use the trivial bound $\|f\|_{L_\infty ([s,t];\B^{\alpha}_\rho)}\le \|f\|_{L_\infty ([S,T];\B^{\alpha}_\rho)}$. We see that in this case  \eqref{SSLcond1} also holds since $1+\alpha H-\frac1q=1+\alpha H>\frac12$ by \eqref{ahcond}.

Next, note that for any $(s,t)\in\Delta_{[S,T]}$, $u\in[s,t]$, we have
\begin{equation}\label{Part02}
\E^s \delta A_{s,u,t}=\E^s \Bigl(\E^s \int_u^t f(r,G_r)\, dr-\E^u \int_u^t f(r,G_r)\, dr\Bigr)=0.
\end{equation}
Thus, \eqref{SSLcond2} holds.
	
Finally, let us verify condition \eqref{SSLcond3}. Fix $(s,t)\in\Delta_{[S,T]}$ and consider a partition $\Pi_N:=\{t_0=s,t_1,\hdots,t_N=t\}$ of the interval $[s,t]$.  It is easy to see that 
$$
\E^{t_i}(\A_{t_{i+1}}-\A_{t_{i}}-A_{t_i,t_{i+1}})=0,\quad i=0,\hdots,N-1,
$$
and $\A_{t_{i+1}}-\A_{t_{i}}-A_{t_i,t_{i+1}}$ is $\F_{t_{i+1}}$--measurable.   Therefore,  $\A_t-\A_s-\sum_{i=0}^{N-1} A_{t_i,t_{i+1}}$ is a sum of martingale differences. Then, by the martingale identity 
\begin{align*}
	\Bigl\|\sum_{i=0}^{N-1} (\A_{t_{i+1}}-\A_{t_i}-A_{t_i,t_{i+1}})\Bigr\|_{L_2(\Omega)}^2
	&= \sum_{i=0}^{N-1}\|\A_{t_{i+1}}-\A_{t_{i}}- A_{t_i,t_{i+1}}\|_{L_2(\Omega)}^2\\
    &\les \sum_{i=0}^{N-1}\left( \int_{t_i}^{t_{i+1}}\|f(r,\cdot)\|_{B(\Rd)}dr   \right) ^2  \\
	&\les \sup_{i}\int_{t_i}^{t_{i+1}}\|f(r,\cdot)\|_{B(\Rd)}dr  \int_{0}^{T}\|f(r,\cdot)\|_{B(\Rd)}dr.
\end{align*}
This implies that $\sum_{i=0}^{N-1} A_{t_i,t_{i+1}}$ converges to $\A_{t}-\A_s$ in $L_2(\Omega)$, hence in probability, as $|\Pi_N|\to0$, and \eqref{SSLcond3} holds.

Thus, all the conditions of \cref{l:SSL} are satisfied with $m=\rho$, $n=\infty$ and we derive for $(s,t)\in\Delta_{[S,T]}$
	\begin{align}\label{bef2sing}
		&\Bigl\|\,\Bigl\|\int_{s}^{t}  f(r,G_r)\,dr\Bigr\|_{L_\rho(\Omega)|\F_{S_0}}\Bigr\|_{L_\infty(\Omega)}\nn\\
		&\,\,\le 
		C \|f\|_{L_q ([S,T];\B^{\alpha}_\rho)}(t-s)^{1+\alpha H-\frac1q}(s-S_0)^{-\frac{Hd}{\rho}} 
		\bigl(1 +\I_{\{H<1/2\}}
		 \frac{(t-s)^{ \frac{d}{2\rho}-\frac{Hd}{\rho}}}{(s-S_0)^{\frac{d}{2\rho}-\frac{Hd}{\rho}}}\bigr),
	\end{align}
where $C=C(\alpha,\rho)$.

\textbf{Step 2}. Now we simplify the right hand side  of \eqref{bef2sing} by applying \cref{lem:singularity}\ref{lsone}. Let $E$ be the space of all random variables $\xi\colon\Omega\to\R^d$ such that 
$$
\|\,\|\xi\|_{L_\rho(\Omega)|\F_{S_0}}\|_{L_\infty(\Omega)}<\infty.
$$
 For $\xi_1,\xi_2\in E$ define $d(\xi_1,\xi_2):=
\|\,\|\xi_1-\xi_2\|_{L_\rho(\Omega)|\F_{S_0}}\|_{L_\infty(\Omega)}$. For $S-S_0\le t\le T-S_0$ put
$$
Y_t:=\int_{S_0}^{S_0+t} f(r,G_r)\,dr.
$$
Then for $(s,t)\in\Delta_{[S-S_0,T-S_0]}$ we can rewrite \eqref{bef2sing} as 
\begin{align*}
d(Y_s,Y_t)&=\Bigl\|\,\Bigl\|\int_{S_0+s}^{S_0+t}f(r,G_r)\,dr\Bigr\|_{L_\rho(\Omega)|\F_{S_0}}
\Bigr\|_{L_\infty(\Omega)}\\
&\le 	C \|f\|_{L_q ([S,T];\B^{\alpha}_\rho)}(t-s)^{1+\alpha H-\frac1q}s^{-\frac{Hd}{\rho}} 
		\bigl(1 +\I_{\{H<1/2\}}
		 \frac{(t-s)^{ \frac{d}{2\rho}-\frac{Hd}{\rho}}}{s^{\frac{d}{2\rho}-\frac{Hd}{\rho}}}\bigr),
\end{align*}
where $C=C(\alpha,\rho)$. By \eqref{restame1}, this implies 
\begin{align*}
\Bigl\|\,\Bigl\|\int_{S}^{T}f(r,G_r)\,dr\Bigr\|_{L_\rho(\Omega)|\F_{S_0}}
\Bigr\|_{L_\infty(\Omega)}&=d(Y_{S-S_0},Y_{T-S_0})\\
&\le  	C \|f\|_{L_q ([S,T];\B^{\alpha}_\rho)}(T-S)^{1+\alpha H-\frac1q}(S-S_0)^{-\frac{Hd}{\rho}},
\end{align*}
where $C=C(\alpha,\rho)$. This yields \eqref{mainbound1}.

\textbf{Step 3}. We now prove \eqref{mainbound2}.
By \eqref{mainbound1}, there exists $C=C(\alpha,\rho)$, such that for any $S-S_0< s \le t\le T-S_0$ we have 
\begin{equation*}
d(Y_s,Y_t)=\Bigl\|\,\Bigl\|\int_{S_0+s}^{S_0+t}f(r,G_r)\,dr\Bigr\|_{L_\rho(\Omega)|\F_{S_0}}
\Bigr\|_{L_\infty(\Omega)}\le  C \|f\|_{L_q ([S,T];\B^{\alpha}_\rho)}(t-s)^{1+\alpha H-\frac1q}
s^{-\frac{Hd}\rho}.
\end{equation*}
We see that the process $Y$ is continuous (because $f$ is bounded).
By \eqref{condnew}, $1+\alpha H-\frac1q>\frac{Hd}\rho$. Thus,  all the conditions of \cref{lem:singularity}\ref{lstwo} are satisfied for $\tau=1+\alpha H-\frac1q$, $\eta=\frac{Hd}\rho$ and  \eqref{restame} implies that
\begin{equation}\label{Hlargets}
	d(Y_s,Y_t)\le  C \|f\|_{L_q ([S,T];\B^{\alpha}_\rho)}(t^{\gamma}-s^\gamma)^{1+\alpha H-\frac1q},\quad  S-S_0\le s\le t\le T-S_0,
\end{equation}
where $C=C(\alpha,\rho)$ and $\gamma$ is  as in \eqref{gammag}. 
We also have
\begin{equation*}
\Bigl\|\,\Bigl\|\int_{S}^{T}f(r,G_r)\,dr\Bigr\|_{L_\rho(\Omega)|\F_{S_0}}\Bigr\|_{L_\infty(\Omega)}=d(Y_{S-S_0},Y_{T-S_0}).
\end{equation*}
Substituting this into \eqref{Hlargets}, we derive \eqref{mainbound2}. 

\textbf{Step 4}.
To obtain \eqref{mainbound2ob} we use the John--Nirenberg inequality (\cref{prop:vmo_john_nirenberg}\ref{jnone}). We denote 
$$
\A_t:=\int_0^{t} f(r,G_r)\,dr,\quad t\in[S,T],
$$
which is a continuous process because $f$ is bounded. Furthermore, by taking  $S=S_0$ in \eqref{mainbound2}, we deduce that for any $(s,t)\in\Delta_{[S,T]}$
\begin{equation*}
\E^s|\A_t-\A_s|= \E^s \bigl| \int_{s}^{t} f(r,G_r)\,dr\bigr|\le  C \|f\|_{L_q ([S,T];\B^{\alpha}_\rho)}(T-S)^{1+\alpha H-\frac1q-\frac{Hd}{\rho}}.
\end{equation*}
Thus, condition \eqref{eq:vmo_gamma} holds. Therefore, all conditions of \cref{prop:vmo_john_nirenberg}\ref{jnone} are satisfied and \eqref{jnres} implies
\begin{equation*}
\|\A_T-\A_S\|_{L_m(\Omega)}\le  C \|f\|_{L_q ([S,T];\B^{\alpha}_\rho)}(T-S)^{1+\alpha H-\frac1q-\frac{Hd}{\rho}},
\end{equation*}
which is  \eqref{mainbound2ob}.
\end{proof}


In the next result, we estimate moments of the variation of $U^f$.
\begin{proposition}\label{lem:integral_along_fbm}
  Let $f$ be an element of $L_q ([0,T];\B^{\alpha}_\rho(\R^d))$, where $q\in[1,\infty]$, $\rho\in[1,\infty]$, $\alpha< 0$.
  Assume  that 
\begin{equation}\label{con.for:U_var}
  \alpha H>\max\left(-\frac12+Hd\bigl(\frac1\rho-\frac12\bigr)^+,-1+\frac1q+\frac{Hd}{\rho}\right).
\end{equation}
Then $U^f$ is well-defined by \eqref{eq:Uffn}. Additionally, for any $\ell\in \left(\frac1{1+ \alpha H-Hd(\frac{1}{\rho}-\frac{1}{2})^+},\ 2\right)$,
 $m\ge1$ there exists a constant $C=C(\alpha,\ell,\rho,m)$ such that for any $(s,t)\in\Delta_{[0,T]}$ we have 
\begin{equation}\label{eq:U_rho_var}
\bigl\|\| U^f\|_{\ell\var;[s,t]}\bigr\|_{L_m(\Omega)|\F_s}\le C \|f\|_{L_q ([s,t];\B^{\alpha}_\rho)}(t-s)^{1+\alpha H-\frac1q-\frac{Hd}\rho}.
\end{equation}
\end{proposition}
\begin{proof}
  When $\rho<2$, we use the embedding $\B^\alpha_\rho\hookrightarrow\B^{\bar \alpha}_2$ with $\bar{\alpha}=\alpha-\frac{d}{\rho}+\frac{d}{2}$. Condition \eqref{con.for:U_var} implies that 
  \begin{align*}
    \bar{\alpha} H&>\max\left(-\frac12,-1+\frac1q+\frac{Hd}{2}\right).
  \end{align*} 
  Hence, the case $\rho<2$ follows from the case $\rho\ge2$ by replacing $\alpha$ with $\bar{\alpha}$ and $\rho$ with $2$. 
  We can thus assume that $\rho\ge2$. Furthermore, in the first two steps, we assume additionally that $f\in L_1([0,T];B(\Rd))$, in which case, \eqref{eq:Uffn} deduces that $U^f_t=\int_0^t f(r,G_r)dr$.

  \textbf{Step 1}. We first prove  \eqref{eq:U_rho_var} with $m=\rho$. Let $\eps>0$. Fix $S_0\in[0,T]$. 
  If $q<\infty$,
  we put, with $\gamma$ defined in \eqref{gammag},
  \begin{equation*}
  w_1(s,t):=\|f\|_{L_q ([s,t];\B^{\alpha}_\rho)}^q;\qquad w_2(s,t):=(t-S_0)^{\gamma}-(s-S_0)^{\gamma},\quad (s,t)\in\Delta_{[S_0,T]}.
  \end{equation*}
  We see that both $w_1$ and $w_2$ are controls. Hence, the function
  \begin{equation*}
    w_3:=w_1^{\frac{1/q}{1+\alpha H}}w_2^{\frac{1+\alpha H-1/q}{1+\alpha H}}
\end{equation*}
  is also a control by \eqref{frizcontrol}. Inequality \eqref{mainbound2} implies that
  \begin{equation}\label{tmp.estUfst}
    \bigl\|\,\|U^f_t-U^f_s\|_{L_\rho(\Omega)|\F_{S_0}}\bigr\|_{L_\infty(\Omega)}\le C w_3(s,t)^{1+\alpha H},\quad (s,t)\in\Delta_{[S_0,T]},
  \end{equation}
  for $C=C(\alpha,\rho)>0$. 
  If $q=\infty$, then  inequality \eqref{tmp.estUfst}
 also holds with 
$$
w_3(s,t):=\|f\|_{L_\infty([S_0,T];\B^{\alpha}_\rho)}^{\frac1{1+\alpha H}}((t-S_0)^{\gamma}-(s-S_0)^{\gamma}).
$$

  We apply now the Kolmogorov continuity theorem in the form of \cref{lem:variation-kolmogorov} to the process $U^f$. 
    We see by above that \eqref{condkol} is satisfied with $\G:=\F_{S_0}$, $w_3$ in place of $w$, $\mu:=1+ \alpha H$, $m=\rho$ (if $\rho=\infty$, with a finite $m$ such that $1+\alpha H>\frac1m$, see \cref{L:firstb}). We see that $1+\alpha H>\frac12$ by \eqref{condnew} and $\rho\ge2$. Thus, $1+\alpha H\in(\frac1\rho;1]$. Therefore, all conditions of  \cref{lem:variation-kolmogorov} are satisfied and \eqref{boundkolmi} implies that for any $\ell\in (1/(1+\alpha H),\rho)$,  there is a deterministic constant $C=C(\alpha,\ell,\rho)$ such that
  \begin{align}\label{condvar}
  \bigl\|\,\| U^f\|_{\ell\var;[S_0,T]}\bigr\|_{L_\rho(\Omega)|\F_{S_0}}
  &\le C  w_3(S_0,T)^{1+\alpha H}\nn\\
  &= C \|f\|_{L_q ([S_0,T];\B^{\alpha}_\rho)}(T-S_0)^{\gamma(1+\alpha H-\frac1q)}\nn\\
  &= C \|f\|_{L_q ([S_0,T];\B^{\alpha}_\rho)}(T-S_0)^{1+\alpha H-\frac1q-\frac{Hd}\rho}.
  \end{align}
  By taking expectation in the above inequality and recalling that $C$ does not depend on $\omega$,  we get \eqref{eq:U_rho_var} with $m=\rho$.

  \textbf{Step 2.} To obtain  \eqref{eq:U_rho_var} for arbitrary $m\ge1$, we apply the John--Nirenberg inequality (\cref{prop:vmo_john_nirenberg}\ref{jnone}).
    Fix $(\bar s,\bar t)\in\Delta_{[0,T]}$ and denote
  \begin{equation*}
    \A_t:=\| U^f\|_{\ell\var;[\bar s,t]},\quad t\in[\bar s,T].
  \end{equation*}
Then using  \eqref{condvar} and recalling \eqref{eq:variation_ineq}, we get   for any $(s,t)\in\Delta_{[\bar s,\bar t]}$
  \begin{equation}\label{eq:cA_vmo}
  \E^s|\A_t-\A_s|\le 
  \E^s\| U^f\|_{\ell\var;[s,t]}\le C \|f\|_{L_q ([\bar s,\bar t];\B^{\alpha}_\rho)}(\bar t-\bar s)^{1+\alpha H-\frac1q-\frac{Hd}\rho},
  \end{equation}
   where  $C=C(\alpha,\eps,\rho)$ is deterministic.  Note that the process $\A$ is continuous since $f$ is bounded. Therefore, applying \cref{prop:vmo_john_nirenberg}\ref{jnone}, we get for any $m\ge1$
  \begin{equation*}
  \bigl\|\| U^f\|_{\ell\var;[\bar s,\bar t]}\bigr\|_{L_m(\Omega)|\F_{\bar s}}
  \le\|\A_{\bar t}-\A_{\bar s}\|_{L_m(\Omega)|\F_{\bar s}}\le C \|f\|_{L_q ([\bar s,\bar t];\B^{\alpha}_\rho)}(\bar t-\bar s)^{1+\alpha H-\frac1q-\frac{Hd}\rho},
  \end{equation*}
 for a constant $C=C(\alpha,\eps,\rho,m)$,  which is \eqref{eq:U_rho_var}.

  \textbf{Step 3.} We remove the condition that $f\in L_1([0,T];B(\Rd))$. Let $(f^n)$ be any sequence in $L_1([0,T];B(\Rd))$ such that $\lim_nf^n= f$ in $L_q([0,T];\B^{\alpha-\varepsilon}_\rho(\Rd))$ for some $\varepsilon\in(0,1)$. By choosing $\varepsilon$ small enough, we can ensure that \eqref{con.for:U_var} holds with  $\alpha-\varepsilon$ in place of $\alpha$, and also that $\ell$ belongs to the range of \eqref{eq:U_rho_var} for $\alpha-\varepsilon$.
  By \eqref{eq:U_rho_var} with $\alpha-\varepsilon$ in place of $\alpha$, we have
  \begin{align*}
    &\bigl\|\| \int_0^\cdot (f^n-f^k)(r,G_r)dr \|_{\ell\var;[0,T]}\bigr\|_{L_m(\Omega)}
    \\&=\bigl\|\| U^{f^n}-U^{f^k}\|_{\ell\var;[0,T]}\bigr\|_{L_m(\Omega)}
    \les \|f^n-f^k\|_{L_q ([0,T];\B^{\alpha-\varepsilon}_\rho)}.
  \end{align*}
    Hence $(\int_0^\cdot f^n(r,G_r)dr)_n$ is a Cauchy sequence in $L_m(\Omega;C^{\ell\var}([0,T]))$, and its limit does not depend on the choice of $(f^n)$ (interlace two such sequences). This implies that $U^f$ is well-defined by \eqref{eq:Uffn}.  
  By Fatou's lemma and lower semi-continuity of variation, we have
 \begin{align*}
  \bigl\|\| U^{f}\|_{\ell\var;[s,t]}\bigr\|_{L_m(\Omega)|\F_s}
  &\le \liminf_n \bigl\|\| U^{f^n}\|_{\ell\var;[s,t]}\bigr\|_{L_m(\Omega)|\F_s}
  \\&\les \liminf_n\|f^n\|_{L_q ([s,t];\B^{\alpha-\varepsilon}_\rho)}(t-s)^{1+\alpha H-\frac1q-\frac{Hd}\rho}.
 \end{align*}
  Here we may take $f^n(r,\cdot)\defby P_{1/n}f(r,\cdot)$. Then $f^n$ belongs to $L_1([0,T];B(\Rd))$, $\lim_n f^n=f$ in $L_q([0,T];\B^{\alpha-\varepsilon}_\rho)$ (see the proof of \cref{lem:mollifying}(i)), and $\|f^n\|_{L_q ([s,t];\B^{\alpha}_\rho)}\le \|f\|_{L_q ([s,t];\B^{\alpha}_\rho)}$ for all $(s,t)$, since $P_{1/n}$ is a contraction on $\B^\alpha_\rho$. Applying \eqref{eq:U_rho_var} to $f^n$ with $\alpha$ itself and using Fatou's lemma as above, we obtain \eqref{eq:U_rho_var} for $f$.  
\end{proof}
\begin{remark}\label{rem:var_exp_bd}
 In view of \eqref{eq:cA_vmo} and \cref{prop:vmo_john_nirenberg}-\ref{jntwo}, we indeed have 
 \begin{align*}
   \E \left[ \exp \left( c \left( \frac{\norm{U^f}_{\ell\vr;[s, t]}}{\norm{f}_{L_q([s, t]; \mathcal{B}_{\rho}^{\alpha})} (t-s)^{1+\alpha H - \frac{1}{q} - \frac{Hd}{\rho}}}  \right)^{\frac{1}{-\alpha H + H d/\rho}}  \right)    \right] \leq C
 \end{align*}
 for some finite positive constants $c$ and $C$ depending only on $H, d, q, \rho, \alpha, \varepsilon$.
\end{remark}

\begin{remark}
\label{r:explanation}
We now see why it was crucial in the proof of \cref{lem:integral_along_fbm} to first apply the Kolmogorov theorem in $p$-variations and then the John--Nirenberg inequality, rather than the other way round. Indeed, bound \eqref{mainbound2} shows that for any fixed $s$ and $(u,t) \in \Delta_{(s,T)}$,
\begin{equation}\label{condAA}
(\E^s|U^f_t-U^f_u|^2)^{1/2}\le C \|f\|_{L_q ([u,t];\B^{\alpha}_\rho)}((t-s)^{\gamma}-(u-s)^{\gamma})^{1+\alpha H-\frac1q}=: w_s(u,t)^{1+\alpha H},
\end{equation}
where the control $w_s$ depends on $s$. Nevertheless, this is sufficient to show that the process $U^f$ has a.s. a finite $(\frac1{1+\alpha H}+\eps)$-variation over the interval $[s,T]$. From this, an application of the John--Nirenberg inequality (\cref{prop:vmo_john_nirenberg}) yields that all moments of this variation are finite.
On the other hand, if we applied the John--Nirenberg inequality first, then to verify condition \eqref{eq:vmo_gamma}, we would need to take $u=s$ in \eqref{condAA}, resulting in
\begin{equation*}
\E^s|U^f_t-U^f_s|\le C \|f\|_{L_q ([s,t];\B^{\alpha}_\rho)}(t-s)^{\gamma(1+\alpha H-\frac1q)}= \wt w(s,t)^{1+\alpha H-\frac{Hd}{\rho}},
\end{equation*}
for a control $\wt w$. This would give the bound
\begin{equation*}
	\|U^f_t-U^f_s\|_{L_m(\Omega)}\le C  \wt w(s,t)^{1+\alpha H-\frac{Hd}{\rho}},
\end{equation*}
and would ultimately yield the finiteness of $(\frac1{1+\alpha H-\frac{Hd}{\rho}}+\eps)$-variation of $U^f$, which is worse than the result of \cref{lem:integral_along_fbm}.
\end{remark}
Under \eqref{eq:parameters}, we can take $\alpha=-1$ in \cref{lem:integral_along_fbm} as shown in the following lemma.
\begin{lemma}\label{lem.minusoneOK}
    For all values of parameters $H,d,p,q$ satisfying \eqref{eq:parameters}, there exists $\eta_0\in(0,1)$ such that $(\alpha,\rho,d,q)=(-1-\eta_0,p,d,q)$ satisfies \eqref{con.for:U_var}.
  Since \eqref{con.for:U_var} is a lower bound on $\alpha H$, the same holds with $\eta_0$ replaced by any $\eta\in[0,\eta_0]$, and, for $\eta>0$, also for $(\alpha,\rho)=(-\eta,p)$.
  \end{lemma}
  \begin{proof}
    It suffices to show verify \eqref{con.for:U_var} with $(\alpha,\rho)=(-1,p)$ and choose $\eta_0\in(0,1)$ such that 
    \begin{align*}
      (-1-\eta_0)H>\max\left(-\frac12+Hd\bigl(\frac1\rho-\frac12\bigr)^+,-1+\frac1q+\frac{Hd}{\rho}\right).
    \end{align*}
    Under \eqref{eq:parameters}, the condition  \eqref{con.for:U_var} with $(\alpha,\rho)=(-1,p)$ is equivalent to   
    \begin{equation}\label{tmp.minusoneOK}
    H+\frac{Hd}{p}< \frac12+\frac{Hd}{p\vee2}.
    \end{equation}
  
  If $p\ge2$, there is nothing to prove: \eqref{tmp.minusoneOK} becomes $H<\frac12$ which is true by \eqref{eq:parameters_range}.
  
  If $p\le2$ and $q\le2$, then \eqref{eq:parameters_subcritical} gives
  \begin{equation*}
  H+\frac{Hd}p<1-\frac1q\le\frac12,
  \end{equation*}
  which implies \eqref{tmp.minusoneOK}.
  
  If $p\le2$, $q>2$, and $d\ge3$, then we have $H< Hd/2$ and 
  $\frac{Hd}{p}\le\frac12$ by \eqref{eq:parameters_girsanov}; thus \eqref{tmp.minusoneOK} also holds. 
  
  If $p\in[(1-H)^{-1},2]$, $q>2$, $d=1,2$, then
  $$
  H+Hd(\frac1p-\frac12)\le H+2 H(1-H-\frac12)=2H(1-H)<\frac12,
  $$
  since $H<\frac12$.  Therefore, condition \eqref{tmp.minusoneOK}  is satisfied in this case as well.
  
    Finally, if $p\in[1,(1-H)^{-1}]$, $q>2$, $d=1,2$, then \eqref{eq:parameters_cameron_martin} implies
  \begin{equation*}
  \frac{Hd}{p} < \frac{p^{-1} }{p^{-1} - 2^{-1}} \bigl( \frac{1}{2} - H \bigr). 
  \end{equation*}
  Therefore, 
  \begin{equation*}
  Hd(p^{-1} - 2^{-1}) < \frac{1}{2} - H, 
  \end{equation*}
  which is \eqref{tmp.minusoneOK}.
  \end{proof}

Combining \cref{lem:integral_along_fbm,lem.minusoneOK}, we obtain the following estimates, which are used repeatedly below.
For a function (or distribution) $f$ on $[0,T]\times\R^d$ and $x\in\R^d$, we write $f(\cdot+x)\defby f(\cdot,\cdot+x)$.
\begin{corollary}\label{cor:integral_estimates}
  Assume \eqref{eq:parameters} and let $G$ satisfy \eqref{con:local_nondeterminism}. Let $\eta_0$ be as in \cref{lem.minusoneOK}, $\eta\in[0,\eta_0]$, $m\ge1$ and $f\in L_q([0,T];\B^{-\eta}_p(\R^d))$.
  Then the following holds for every $(s,t)\in\Delta_{[0,T]}$.
  \begin{enumerate}[\rm{(}i\rm{)}]
    \item For any $\ell\in\big((1-H-\eta H-Hd(\frac1p-\frac12)^+)^{-1},2\big)$, there exists $C=C(\ell,m,\eta)$ such that
    \begin{equation}\label{eq:int_est_grad}
      \bigl\|\,\| U^{\nabla f}\|_{\ell\var;[s,t]}\bigr\|_{L_m(\Omega)|\F_s}\le C \|f\|_{L_q ([s,t];\B^{-\eta}_p)}(t-s)^{1-H-\frac1q-\frac{Hd}p-\eta H},
    \end{equation}
    where $U^{\nabla f}$ is understood componentwise.
        \item If $\eta>0$, then for any $\ell\in\big((1-\eta H-Hd(\frac1p-\frac12)^+)^{-1},2\big)$, there exists $C=C(\ell,m,\eta)$ such that
    \begin{equation}\label{eq:int_est_f}
      \bigl\|\,\| U^{f}\|_{\ell\var;[s,t]}\bigr\|_{L_m(\Omega)|\F_s}\le C \|f\|_{L_q ([s,t];\B^{-\eta}_p)}(t-s)^{1-\frac1q-\frac{Hd}p-\eta H}.
    \end{equation}
    \item There exists $C=C(m,\eta)$ such that for all $x,y\in\R^d$,
    \begin{equation}\label{eq:int_est_xy}
      \bigl\|U^{f(\cdot+x)}_t-U^{f(\cdot+x)}_s-U^{f(\cdot+y)}_t+U^{f(\cdot+y)}_s\bigr\|_{L_m(\Omega)|\F_s}\le C|x-y| \|f\|_{L_q ([s,t];\B^{-\eta}_p)}(t-s)^{1-H-\frac1q-\frac{Hd}p-\eta H}.
    \end{equation}
  \end{enumerate}
\end{corollary}
\begin{proof}
  (i) Since $\|\nabla g\|_{\B^{-1-\eta}_p}\le C\|g\|_{\B^{-\eta}_p}$ (see \cite[Chapter~2]{BCD11}), we have $\|\nabla f\|_{L_q([s,t];\B^{-1-\eta}_p)}\le C\|f\|_{L_q([s,t];\B^{-\eta}_p)}$. By \cref{lem.minusoneOK}, $(\alpha,\rho)=(-1-\eta,p)$ satisfies \eqref{con.for:U_var}, and \eqref{eq:int_est_grad} follows from \cref{lem:integral_along_fbm}.
  (ii) By \cref{lem.minusoneOK}, $(\alpha,\rho)=(-\eta,p)$ satisfies \eqref{con.for:U_var}, and \eqref{eq:int_est_f} follows from \cref{lem:integral_along_fbm}.
  (iii) Put $g\defby f(\cdot+x)-f(\cdot+y)$. By linearity of \eqref{eq:Uffn}, $U^{f(\cdot+x)}-U^{f(\cdot+y)}=U^g$. Writing $g=\int_0^1\nabla f(\cdot+y+\theta(x-y))\cdot(x-y)\,d\theta$ and using translation invariance of Besov norms, we get $\|g\|_{L_q([s,t];\B^{-1-\eta}_p)}\le C|x-y|\|f\|_{L_q([s,t];\B^{-\eta}_p)}$.
  Hence \eqref{eq:int_est_xy} follows from \cref{lem:integral_along_fbm} with $(\alpha,\rho)=(-1-\eta,p)$, since $|U^g_t-U^g_s|\le\|U^g\|_{\ell\var;[s,t]}$.
\end{proof}

When $f$ belongs to $L_q([0,T];L_\rho(\Rd))$, we have similar estimates. 
\begin{lemma}\label{P:lpq}
  Let $f\colon[0,T]\times\R^d\to\R$ be a measurable function in $L_q ([0,T];L_\rho(\Rd))$, where $q\in[1,\infty]$, $\rho\in[1,\infty]$.
  Assume that
  \begin{equation}\label{con.for:U_varlqlp}
    \frac1q+\frac{Hd}{\rho}<1.
  \end{equation}
  Then 
  \begin{align}\label{eq:EsU_lqlp}
    \E^s\abs[\big]{\int_s^t f(r,G_r)dr}\le C \norm{f}_{L_q([s,t];L_\rho(\Rd))} (t-s)^{1-\frac1q-\frac{Hd}{\rho}}.
  \end{align}
  Furthermore, $U^f=\int_0^\cdot f(r,G_r)dr$ a.s. and for any $\ell\in(\frac{1}{1-\frac{Hd}{\rho}},\infty)$, any $m\ge1$, there exists a constant $C=C(m,q,\rho,\ell)$ such that
  \begin{equation}\label{eq:U_varLqLp}
    \bigl\|\norm[\big]{\int_0^\cdot f(r,G_r)dr}_{\ell\var;[s,t]}\bigr\|_{L_m(\Omega)|\F_s}\le C \|f\|_{L_q ([s,t];L_\rho)}(t-s)^{1-\frac1q-\frac{Hd}\rho} \quad \forall s\le t\le T.
    \end{equation} 
\end{lemma}
\begin{proof}
    We have $\E^s |f(r,G_r)|=P_{\Sigma_{s,r,r}} |f|(r,\E^s G_r)$ for $r>s$. Thus, by Young's inequality, with $\rho'=\frac{\rho}{\rho-1}$, we have
  \begin{align*}
    \E^s|f(r,G_r)|
    \le \norm{p_{\Sigma_{s,r,r}}}_{L_{\rho'}(\Rd)} \norm{f(r,\cdot)}_{L_\rho(\Rd)}
    \le C (r-s)^{-\frac{Hd}{\rho}}\norm{f(r,\cdot)}_{L_\rho(\Rd)}.
  \end{align*}
  From \eqref{con:local_nondeterminism}, we have $\Sigma_{s,r,r}\ge \lambda^{-1}(r-s)^{2H}\I$, which implies that $\norm{p_{\Sigma_{s,r,r}}}_{L_{\rho'(\Rd)}} \les (r-s)^{-\frac{Hd}{\rho}}$. From here, since $\E^s|\int_s^tf(r,G_r)dr|\le\int_s^t\E^s|f(r,G_r)|dr$, we obtain \eqref{eq:EsU_lqlp} by H\"older's inequality in time, since $\frac{Hd}{\rho}\frac{q}{q-1}<1$ by \eqref{con.for:U_varlqlp}.

  A consequence of \eqref{eq:EsU_lqlp} is that $\int_0^T|f(r,G_r)|dr$ is finite a.s., which implies that $t\mapsto U^f_t$ is continuous a.s.. 
  Applying \cref{prop:vmo_john_nirenberg}, we have 
  \begin{align*}
    \norm[\big]{\int_s^t f(r,G_r)dr}_{L_m(\Omega)|\F_s}
    &\le C \norm{f}_{L_q([s,t];L_\rho(\Rd))} (t-s)^{1-\frac1q-\frac{Hd}{\rho}}
  \end{align*}
  for any $m\ge1$ and $s\le t\le T$. Applying \cref{lem:variation-kolmogorov}, we obtain \eqref{eq:U_varLqLp}. 

  To see that $U^f_t=\int_0^tf(r,G_r)dr$, we consider two cases. When $\rho<\infty$, we take $(f^n)$ as a sequence of spatial mollifiers of $f$. Such sequence belongs to $L_1B$ and satisfies $\lim_n f^n= f$ in $L_q([0,T];\B^{-\varepsilon}_\rho(\R^d))$ and in $L_q([0,T];L_\rho(\R^d))$. Then, by \eqref{eq:U_varLqLp}, we have $\int_0^\cdot f^n(r,G_r)dr$ converges uniformly to $\int_0^\cdot f(r,G_r)dr$ a.s. In view of \eqref{eq:Uffn}, this implies that $U^f_t=\int_0^tf(r,G_r)dr$. 
When $\rho=\infty$, we can find a function $\bar f$ in $L_qB$ such that $f(r,x)=\bar f(r,x)$ for a.e. $(r,x)$. Due to \eqref{eq:EsU_lqlp}, we have $\int_0^\cdot f(r,G_r)dr=\int_0^\cdot \bar f(r,G_r)dr$ a.s.. On the other hand, since $\bar f\in L_1B$, $\int_0^\cdot \bar f(r,G_r)dr=U^{\bar f}$ by definition. Due to \eqref{eq:U_rho_var}, $U^{\bar f}=U^{f}$.
  Thus, $U^f_t=\int_0^tf(r,G_r)dr$.
\end{proof}

\subsection{Exponential integrability }\label{subsec:girsanov}



Throughout this subsection, $H\in(0,\frac12)$ and $G$ is either the fractional Brownian motion $B^H$ or the Riemann--Liouville fractional Brownian motion $\tilde{B}^H$ from \eqref{eq:RL_def}, with Cameron--Martin space $\H^G$ as in \cref{subsec:girsanov_prelim}.
By \cref{lem:LND_examples}, the results of \cref{subsec:add_fc_fbm} apply to $G$, and we recall the embedding \eqref{GGbound} of Sobolev spaces into $\H^G$.

For a given $f$ in $L_q([0,T];L_p(\R^d))$, let $X$ be a strong solution to equation 
\begin{align}
  dX_t=f(t,X_t)\,dt+dG_t \quad\text{on}\quad [0,T], \label{eqn.XfG}
\end{align}
This means that $X$ is adapted to $\F^G$, $\int_0^T|f(r,X_r)|dr$ is finite and $ X_t=X_0+\int_0^t f(r,X_r)\,dr+G_t$ for all $t\in[0,T]$ a.s. We are going to estimate the exponential moments of $\norm {\int_0^\cdot f(r,X_r)dr}_{\H^G}$.

\begin{lemma}[Exponential bounds]\label{lem:drift_cameron_martin}
	Let $T>0$, $H\in(0,1)$, $p\in[1,\infty]$, $q\in[1,\infty]$, and 
	\[
	\frac{Hd}{p}+\frac1q<1.
	\]
	Let $f\in L_q([0,T],L_p(\R^d))$ and consider the additive functional $\int_0^\cdot f(r,G_r)\, dr$, which is identical to $U^f(G)$ by \cref{P:lpq}. 
	Then the following statements hold. 
	\begin{enumerate}[\rm{(}i\rm{)}]
		\item For any $\eps\in(0,1-\frac{Hd}p-\frac1q)$ there exist positive constants $C=C(\eps,T)$ and $\mu=\mu(\eps,T)$ such that 
		\begin{align}\label{whatwewantS1}
			&\E \exp\left( \mu\|f\|^{-\frac{p}{Hd}}_{L_qL_p}\|U^f\|_{\W^{1-\frac{Hd}p-\eps}_q([0,T])}^{\frac{p}{Hd}}\right)\le C \quad \text{if}\quad p<\infty
      \\\shortintertext{and}
      &\|U^f\|_{\W^{1-\frac{Hd}p-\eps}_q([0,T])}\le C\|f\|_{L_qL_p}\quad \text{a.s. if}\quad p=\infty.\label{whatwewantS11}
		\end{align}

		\item Let $u$ be such that 
		\begin{equation}\label{condr}
      1\le u\le \min(p,q), \tand \frac{1}{u}> \frac{Hd}{p}+\frac1q
		\end{equation}
		Then there exist positive constants $C=C(T)$ and $\mu=\mu(T)$ 
		such that
		\begin{align}\label{whatwewantS2}
			&\E \exp\left( \mu\|f\|^{-\frac{p}{Hd}}_{L_qL_p} \|U^f\|_{\W^{1}_u([0,T])}^{\frac{p}{Hd}}\right)\le C 
      \quad \text{if}\quad p<\infty
      \\\shortintertext{and}
      &\|U^f\|_{\W^{1}_u([0,T])}\le C \|f\|_{L_qL_p},\quad  \text{a.s. if}\quad p=\infty.\label{whatwewantS22}
		\end{align}

		\item\label{cmthree} Assume additionally, that conditions \eqref{eq:parameters_range}, \eqref{eq:parameters_subcritical}, \eqref{eq:parameters_cameron_martin} hold. Then
		there exist positive constants $C=C(T)$, $\mu=\mu(T)$ and $\eps=\varepsilon(p,q,H,d)$ such that
		\begin{align}\label{whatwewantS3}
			&\E \exp\left(\mu\|f\|^{-\frac{p}{Hd}}_{L_qL_p} \|U^f\|_{\W^{H+\frac12+\eps}_2([0,T])}^{\frac{p}{Hd}}\right)\le C 
      \quad \text{if}\quad p<\infty
		  \\\shortintertext{and}\label{whatwewantS4}
			&\|U^f\|_{\W^{H+\frac12+\eps}_2([0,T])}\le C \|f\|_{L_qL_p},\quad  \text{a.s. if}\quad p=\infty.
		\end{align}
	\end{enumerate}
\end{lemma}

\begin{proof} 
  We note that the relation $U^f_t=\int_0^t f(r,G_r)dr$ is used for notational convenience, and the presentation \eqref{eq:Uffn} is not necessary for the proof.

(i) By \cref{P:lpq}, there exists a constant $C>0$ such that for any $(s,t)\in\Delta_{[0,T]}$
	\begin{equation*}
		\E^s|U_t^f-U_s^f|\le C (t-s)^{1-\frac{Hd}{p}-\frac1q} \|f\|_{L_q([s,t];L_p(\R^d))}.
	\end{equation*}
		Recall from \eqref{controlmaps} that $(s,t)\mapsto \|f\|_{L_q([s,t];L_p(\R^d))}^q$ is a control. If $p<\infty$, \cref{lem:vmo_sobolev} with $\beta=1-\frac{Hd}p\in[\frac1q,1)$ yields \eqref{whatwewantS1}. If $p=\infty$, we have the pathwise bound $|U^f_t-U^f_s|\le\int_s^t\|f(r)\|_{L_\infty}dr\le\|f\|_{L_q([s,t];L_\infty)}(t-s)^{1-\frac1q}$, which directly gives \eqref{whatwewantS11}.
	
	(ii) By definition, we have
	\begin{align}\label{normw1}
    \|U^f\|_{\W^{1}_u([0,T])}^u
  =
	\int_0^T |f|^u(r,G_r)\,dr
  =U^{|f|^u}_T.
  \end{align}
	
  It is evident that the function $|f|^{u}\in L_{\frac{q}{u}}([0,T], L_{\frac{p}{u}}(\R^d))$ (note that $p/u\ge1$ and $q/u\ge1$ by \eqref{condr}).
  Additionally, $(Hd/p+1/q)u<1$ by \eqref{condr}
	Hence, \cref{P:lpq} implies that
	\begin{equation}\label{boundf}
		\E^s|U_t^{|f|^u}-U_s^{|f|^u}|\le C (t-s)^{1-\frac{Hdu}{p}-\frac{u}q} \| |f|^u\|_{L_{\frac{q}u}([s,t],L_{\frac{p}u}(\R^d))}=:w(s,t)^{1-\frac{Hdu}{p}},
	\end{equation}
	for some constant $C>0$. As previously, one sees that $w$, as defined above, is a  control.  Thus,  
	\cref{prop:vmo_john_nirenberg}\ref{jntwo} implies that for some constants $C>0$, $c>0$ we have
	\begin{align*}
		&\E\Bigl[\exp\Bigl( \frac{c|U_T^{|f|^u}|^{\frac{p}{Hdu}}}{w(0,T)^{\frac{p}{Hdu}(1-\frac{Hdu}p)}}\Bigr)\Bigr]\le C \quad \text{if}\quad p<\infty
    \\\shortintertext{and}
    &\|U_T^{|f|^u}\|_{L_{\frac{p}{Hdu}}(\Omega)}\le C w(0,T)^{1-\frac{Hdu}p} \quad \text{a.s. if}\quad p=\infty.
	\end{align*}
	Note that 
	\[
	w(0,T)^{1-\frac{Hdu}p} = C(T)\| |f|^u\|_{L_{\frac{q}u}([0,T],L_{\frac{p}u}(\R^d))}=
	C(T)\| f\|_{L_{q}([0,T],L_{p}(\R^d))}^u.
	\]
	Combining with \eqref{normw1}, we obtain \eqref{whatwewantS2} and \eqref{whatwewantS22}.
	
	(iii) We deduce this result from parts (i) and (ii) using embedding inequalities. We consider different cases depending on the relationship between the parameters.
	
	\textbf{Case~1}: $q\le2$. In this case, by the Sobolev embedding we have 
  \begin{align*}
    \|U^f\|_{\W^{1-\frac{Hd}p-\frac1q+\frac12-\eps}_2([0,T])}\le C \|U^f\|_{\W^{1-\frac{Hd}p-\eps}_q([0,T])}.
  \end{align*}
  Our assumption \eqref{eq:parameters_subcritical} implies 
	\[
	1-\frac{Hd}p-\frac1q+\frac12>H+\frac12.
	\]
	Therefore, by taking $\eps>0$ small enough, we get part (iii) from part (i).
	
	\textbf{Case~2}: $q>2$ and $p>\frac1{1-H}$. 
	If $u$ can be chosen so that $u\ge2$, then obviously 
  \begin{align*}
    T^{\frac{1}{u}-\frac{1}{2}  } \|U^f\|_{\W^{1}_2([0,T])}\le \|U^f\|_{\W^{1}_u([0,T])},
  \end{align*}
	which gives us \eqref{whatwewantS3} from (ii). Else, when $u<2$, then by Sobolev embedding, we have
  \begin{align*}
    \|U^f\|_{\W^{1-\frac1u+\frac12}_2([0,T])}\le C\|U^f\|_{\W^{1}_u([0,T])}.
  \end{align*}
  In the current case, it is possible to choose $u$ so that $1-\frac1u+\frac12\ge H+\frac12+\eps$ and at the same time \eqref{condr} holds.  This yields \eqref{whatwewantS3} and \eqref{whatwewantS4} in view of (ii).
	
	\textbf{Case~3}: $q>2$ and $p\le\frac1{1-H}$. This is the most challenging case and we  rely on the interpolation inequality \eqref{mainresGG}.
	
	First, we note that in this regime we can and choose $u=p$.  
    We apply \eqref{mainresGG} with $p_1=q$, $p_2=p$, $\gamma=\frac{Hd}{p}+\varepsilon'$ for a small $\varepsilon'>0$. We see that we must have
	$$
	\theta=\frac{\frac1p-\frac12}{\frac1p-\frac1q}
	$$
	and condition \eqref{GNcond} reads
	$$
	H+\frac{Hd}{p}\theta<\frac12,
	$$
		which, for $\varepsilon'=0$, is exactly the imposed assumption \eqref{eq:parameters_cameron_martin}. Since \eqref{eq:parameters_cameron_martin} is a strict inequality, \eqref{GNcond} also holds for $\varepsilon'>0$ small enough. By \cref{lem:GG} we get that there exist $\eps>0$ and a constant $C$ such that
  \begin{align*}
    \|U^f\|_{\W^{H+\frac12+\eps}_2([0,T])}\le C\|U^f\|_{\W^{1-\frac{Hd}p-\varepsilon'}_q([0,T])}+C \|U^f\|_{\W^{1}_p([0,T])}.
  \end{align*}
  Combining with parts (i) (with $\varepsilon'$ in place of $\eps$) and (ii) (with $u=p$), and using $(a+b)^{r}\le2^{r-1}(a^r+b^r)$ for $r=\frac{p}{Hd}\ge1$ together with the Cauchy--Schwarz inequality, we obtain \eqref{whatwewantS3} with a smaller constant $\mu$, and \eqref{whatwewantS4}.
\end{proof}

Our next step is to remove the restriction on $\lambda$ in \eqref{whatwewantS3}. This is done by
applying a truncation argument similar to the one used in \cite[Lemma C.3]{galeatigerencser22}. 
Recall the definition of $\Lambda_f$ introduced in \eqref{eq:frak_c}.

\begin{corollary}\label{c:alllambda}
Assume that the parameters $H,d,p,q$ satisfy  \eqref{eq:parameters} and that $p<\infty$ (if $q=\infty$, recall the convention of \cref{rem:Lambda_q_infty}).
Let $T>0$, $\lambda>0$, $f\in L_q([0,T],L_p(\R^d))$. 
Then there exist finite positive  constants $\eps=\eps(p,q,H,d)$, $\eta =\eta (T,\lambda)$, $c=c(T,\lambda)$ and $C=C(T,\lambda)$ 
such that
\begin{equation}\label{whatwewantS51}
	\E \exp\bigl( \lambda \|\int_0^\cdot f(r,G_r)dr\|_{\W^{H+\frac12+\eps}_2([0,T])}^{2}\bigr)\le
  C e^{c\Lambda_f^2(\eta )}.
\end{equation}
If, in addition, $p>2Hd$, then there exists a finite constant $C'=C'(T,\lambda,\|f\|_{L_qL_p})$, increasing in $\|f\|_{L_qL_p}$, such that
\begin{equation}\label{whatwewantS52}
	\E \exp\bigl( \lambda \|\int_0^\cdot f(r,G_r)dr\|_{\W^{H+\frac12+\eps}_2([0,T])}^{2}\bigr)\le C'.
\end{equation}

\end{corollary}
\begin{proof}
  As in the previous result, we can write $U^f=\int_0^\cdot f(r,G_r)dr$ for notational convenience. 
Fix $\lambda>0$ and take $\eps$ as in \eqref{whatwewantS3}. If necessary, we decrease $\eps$ so  that $H+\frac12+\eps\le1$.

We first prove \eqref{whatwewantS52}. Let $p>2Hd$, put $r\defby\frac{p}{Hd}>2$ and assume $f\neq0$. Let $\mu$ be as in \eqref{whatwewantS3} and $\lambda_*\defby\mu\|f\|^{-r}_{L_qL_p}$. By Young's inequality with exponents $\frac r2$ and $\frac{r}{r-2}$,
\begin{equation*}
  \lambda x^2\le \lambda_* x^{r}+ \lambda^{\frac{r}{r-2}}\lambda_*^{-\frac{2}{r-2}},\quad x\ge0.
\end{equation*}
Hence \eqref{whatwewantS3} implies
\begin{equation*}
\E \exp\bigl( \lambda \|U^f\|_{\W^{H+\frac12+\eps}_2([0,T])}^{2}\bigr)\le
C \exp\bigl( \lambda^{\frac{r}{r-2}}\mu^{-\frac{2}{r-2}}\|f\|_{L_qL_p}^{\frac{2r}{r-2}}\bigr),
\end{equation*}
which is \eqref{whatwewantS52}.

We now prove \eqref{whatwewantS51}, which holds for every $p\in[2Hd,\infty)$.
We decompose $f=f_\eta +(f-f_\eta )$ where $f_\eta :=f\I_{(|f|>\Lambda_f(\eta ))}$ for some $\eta>0$ chosen later.  By definition of $\Lambda_f$  in \eqref{eq:frak_c} (see also \cref{rem:Lambda_q_infty}), we have
\begin{equation}\label{nprop}
\|f_\eta \|_{L_q([0,T];L_p(\R^d))}\le \eta .
\end{equation}
By definition of $U^f$, we have $U^f=U^{f_\eta }+U^{f-f_\eta }$. Note that by \eqref{eq:trivialW}, we have
\begin{equation*}
\|U^{f-f_\eta }\|_{\W^{H+\frac12+\eps}_2([0,T])}\le C 
\|U^{f-f_\eta }\|_{\W^{1}_2([0,T])}\le C \sqrt{T} \|f-f_\eta \|_{L_\infty}=C\sqrt{T}\Lambda_f(\eta ).
\end{equation*}
Therefore,
\begin{align}\label{biglambda1}
&\E\exp\bigl( \lambda \|U^{f}\|_{\W^{H+\frac12+\eps}_2([0,T])}^{2}\bigr)\nn\\
&\quad\le
\E\exp\bigl( 4\lambda \|U^{f_\eta }\|_{\W^{H+\frac12+\eps}_2([0,T])}^{2}+
 4\lambda \|U^{f-f_\eta }\|_{\W^{H+\frac12+\eps}_2([0,T])}^{2}\bigr)\nn\\
 &\quad\le \exp(4\lambda C^2 T \Lambda_f(\eta )^{2})
\E\exp\bigl( 4\lambda \|U^{f_\eta }\|_{\W^{H+\frac12+\eps}_2([0,T])}^{2}\bigr).
\end{align}
We take $\eta >0$ such that $4 \lambda \eta ^{2} = \mu$,
where $\mu$ is as in \eqref{whatwewantS3}. Then, recalling \eqref{nprop}, we get 
$$
4 \lambda \|f_\eta \|^{2}_{L_q([0,T];L_p(\R^d))} \le \mu.
$$
It follows that 
\begin{align*}
  \E\exp\bigl( 4\lambda \|U^{f_\eta }\|_{\W^{H+\frac12+\eps}_2([0,T])}^{2}\bigr)&\le
  \E\exp\bigl( \mu \|f_\eta \|^{-2}_{L_q([0,T];L_p(\R^d))} \|U^{f_\eta }\|_{\W^{H+\frac12+\eps}_2([0,T])}^{2}\bigr).
\end{align*}
Applying the elementary inequality $x^2\le x^{\frac{p}{Hd} }+1$, the expectation above is bounded by
\begin{align*}
  e^\mu\E\exp\bigl( \mu \|f_\eta \|^{-\frac{p}{Hd}}_{L_q([0,T];L_p(\R^d))} \|U^{f_\eta }\|_{\W^{H+\frac12+\eps}_2([0,T])}^{\frac{p}{Hd}}\bigr).
\end{align*}
which is finite by \eqref{whatwewantS3}.
\end{proof}

Now we are ready to prove the main result of this subsection: quantitative bounds for the drift of solutions of \eqref{eqn.XfG}. 
\begin{lemma}[Exponential bound for the integral along the solution]\label{lem:girsanov_sol}
Assume that \eqref{eq:parameters} holds. Let $f\in L_q([0,T];L_p(\R^d))$, $x\in\R^d$.   Let $X$ be a strong solution to \eqref{eqn.XfG} with $X_0=x$.
Set 
 \begin{equation*}
 \psi_t:= \int_0^t f(r, X_r) \,d r,\quad t\in[0,T].
\end{equation*}
Assume that
\begin{equation}\label{girsol}
\E\exp\bigl(\lambda\|\psi\|^2_{\H^G}\bigr)<\infty
\quad\text{for any}\quad \lambda>0.
\end{equation}
Then 
there exists a constant  $\eps=\eps(p,q,H,d)>0$, such that for any $\lambda>0$, there exist finite positive  constants $\eta =\eta (T,\lambda)$, $c=c(T,\lambda)$, $C=C(T,\lambda)$ and $\bar C=\bar C(T)$ such that
\begin{align}
    &\E\exp\bigl( \lambda \|\psi\|_{\W^{H+\frac12+\eps}_2([0,T]) }^2 \bigr)\le C e^{c\Lambda_f^2(\eta )},\twhen p<\infty,\label{novikovbound1}\\
  &\E\exp\bigl( \lambda \|\psi\|_{\W^{H+\frac12+\eps}_2([0,T]) }^2 \bigr)\le C'(T,\lambda,\|f\|_{L_qL_p}),\twhen 2Hd<p<\infty,\label{novikovbound1b}\\
  &\|\psi\|_{\W^{H+\frac12+\eps}_2([0,T]) }\le \bar C \|f\|_{L_qL_p},\quad a.s.,\twhen p=\infty.\label{novikovbound2}
\end{align}
\end{lemma}

\begin{remark}\label{r:boundedisok}
The a priori assumption \eqref{girsol} is only used to justify the change of measure in the proof below; the point of \cref{lem:girsanov_sol} is that the resulting bounds \eqref{novikovbound1}--\eqref{novikovbound2} do not depend on the quantity in \eqref{girsol}.
The assumption is always satisfied if $f$ is bounded: by \eqref{GGbound} with $\eps=\frac12-H$,
\begin{equation*}
  \|\psi\|_{\H^G}\le C \norm{\psi}_{\W^{1}_2([0,T])} = C\Big(\int_0^T|f(r,X_r)|^2\,dr\Big)^{\frac12} \leq C T^{\frac12} \sup_{r,x}|f(r,x)|,
\end{equation*}
so that \eqref{girsol} holds for any $\lambda>0$. This is how \cref{lem:girsanov_sol} is applied to smooth approximations of the drift.
\end{remark}

\begin{proof}
Let $\eps>0$ be as in \cref{c:alllambda} and fix $\lambda>0$. 
Observe that thanks to assumption \eqref{girsol}, all conditions of \cref{lem:girsanov} are satisfied. Therefore, there exists a probability measure $\wt \P$ equivalent to $\P$ such that the process $X-x=G+\psi$ has the law of $G$ under this measure.
When $p=\infty$, we easily deduce \eqref{novikovbound2} from \eqref{whatwewantS4}. Henceforth, we assume that $p<\infty$. 
Recalling the definition of the process $\xi=\xi^G$ in \eqref{eq:xi_def}, and writing $\|\cdot\|_{\W}\defby\|\cdot\|_{\W^{H+\frac12+\eps}_2([0,T])}$ and $C_0$ for the constant in \eqref{GGbound}, we deduce that
\begin{align}
\E_{\P}\bigl[e^{\lambda\|\psi\|^2_{\W}}\bigr]&=\E_{\wt \P}\bigl[e^{\lambda\|\psi\|^2_{\W}}\xi_T(\psi)^{-1}]\nn\\
&\le
\bigl(\E_{\wt \P} e^{2\lambda\|\psi\|^2_{\W}}\bigr)^{1/2}
\bigl(\E_{\wt \P} [\xi_T(\psi)^{-2}]\bigr)^{1/2}\nn\\
&\le \bigl(\E_{\wt \P} e^{2\lambda\|\psi\|^2_{\W}}\bigr)^{1/2}\bigl(\E_{\wt \P} e^{6\|\psi\|^2_{\H^G}}\bigr)^{1/4}
\le \bigl(\E_{\wt \P} e^{(2\lambda+6C_0^2)\|\psi\|^2_{\W}}\bigr)^{3/4} \label{girsol1}\\
&=\bigl(\E_{ \P} e^{(2\lambda+6C_0^2)\left\|\int_0^\cdot f(r,G_r+x)\,dr\right\|^2_{\W}}\bigr)^{3/4}\label{almosthere}
\end{align}
where \eqref{girsol1} follows from  \eqref{eq:gir_density_mom2} (with $\lambda=-2$), \eqref{GGbound} and H\"older's inequality, and in \eqref{almosthere}
we used that $\Law_{\wt \P}(X)=\Law_{\P}(x+G)$.
Applying \cref{c:alllambda} to estimate \eqref{almosthere}, noting that
$\|f(\cdot,x+\cdot)\|_{L_qL_p}=\|f\|_{L_qL_p}$ and $\Lambda_{f(\cdot,x+\cdot)}=\Lambda_f$, we  obtain the desired bounds \eqref{novikovbound1} and \eqref{novikovbound1b}.
\end{proof}
\begin{remark}\label{rem:int_f_sol}
    In the setting of \cref{lem:girsanov_sol}, changing the measure as in its proof and applying \cref{P:lpq} under $\wt\P$, we have for any $m\ge1$
  \begin{equation*}
    \norm{\psi_t - \psi_s}_{L_m(\Omega)} \leq C \norm{f}_{L_qL_p} 
    (t-s)^{1 - \frac{Hd}{p} - \frac{1}{q}},
  \end{equation*}
  where the constant $C$ is chosen as above.
\end{remark}

\subsection{Linearization identity}\label{subsec:linearization}
As in \cref{subsec:girsanov}, $H\in(0,\frac12)$ and $G$ is either $B^H$ or $\tilde B^H$. Let $\psi^1,\psi^2$ be two adapted processes and let $f$ be a measurable function. When $f$ is regular enough so that the mean value theorem is applicable, we have the following linearization identity
\begin{multline*}
  \int_{0}^t f(r,G_r+\psi^1_r)\,dr-\int_{0}^t f(r,G_r+\psi^2_r)\,dr
  \\=\int_0^1\int_0^t \nabla f(r,G_r+\theta\psi^1_r+(1-\theta)\psi^2_r) \cdot (\psi^1_r-\psi^2_r)\,dr\,d\theta .
\end{multline*}
In the next result, we show that this identity can be extended to integrable function $f$ under some suitable conditions. 
\begin{proposition}\label{prop:linearization}
  Let $f\in L_q([0,T],L_p(\Rd))$ for $p,q$ satisfying condition \eqref{eq:parameters}. Let $\psi^1,\psi^2$ be two adapted processes
  such that for each $\lambda>0$
\begin{align*}
  M(\lambda):=\sup_{i=1,2} \E e^{\lambda \norm{\psi^i}_{\H^G}^2} <\infty.
\end{align*}
Then there exists a continuous adapted process $V=V^{f,G+\psi^1,G+\psi^2}$ with values in $\R^{d\times d}$ such that the following properties hold.
\begin{enumerate}[\rm{(}i\rm{)}]
    \item For any $\ell\in \left(\frac1{1- H-Hd(\frac{1}{p}-\frac{1}{2})^+},\ 2\right)$, $V$ has finite $\ell$-variation.
  \item With probability one, 
  \begin{align}\label{eq:linearization}
    \int_{0}^t f(r,G_r+\psi^1_r)\,dr-\int_{0}^t f(r,G_r+\psi^2_r)\,dr
    =\int_0^t (\psi^1_r-\psi^2_r)\cdot d V_r, \quad \forall t\in[0,T],
  \end{align} 
  where the integral with respect to $dV$ is defined in the sense of Young integral.
    \item For $\eta_0$ as in \cref{lem.minusoneOK},  for any $m\ge1$ and any $\eta\in[0,\eta_0]$ such that $\ell>(1-H-\eta H-Hd(\frac1p-\frac12)^+)^{-1}$, there exists a constant $C=C(m,\ell,\eta)M(12)^{1/{(2m)}}$ such that for any $(s,t)\in \Delta_{[0,T]}$,
  \begin{equation}\label{eq:V_ellvar}
    \bigl\|\| V\|_{\ell\var;[s,t]}\bigr\|_{L_m(\Omega)}\le C \|f\|_{L_q ([s,t];\B^{-\eta}_p(\Rd))}(t-s)^{1- H-\frac1q-\frac{Hd}p -\eta H}.
    \end{equation}
\end{enumerate}
\end{proposition}
\begin{remark}\label{rem:linearization_shift}
  \Cref{prop:linearization} remains valid, with the same constants, if $\psi^i$ is replaced by $y^i+\psi^i$ for deterministic $y^i\in\Rd$, $i=1,2$ (the exponential moment condition being imposed on $\psi^i$ only). Indeed, in the proof, the Girsanov transform turns $G+\theta(y^1+\psi^1)+(1-\theta)(y^2+\psi^2)$ into $G+\theta y^1+(1-\theta)y^2$, and \cref{lem:integral_along_fbm} is applied to the translated function $\nabla f^n(\cdot,\cdot+\theta y^1+(1-\theta)y^2)$, whose Besov norms coincide with those of $\nabla f^n$. In particular, for strong solutions $X^i=x^i+B^H+\int_0^\cdot b^i(r,X^i_r)\,dr$, $i=1,2$, we write $V^{b^1,X^1,X^2}$ for the process obtained with $f=b^1$, $y^i=x^i$ and $\psi^i=\int_0^\cdot b^i(r,X^i_r)\,dr$; this is the process used in \eqref{eq:def_V_and_R}.
\end{remark}
\begin{proof}
  Because of part (i) and the embedding $\H^G\hookrightarrow C^{1/\delta\var}$ for any $\delta\in(1/2,1/2+H)$ (for $G=B^H$ see \cite[Remark 15.11]{friz_victoir_2010}; for $G=\tilde B^H$ this follows from \eqref{GGbound}, \cite[Theorem~3.1]{Gorenflo:2015vr} and \cite[Corollary~A.3]{friz_victoir_2010}), we see that the Young integral with respect to $dV$ in \eqref{eq:linearization} is well-defined.

  Let $(f^n)$ be the mollifying sequence of $f$ as in \cref{lem:mollifying}. For each $n$, we can apply the mean value theorem to obtain the linearization identity
  \begin{align*}
    \int_{0}^t f^n(r,G_r+\psi^1_r)\,dr-\int_{0}^t f^n(r,G_r+\psi^2_r)\,dr
    =\int_0^t (\psi^1_r-\psi^2_r)\cdot d V^n_r,
  \end{align*}
  where 
  \begin{align*}
    V^n_t=\int_0^1\int_0^t \nabla f^n(r,G_r+\theta\psi^1_r+(1-\theta)\psi^2_r) \,dr\,d\theta.
  \end{align*}
  For each $i=1,2$, the law of $G+\psi^i$ is absolutely continuous with respect to the law of $G$ by \cref{lem:girsanov}. It follows from \cref{P:lpq} that almost surely
  \begin{align*}
    \lim_n\int_{0}^\cdot  f^n(r,G_r+\psi^i_r)\,dr
    =\int_{0}^\cdot  f(r,G_r+\psi^i_r)\,dr,\quad \text{uniformly on }[0,T].
  \end{align*}

  In view of the continuity of Young integrals in variation norms, to obtain \eqref{eq:linearization}, it suffices to show that with probability one, $V^n$ converges to a continuous adapted process $V$ in $\ell$-variation for an $\ell$ as in the statement.
  Indeed, we have
\begin{equation}\label{genboundstep0}
\| \, \|V^n\|_{\ell\var;[s,t]} \|_{L_m(\Omega)}\le \int_0^1 \bigl\|\,\bigl\| \int_0^{\cdot} \nabla f^n
		(r,G_r+\theta \psi^1_r+ (1 - \theta) \psi^2_r)\,dr\bigr\|
_{\ell\var;[s,t]}\bigr\|_{L_m(\Omega)}\, d \theta.
\end{equation}
For any $\lambda>0$ and  $\theta\in[0,1]$, we have
\begin{equation*}
\E e^{\lambda \|\theta \psi^{1} + (1-\theta) \psi^{2}\|_{\H^G}^2}\le \E e^{2\lambda( \|\psi^1\|_{\H^G}^2+\|\psi^2\|_{\H^G}^2)}
\end{equation*}
which is finite by assumption.
Hence, Girsanov's theorem (\cref{lem:girsanov}) is applicable: for any $\theta\in[0,1]$ there exists an equivalent probability measure $\wt\P_\theta$ such that the process 
$
G+\theta \psi^{1} + (1-\theta) \psi^{2}
$
has under $\wt\P_\theta$ the same law as $G$ under $\P$.
We put
\begin{align*}
  F_{s,t}=\bigl\| \int_0^{\cdot} \nabla f^n
(r,G_r+\theta \psi^1_r+ (1 - \theta) \psi^2_r)\,dr\bigr\|
_{\ell\var;[s,t]}^m.
\end{align*} 
We have, by changes of measures and Cauchy--Schwarz inequality, 
\begin{align}\label{genboundstep1}
  \E_{\P} F_{s,t}=\E_{\wt \P_\theta} F_{s,t} \xi_{T}^{-1}
  \le 
  \left(\E_{\wt \P_\theta}F_{s,t}^2 \right)^{\frac{1}{2}} \left(\E_{\wt \P_\theta} \xi_{T}^{-2}\right)^{\frac{1}{2}}.
\end{align}
We note that by \cref{lem:girsanov}(i), under $\wt \P_\theta$ the process $G+\theta \psi^{1} + (1-\theta) \psi^{2}$ is of the form $\int_0^\cdot K^G(\cdot,r)\,d\wt W_r$ for an $(\F_t)$-Brownian motion $\wt W$; hence, by \cref{lem:LND_examples}, the results of \cref{subsec:add_fc_fbm} apply to it under $\wt\P_\theta$ with respect to $(\F_t)$. Hence, applying \cref{lem:integral_along_fbm} to $\nabla f^n\in L_q([0,T];\B^{-1-\eta}_p)$ for $\eta\in[0,\eta_0]$, which is applicable under \eqref{eq:parameters} thanks to \cref{lem.minusoneOK}, noting that $\|\nabla f^n\|_{L_q([s,t];\B^{-1-\eta}_{p})}\le \| f^n\|_{L_q([s,t];\B^{-\eta}_{p})}$, we can find a constant $C=C(\ell,m)>0$ so that
\begin{equation}\label{eq:nabla_b_1_B_theta}
  \E_{\wt \P_\theta}F_{s,t}^2   
    \le C (t-s)^{2m(1-H-\frac{Hd}{p} - \frac{1}{q}-\eta H)}
  \| f^n\|_{L_q([s,t];\B^{-\eta}_{p})}^{2m}.
\end{equation}
Next, we treat the second factor on the right-hand side of \eqref{genboundstep1}. We apply \eqref{eq:gir_density_mom} 
to get
\begin{align}
  \E_{\wt \P_\theta} \xi_{T}^{-2}
  =\E_{\P}\xi_{T}^{-1}
  &\le \E_{\P}e^{3\|\theta \psi^{1} + (1-\theta) \psi^{2}\|^2_{\H^G}}
  \le\bigl(\E_{ \P}e^{12\| \psi^{1}\|^2_{\H^G}}\bigr)^{\frac12} \bigl(\E_{\P}e^{12\|\psi^{2}\|^2_{\H^G}}\bigr)^{\frac12}.\label{factor2new}
\end{align}
Using these estimates in \eqref{genboundstep1}, we obtain that 
\begin{equation}\label{tmp.vnellvar}
  \bigl\|\| V^n\|_{\ell\var;[s,t]}\bigr\|_{L_m(\Omega)}\les\|f^n\|_{L_q ([s,t];\B^{-\eta}_p)}(t-s)^{1- H-\frac1q-\frac{Hd}p-\eta H}.
\end{equation}
By linearity, we have
\begin{align*}
  \bigl\|\| V^n-V^k\|_{\ell\var;[0,T]}\bigr\|_{L_m(\Omega)}\les\|f^n-f^k\|_{L_q ([0,T];\B^{-\eta}_p)}.
\end{align*}
Since $\eta$ can be chosen to be positive, the previous estimate and \cref{lem:mollifying} imply that $V^n$ converges in $\ell$-variation in $L_m(\Omega)$ to a continuous adapted process $V$. Thus, we have shown \eqref{eq:linearization}. 

Finally, to obtain \eqref{eq:V_ellvar}, it suffices to pass through the limit in \eqref{tmp.vnellvar}.
  The proof is complete.
\end{proof}

\section{Strong existence}\label{sec:strong.existence}
In this section, we prove the existence of a strong solution to \eqref{eq:main-sde} through approximating equations with regular drifts. Along the way, we obtain pathwise stability for solutions satisfying an exponential moment condition (\cref{l:pathwise.stabitliy}) and path-by-path regularization estimates (\cref{subsec:refined_pathwise_estimates}), which are also used in \cref{sec:flow_and_path_by_path_uniqueness}; stability in moment norms is the subject of \cref{sec:stability}. 
\begin{proposition}[Strong existence]\label{l:strongsol}
  Assume that \eqref{eq:parameters} holds.  Let $T>0$, $x\in\R^d$, $b\in L_q([0,T];L_p(\R^d))$.
  There exists a strong solution $X$ to \Gref{x;b} such that for some $\eps>0$,
  for any $\lambda>0$
    \begin{align}\label{nboundsolnew}
   \E \exp ( \lambda \|X-B^H-x\|_{\W^{H+\frac12+\eps}_2([0,T])}^2 )<\infty .
  \end{align}
  More precisely, for every $\lambda>0$ there exist $\eta=\eta(T,\lambda)$ and finite constants $c=c(T,\lambda)$, $C=C(T,\lambda)$ and $C'=C'(T,\lambda,\|b\|_{L_qL_p})$, the latter increasing in $\|b\|_{L_qL_p}$, such that
  \begin{align}\label{nboundsolquant}
   \E \exp ( \lambda \|X-B^H-x\|_{\W^{H+\frac12+\eps}_2([0,T])}^2 )\le C'\I_{\{p>2Hd\}}+Ce^{c\Lambda_b^2(\eta)}\I_{\{p=2Hd\}}.
  \end{align}
  \end{proposition}
  To show \cref{l:strongsol}, we rely on a pathwise stability result for strong solutions satisfying certain exponential moment conditions, which is weaker than the stated pathwise stability in \cref{sec:main_results}. From here, it is relatively straightforward to conclude that the sequence $(X^n)$ converges to a limiting process $X$. 
  However, because $b$ is not continuous, showing that $X$ is indeed a solution to \Gref{x;b} is intricate and requires additional analysis.  This is where the refined pathwise estimates, inspired by \cite{anzeletti2023pathbypath}, are needed. These intermediate results are presented in \cref{sub.pathwise.stability,subsec:refined_pathwise_estimates} respectively. The proof of \cref{l:strongsol} is given in \cref{sub:existence_strong_sol}.

\subsection{Pathwise stability}\label{sub.pathwise.stability}

\begin{proposition}[Pathwise stability]\label{l:pathwise.stabitliy}
  Assume that \eqref{eq:parameters} holds. 
  Let $x^1,x^2\in\R^d$, and $b^1,b^2\in L_q([0,T];L_p(\R^d))$. For $i=1,2$, let $X^i$ be a strong solution to \Gref{x^i;b^i} such that 
  \begin{align}\label{con.expXi}
    M(\lambda):=\sup_{i=1,2}\E e^{\lambda \norm{\int_0^\cdot b^i(r,X^i_r)\,dr}_{\H}^2}<\infty,\quad \forall \lambda>0.
  \end{align}
  
  Then for any $\ell\in(\frac{1}{1-H-\frac{Hd}{p}+\frac{Hd}{p\vee2}},2)$ there exist positive 	constants $c=c(\ell)$, $C = C(\ell)$ such that a.s.
    \begin{equation}\label{eq:pathwise_stability}
      \sup_{t \in [0, T]} |X^1_t - X^2_t|
      \leq 4 e^{c \|V\|_{\ell\var;[0,T]}^{\ell} }   (\abs{x^1 - x^2}+\norm{R}_{\ell\var;[0,T]}),
    \end{equation}
    and
    \begin{align}
      \bra{X^1 - X^2}_{\ell\var;[s,t]}\le C(\sup_{r \in [0, T]} |X^1_r - X^2_r|  \|V\|_{\ell\var;[s,t]} +  \|R\|_{\ell\var;[s,t]}), \label{eq:stability_ellvar}
    \end{align}
  where we have denoted, for $t \in [0, T]$, 
  \begin{equation}\label{eq:def_V_and_R}
    V_t=V^{b^1,X^1,X^2}_t
    \tand  
    R_t:=\int_0^t (b^1-b^2)(r,X_r^2)\,dr,
  \end{equation}	
  for $V^{b^1,X^1,X^2}$ as in \cref{prop:linearization,rem:linearization_shift}.
  
    Furthermore, for $\eta_0$ as in \cref{lem.minusoneOK}, for any $\eta\in[0,\eta_0]$ such that $\ell>(1-H-\eta H-Hd(\frac1p-\frac12)^+)^{-1}$ and any $m\ge1$, there exists a constant $C_0=C_0(\eta,m)M(12)^{1/(2m)}$ such that for each $(s,t)\in \Delta_{[0,T]}$, we have
  \begin{align}\label{Vbound11}
    &	\| \, \|V\|_{\ell\var;[s,t]} \|_{L_m(\Omega)}\le  C_0\| b^1\|_{L_q([s,t];\B^{-\eta}_p)}
    (t-s)^{1-H-\frac{Hd}{p} - \frac{1}{q}-\eta H} ,\\
    &	\| \, \|R\|_{\ell\var;[s,t]} \|_{L_m(\Omega)}\le C_0\|b^1-b^2\|_{L_q([s,t];\B^{-1-\eta}_p)}
    (t-s)^{1-H-\frac{Hd}{p} - \frac{1}{q}-\eta H}\label{Vbound22}.
  \end{align}
  \end{proposition}
  Note that the quantity  $1-H-\frac{Hd}{p}+\frac{Hd}{p\vee2}$ is obviously less than $1$. On the other hand, by \cref{lem.minusoneOK}, this quantity is bigger than $1/2$. Therefore, the interval of $\ell$ satisfying the assumptions of the lemma is not empty and any such $\ell$ is greater than $1$.
 
\begin{proof}[Proof of \cref{l:pathwise.stabitliy}]
Let $Y:=X^1 - X^2$. 
Then for $(s,t)\in\Delta_{[0,T]}$, applying \cref{prop:linearization}, we have 
\begin{align}
	Y_t-Y_s&= \int_s^t (b^1(r,X^1_r) - b^1(r,X^2_r))\,d r + \int_s^t (b^1 - b^2)(r, X^2_r)\,dr\nn \notag\\
	&=\int_s^t Y_r \cdot d V_r+R_t-R_s.  \label{eq:Y_linearization}
\end{align}
It is obvious that $Y$ is continuous and has finite $\ell$-variation on $[0,T]$. An application of \cref{lem:linearYDE} yields \eqref{eq:pathwise_stability} and \eqref{eq:stability_ellvar}.

The estimate \eqref{Vbound11} follows from \cref{prop:linearization}. The estimate for $R$ in \eqref{Vbound22} follows from \cref{lem:integral_along_fbm} applied to $b^1-b^2$ with $(\alpha,\rho)=(-1-\eta,p)$, which is admissible by \cref{lem.minusoneOK}, and Girsanov's theorem, by an argument similar to the one used in \cref{prop:linearization}; hence we omit the details.
\end{proof}

\subsection{Refined pathwise estimates}\label{subsec:refined_pathwise_estimates}
\begin{lemma}[Path-by-path estimate on averaging]\label{lem.regularizing}
  Assume \eqref{eq:parameters}.
        Let $\eta_0$ be the constant in \cref{lem.minusoneOK}, $\eta\in(0,\eta_0]$, $K\subset\Rd$  be a closed ball centred at $0$, $\gamma \in (0,1-H-\frac{Hd}{p}-\frac{1}{q}-\eta H)$ and $\varepsilon \in (0,1)$. There exists a map $\Xi_{T,K}=\Xi_{T,K,\gamma,\varepsilon,\eta}\colon L_q\B^{-\eta}_p\to L_1(\Omega)$ such that
  \begin{enumerate}[\rm{(}i\rm{)}]
      \item \label{(i)}$\Xi_{T,K}$ satisfies the triangle inequality, i.e. $\Xi_{T,K}(f+g)\le \Xi_{T,K}(f)+\Xi_{T,K}(g)$;
                \item\label{it.xif} for every $f\in L_q\B^{-\eta}_p$,
      \begin{align}
          \E[\Xi_{T,K}(f)]\leq C(T, K, \gamma, \varepsilon, \eta) \|f\|_{L_q\B^{-\eta}_p};
          \label{est.Xif}
      \end{align}
      \item \label{(ii)}for every  $f\in L_qL_p$, there exists an event $\Omega_{f,T,K}$ of full measure such that for every $\omega\in \Omega_{f,T,K}$, every  $s\le t\le T$ and every $x,y\in K$, $r\mapsto f(r,B^H_r(\omega)+x)$ is integrable on $[0,T]$ and
      \begin{align}
        \Big|\int_s^tf(r,B^H_r(\omega)+x)dr \Big|
         &\le \Xi_{T,K}(f)(\omega)(t-s)^{\gamma},
          \label{est.asbounded}\\
          \Big|\int_s^t [f(r,B_r^H(\omega)+x)-f(r,B_r^H(\omega)+y)]\dd r\Big|
          &\le \Xi_{T,K}(f)(\omega)|x-y|^{1- \varepsilon}(t-s)^{\gamma}.\label{est.asLip}
      \end{align}
  \end{enumerate}
\end{lemma}
\begin{proof}
    The argument is similar to that of \cite[Lemma 4.1]{anzeletti2023pathbypath} with some adjustments; here $G=B^H$. Put $\gamma_*\defby1-H-\frac{Hd}{p}-\frac1q-\eta H$, so that $\gamma<\gamma_*$.
  Fix $m\ge2$ so large that $\frac2m<\gamma_*-\gamma$ and $\frac{2d-1}{m}<\varepsilon$, and fix
  \begin{equation*}
    \gamma_1\in\Big(\gamma+\frac2m,\ \gamma_*+\frac1m\Big),\qquad \gamma_2\in\Big(1-\varepsilon+\frac{2d}m,\ 1+\frac1m\Big).
  \end{equation*}
  For each $f\in L_q\B^{-\eta}_p$, define
  \begin{multline*}
    \Xi_{T,K}(f)=C_{T,K} \Big(\iint_{[0,T]^2} \Big|\frac{U^f_t-U^f_s}{|t-s|^{\gamma_1}}\Big|^m dsdt
    \\+ \iint_{[0,T]^2\times K^2} \Big|\frac{U^{f(\cdot+x)}_t-U^{f(\cdot+x)}_s-U^{f(\cdot+y)}_t+U^{f(\cdot+y)}_s}{|t-s|^{\gamma_1}|x-y|^{\gamma_2}}\Big|^m dsdtdxdy\Big)^{\frac{1}{m} },
  \end{multline*}
  where the constant $C_{T,K}$ is chosen below.
  It is evident that \ref{(i)} holds, by linearity of $f\mapsto U^f$ and the Minkowski inequality.

  By \cref{cor:integral_estimates}(ii),(iii) (with $G=B^H$),
  \begin{align*}
    \E\Big|U^f_t-U^f_s\Big|^m&\le C\|f\|^m_{L_q\B^{-\eta}_p}|t-s|^{m(\gamma_*+H)},\\
    \E\Big|U^{f(\cdot+x)}_t-U^{f(\cdot+x)}_s-U^{f(\cdot+y)}_t+U^{f(\cdot+y)}_s\Big|^m&\le C\|f\|^m_{L_q\B^{-\eta}_p}|x-y|^m|t-s|^{m\gamma_*}.
  \end{align*}
  Since $m(\gamma_*-\gamma_1)>-1$ and $m(1-\gamma_2)>-1\ge-d$, the integrals $\iint_{[0,T]^2}|t-s|^{m(\gamma_*-\gamma_1)}dsdt$ and $\iint_{K^2}|x-y|^{m(1-\gamma_2)}dxdy$ are finite. Hence $\E[\Xi_{T,K}(f)^m]\le C\|f\|^m_{L_q\B^{-\eta}_p}$, which implies \ref{it.xif} by Jensen's inequality.

  For bounded continuous $f\in L_qL_p$, the field $(t,x)\mapsto U^{f(\cdot+x)}_t=\int_0^tf(r,B^H_r+x)\,dr$ is continuous by dominated convergence, and the multiparameter Garsia--Rodemich--Rumsey inequality \cite[Theorem~3.1]{HuLe13} (applied as in \cite[proof of Lemma~4.1]{anzeletti2023pathbypath}) yields, for a suitable choice of $C_{T,K}$,
  \begin{equation*}
        |U^f_t-U^f_s|\le\Xi_{T,K}(f)|t-s|^{\gamma_1-\frac2m},
  \end{equation*}
  \begin{equation*}
    |U^{f(\cdot+x)}_t-U^{f(\cdot+x)}_s-U^{f(\cdot+y)}_t+U^{f(\cdot+y)}_s|\le \Xi_{T,K}(f)|t-s|^{\gamma_1-\frac2m}|x-y|^{\gamma_2-\frac{2d}m}
  \end{equation*}
  for all $s,t\in[0,T]$ and $x,y\in K$, almost surely. Here we use that $K$ is a ball, so that the Lebesgue measure on $K$ charges every neighbourhood of every point of $K$. Since $\gamma_1-\frac2m>\gamma$, $\gamma_2-\frac{2d}m>1-\varepsilon$ and $0\in K$, this gives \eqref{est.asbounded} and \eqref{est.asLip} for bounded continuous $f$, after enlarging $C_{T,K}$ by a factor depending on $T$ and the diameter of $K$.

  Finally, we show \ref{(ii)} in general, following \cite[proof of Lemma~4.3]{anzeletti2023pathbypath}.
  Note that $L_qL_p\hookrightarrow L_q\B^{-\eta}_p$ so that \ref{(i)} and \ref{it.xif} hold for functions in  $L_qL_p$. Furthermore, by \cref{P:lpq}, $U^f=\int_0^\cdot f(r,B^H_r)dr$ for every $f\in L_qL_p$. Recall that $q<\infty$ by \cref{rem:Lambda_q_infty}, and fix $r_0\in[1,p]\cap[1,\infty)$ with $d(\frac1{r_0}-\frac1p)\le\eta$, so that $L_{r_0}\hookrightarrow\B^0_{r_0}\hookrightarrow\B^{-\eta}_p$.
    In particular, if $|h|\le g$ pointwise with $g\in L_qL_{r_0}$, then $\|h\|_{L_q\B^{-\eta}_p}\le C\|g\|_{L_qL_{r_0}}$.

  \textit{Step 0: localization.} Let $r_K$ be the radius of $K$ and, for $R\in\N$, put $\Omega_R\defby\{\sup_{r\le T}|B^H_r|+r_K<R\}$ and $g_R\defby g\I_{[0,T]\times B_R}$ for a function $g$, where $B_R\defby\{z\in\Rd:|z|<R\}$. For $\omega\in\Omega_R$, $x\in K$ and $s\le t$, we have $\int_s^tg(r,B^H_r(\omega)+x)\,dr=\int_s^tg_R(r,B^H_r(\omega)+x)\,dr$. Moreover, if $g\in L_qL_p$, then by \cref{P:lpq} and Fubini's theorem, almost surely on $\Omega_R$ we have $U^{g(\cdot+x)}=U^{g_R(\cdot+x)}$ for almost every $x\in K$, and hence $\Xi_{T,K}(g)=\Xi_{T,K}(g_R)$. Since $\bigcup_R\Omega_R=\Omega$, it therefore suffices to prove \ref{(ii)} for functions supported in $[0,T]\times B_R$, for each fixed $R\in\N$, which we assume from now on. All open sets below are contained in the bounded set $D_R\defby(-1,T+1)\times B_{R+1}$; for open $U\subset D_R$ we have $\|\I_U\|_{L_qL_{r_0}}\le C(R,T)|U|^{\min(\frac1q,\frac1{r_0})}$ by H\"older's inequality.

    \textit{Step 1: indicators of open sets.} Let $U\subset D_R$ be open. There are continuous functions $g_k$ with compact support in $U$ such that $0\le g_k\uparrow\I_U$ pointwise. By dominated convergence (in $x$ for each $r$, and then in $r$, using that $U$ is bounded), $\|\I_U-g_k\|_{L_qL_{r_0}}\to0$, hence $\|\I_U-g_k\|_{L_q\B^{-\eta}_p}\to0$, and by \ref{(i)} and \ref{it.xif}, along a subsequence, $\Xi_{T,K}(g_k)\to\Xi_{T,K}(\I_U)$ almost surely. By monotone convergence, $\int_s^tg_k(r,B^H_r+x)\,dr\to\int_s^t\I_U(r,B^H_r+x)\,dr$ for every $\omega$, $s\le t$ and $x$. Passing to the limit in \eqref{est.asbounded} and \eqref{est.asLip} for $g_k$, we obtain them for $\I_U$ on an event of full probability.

    \textit{Step 2: null sets.} Let $\mathcal{N}\subset[0,T]\times B_R$ be a Lebesgue null set, and choose open sets $U_j$ with $\mathcal{N}\subset U_j\subset D_R$ and $|U_j|\le2^{-j}$. Then $\|\I_{U_j}\|_{L_q\B^{-\eta}_p}\le C\|\I_{U_j}\|_{L_qL_{r_0}}\to0$, so that $\E\,\Xi_{T,K}(\I_{U_j})\to0$ and, along a subsequence, $\Xi_{T,K}(\I_{U_j})\to0$ almost surely. By Step 1, on an event of full probability,
  \begin{equation*}
    \sup_{z\in K}\int_0^T\I_{\mathcal{N}}(r,B^H_r+z)\,dr\le\sup_{z\in K}\int_0^T\I_{U_j}(r,B^H_r+z)\,dr\le\Xi_{T,K}(\I_{U_j})T^\gamma\to0.
  \end{equation*}
  Hence there is an event $\Omega_{\mathcal{N}}$ of full probability on which $\int_0^T\I_{\mathcal{N}}(r,B^H_r(\omega)+z)\,dr=0$ for all $z\in K$.

    \textit{Step 3: bounded functions.} Let $f\in L_qL_p$ with $|f|\le M$, supported in $[0,T]\times B_R$. By Lusin's theorem, for each $n$ there is a continuous function $g^n$ with $|g^n|\le M$ such that $|\{f\ne g^n\}|\le2^{-n-1}$. Let $\chi$ be a continuous function with $0\le\chi\le1$, $\chi=1$ on $\overline{B_R}$ and $\chi=0$ outside $B_{R+1}$, and put $f^n(r,z)\defby\chi(z)g^n(r,z)$. Since $f=0$ outside $[0,T]\times B_R$, we have $\{f^n\ne f\}\cap([0,T]\times\Rd)\subset\{g^n\ne f\}\cap([0,T]\times B_{R+1})$, and by outer regularity of the Lebesgue measure there is an open set $U^n\subset D_R$ with $|U^n|\le2^{-n}$ such that $f^n=f$ on $([0,T]\times\Rd)\setminus U^n$. Moreover, $f^n$ is continuous, $|f^n|\le M$, and $f^n$ has bounded support. Then $|f-f^n|\le2M\I_{U^n}$, so that $f^n\in L_qL_p$ and $\|f-f^n\|_{L_q\B^{-\eta}_p}\le C\|f-f^n\|_{L_qL_{r_0}}\le 2CM\|\I_{U^n}\|_{L_qL_{r_0}}\to0$. Along a subsequence, $\Xi_{T,K}(f^n)\to\Xi_{T,K}(f)$ and $\Xi_{T,K}(\I_{U^n})\to0$ almost surely. By Step 1, for $x\in K$ and $s\le t$,
  \begin{equation*}
    \Big|\int_s^t(f-f^n)(r,B^H_r+x)\,dr\Big|\le2M\int_0^T\I_{U^n}(r,B^H_r+x)\,dr\le2M\,\Xi_{T,K}(\I_{U^n})T^\gamma.
  \end{equation*}
  Writing $f=f^n+(f-f^n)$, applying \eqref{est.asbounded} and \eqref{est.asLip} to the bounded continuous functions $f^n$ and letting $n\to\infty$, we obtain \eqref{est.asbounded} and \eqref{est.asLip} for $f$ on an event of full probability.

    \textit{Step 4: general $f\in L_qL_p$ supported in $[0,T]\times B_R$.}    To remove the boundedness restriction, we use a truncation procedure. For each integer $n\ge1$, define $f^n=f\1_{|f|\le n}$. 
  Note that $(f^n)$ and $(|f^n|)$ are sequences of bounded functions converging to $f$ and $|f|$, respectively, in $L_qL_p$ and almost everywhere. Using \ref{it.xif}, we can choose a subsequence, still denoted by $n$, such that $\Xi_{T,K}(f^n)\to \Xi_{T,K}(f)$ and $\Xi_{T,K}(|f^n|)\to \Xi_{T,K}(|f|)$  a.s.. Define $\mathcal{N}=\{(r,z)\in[0,T]\times\Rd: \lim_nf^n(r,z)=f(r,z)\text{ and }\lim_n|f^n|(r,z)=|f|(r,z)\}^c$ and note that $|\mathcal{N}|=0$. Let $\Omega^*\defby\Omega_{\mathcal{N}}$ be the event of full probability from Step 2.
  Then for every $\omega\in\Omega^*$, $f^n(r,B^H_r(\omega)+x)\to f(r,B^H_r(\omega)+x)$ and $|f^n|(r,B^H_r(\omega)+x)\to |f|(r,B^H_r(\omega)+x)$ for a.e. $r\in [0,T]$. 
  Using this, \eqref{est.asbounded} and monotone convergence theorem, we then have for every $\omega\in\Omega^*$, $x\in K$,
  \begin{align*}
    \int_0^T |f(r,B^H_r(\omega)+x)|dr
    &=\lim_n\int_0^T |f^n(r,B^H_r(\omega)+x)|dr
    \le \Xi_{T,K}(|f|)(\omega)T^\gamma.
  \end{align*}
  This shows that for a.s. $\omega$, for every $x\in K$, the function $t\mapsto f(t,B^H_t(\omega)+x)$ is integrable on $[0,T]$. Similarly, by dominated convergence theorem, we have for every $\omega\in\Omega^*$, for every $(t,x)\in [0,T]\times K$, 
  \begin{align*}
    \int_0^t f(r,B^H_r(\omega)+x)dr
    &=\lim_n\int_0^t f^n(r,B^H_r(\omega)+x)dr.
  \end{align*}
    We then apply \eqref{est.asbounded} and \eqref{est.asLip} for the bounded functions $f^n$ (Step 3) and pass to the limit as $n\to\infty$ to obtain the desired estimates for $f$, completing the proof of \ref{(ii)}. 
\end{proof}
In what follows, $\norm{\psi}_{1\vr;[s,t]}$ denotes the total variation of a path $\psi$ on $[s,t]$, see \cref{subsec.notation}.
\begin{lemma}[Perturbation by a path of finite variation]\label{lem:fbm_regularizing}
  Assume \eqref{eq:parameters}.
    Let $f:[0,T]\times\Rd \to \mathbb{R}$ be a function in $L_qL_p$, $K\subset\Rd$ be a closed ball centred at $0$, $\gamma \in (0,1-H-\frac{Hd}{p}-\frac{1}{q})$ and $\varepsilon\in(0,\gamma)$.
  Let $\eta\in(0,\eta_0]$ be such that $\gamma<1-H-\frac{Hd}{p}-\frac1q-\eta H$, and let $\Xi_{T,K}=\Xi_{T,K,\gamma,\varepsilon,\eta}$ and $\Omega_{f,T,K}$ be as in \cref{lem.regularizing}, where we assume, enlarging $\Omega_{f,T,K}$ if necessary, that $\Omega_{f,T,K}\subset\Omega_{|f|,T,K}$ (integrability of $r\mapsto f(r,B^H_r+\psi_r)$ requires \eqref{est.asbounded}--\eqref{est.asLip} for $|f|$ as well).
  Then there is a deterministic constant $C=C(\varepsilon,\gamma)$ such that for every $\omega\in \Omega_{f,T,K}$, every $(s, t) \in \Delta_{[0, T]}$ and every function  $\psi:[0,T]\to K$  of finite variation, the map $r \mapsto f(r, B^H_r(\omega) + \psi_r)$ is integrable on $[0,T]$ and
  \begin{gather}
    \Big|\int_s^t f(r,B^H_r(\omega)+\psi_r)\dd r\Big|\le C \Xi_{T,K}(f)(\omega)(1+\norm{\psi}_{1\vr;[s,t]}^{1- \varepsilon})(t-s)^{\gamma}\label{est.intfome}
    \\\shortintertext{and}
    \Big|\int_s^t f(r,B^H_r(\omega)+\psi_r)\dd r-\int_s^t f(r,B^H_r(\omega)+\psi_s)\dd r\Big|\le C \Xi_{T,K}(f)(\omega)\norm{\psi}_{1\vr;[s,t]}^{1- \varepsilon}(t-s)^{\gamma}. \label{est.young}
  \end{gather}
\end{lemma}
\begin{proof}
  The stated result is an analogue of \cite[Lemma~4.3]{anzeletti2023pathbypath}, whose proof can be applied here with the necessary estimates provided by \cref{lem.regularizing}. 
\end{proof}
\begin{corollary}\label{cor:continuity_along_B}
  Assume \eqref{eq:parameters}.
  Let  $f:[0,T]\times\Rd \to \mathbb{R}$ be a function in $L_qL_p$.  
  Then there exists an event $\Omega_{f,T}$ of full probability such that for every $\omega\in \Omega_{f,T}$
  and for a sequence $(\psi^m_t)_{t \in [s, T]}$ of paths such that 
  \begin{itemize}
    \item $\liminf_m \norm{\psi^m}_{1\vr; [s, T]} < \infty$ and 
    \item $\lim_{m \to \infty} \psi^m = \psi$ in $C([s, T])$,
  \end{itemize}
  we have $\int_s^T \abs{f(r, B_r^H(\omega) + \psi_r)} \dd r < \infty$ and 
  \begin{align*}
    \liminf_{m \to \infty} \sup_{t \in [s, T]} \abs[\Big]{ \int_s^t f(r, B_r^H(\omega)+ \psi^m_r) \dd r 
    - \int_s^t f(r, B_r^H(\omega) + \psi_r) \dd r } = 0.
  \end{align*}
\end{corollary}
\begin{proof}
  The proof is similar to \cite[Corollary~4.4]{anzeletti2023pathbypath}. Firstly, we note that for every partition $\pi$ of $[s, T]$ we have 
  \begin{align*}
    \sum_{[t_1, t_2] \in \pi} \abs{\psi_{t_2} - \psi_{t_1}} = \liminf_{m} \sum_{[t_1, t_2] \in \pi} 
    \abs{\psi^m_{t_2} - \psi^m_{t_1}} 
    \leq \liminf_m \norm{\psi^m}_{1\vr; [s, T]}.
  \end{align*}
  Therefore, 
  \begin{align}\label{eq:psi_var_by_psi_m_var}
    \norm{\psi}_{1\vr; [s, T]} \leq \liminf_m \norm{\psi^m}_{1\vr; [s, T]}<\infty.
  \end{align}

    Let $K^i\defby\{x\in\Rd:|x|\le i\}$, $i\ge1$. 
    Fix $\gamma$, $\varepsilon$ and $\eta$ as in \cref{lem:fbm_regularizing}, and let $\Xi_{T,K}=\Xi_{T,K,\gamma,\varepsilon,\eta}$. For each $g\in L_qL_p$, let $\Omega_{g, T, K^i}$ be the event of full probability such that the statements of \cref{lem.regularizing,lem:fbm_regularizing} hold.
    Let $(f^n)_{n=1}^{\infty}$ be a subsequence of the mollifying sequence of $f$ from \cref{lem:mollifying}; each $f^n$ is bounded and continuous in the space variable, and $f^n\to f$ in $L_q\B^{-\eta}_p$. By \cref{lem.regularizing} \ref{it.xif} and a diagonal argument, we may choose the subsequence such that
  for some event $\tilde{\Omega}$ of full probability and for every $i$ we have
  \begin{align}\label{eq:Xi_continuity}
    \lim_{n \to \infty} \Xi_{T, K^i}(f - f^n) = 0 \quad \text{on } \tilde{\Omega}.
  \end{align}
  We set 
  \begin{align*}
    \Omega_{f, T} \defby \tilde{\Omega} \cap \Big( \cap_i \Omega_{f, T, K^i} \Big) \cap \Big( \cap_{i, n} \Omega_{f-f^n, T, K^i} \Big).
  \end{align*}
    Let $K = K^i$ be such that $\psi_t, \psi^n_t \in K$ for every $t \in [s, T]$ and $n$.
  In view of \eqref{eq:psi_var_by_psi_m_var} and \eqref{est.intfome}, for $\omega \in \Omega_{f, T}$ we have 
  \begin{align*}
    \int_s^T \abs{f(r, B^H_r(\omega) + \psi_r)} \dd r < \infty.
  \end{align*}
  Furthermore, for $\omega \in \Omega_{f, T}$, by the triangle inequality 
  \begin{align*}
    \MoveEqLeft[7]
    \abs[\Big]{\int_s^t [f(r, B^H_r(\omega) + \psi_r) - f(r, B^H_r(\omega) + \psi^m_r)] \dd r} \\
    &\leq \abs[\Big]{\int_s^t [f(r, B^H_r(\omega) + \psi_r) - f^n(r, B^H_r(\omega) + \psi_r)] \dd r} \\
    &\phantom{\leq} + \abs[\Big]{\int_s^t [f(r, B^H_r(\omega) + \psi^m_r) - f^n(r, B^H_r(\omega) + \psi^m_r)] \dd r} \\
    &\phantom{\leq} + \abs[\Big]{\int_s^t [f^n(r, B^H_r(\omega) + \psi_r) - f^n(r, B^H_r(\omega) + \psi^m_r)] \dd r}. 
  \end{align*}
  To estimate the first two terms, applying \eqref{est.intfome} and \eqref{eq:psi_var_by_psi_m_var}, we get 
  \begin{multline*}
   \sup_{t \in [s, T]} \abs[\Big]{\int_s^t [f(r, B^H_r(\omega) + \psi_r) - f^n(r, B^H_r(\omega) + \psi_r)] \dd r} \\
   \leq C \Xi_{T, K}(f-f^n, \omega) (1 + \liminf_m \norm{\psi^m}_{1\vr; [s, T]}^{1-\varepsilon}),
  \end{multline*}
  \begin{multline*}
   \liminf_{m} \sup_{t \in [s, T]} \abs[\Big]{\int_s^t [f(r, B^H_r(\omega) + \psi_r^m) - f^n(r, B^H_r(\omega) + \psi^m_r)] \dd r} \\
   \leq C \Xi_{T, K}(f-f^n, \omega) (1 + \liminf_m \norm{\psi^m}_{1\vr; [s, T]}^{1-\varepsilon}).
  \end{multline*}
    To estimate the third term, since $f^n$ is bounded and continuous in the space variable, dominated convergence gives, for each fixed $n$,  
  \begin{align*}
    \lim_{m} \sup_{t \in [s, T]} \abs[\Big]{\int_s^t [f^n(r, B^H_r(\omega) + \psi_r) - f^n(r, B^H_r(\omega) + \psi^m_r)] \dd r} = 0.
  \end{align*}
  It remains to apply \eqref{eq:Xi_continuity}.
\end{proof}

\begin{lemma}(A priori bound)\label{lem:a_priori_bound}
  Assume \eqref{eq:parameters}.
    Let  $K\subset\Rd$ be a closed ball centred at $0$, $\gamma \in (0,1-H-\frac{Hd}{p}-\frac{1}{q})$,  
    $\varepsilon \in (0, \gamma)$ and $b \in L_qL_p$.
    Let $\Xi_{T, K}$ and the events $\Omega_{|b|,T,K}$ be as in \cref{lem:fbm_regularizing}.
  Then, there exists a finite constant $C = C(\gamma, \varepsilon)$ with the following property. Let $\omega\in\Omega_{|b|,T,K}$, $0\le u<v\le T$, and let $Y\colon[u,v]\to\Rd$ be a solution to \eqref{fSDE} on $[u,v]$ along $B^H(\omega)$, that is, $\int_u^v|b(r,Y_r)|\dd r<\infty$ and $Y_t=Y_u+\int_u^tb(r,Y_r)\dd r+B^H_t(\omega)-B^H_u(\omega)$ for $t\in[u,v]$. If $Y_t - B^H_t(\omega) \in K$ for all $t\in[u,v]$, then
  \begin{equation}\label{eq:a_priori_bound}
        \int_u^v \abs{b(r, Y_r)} \dd r \leq C\Xi_{T, K}(\abs{b})(\omega) (v-u)^{\gamma} + C
    (\Xi_{T, K}(\abs{b})(\omega) (v-u)^{\gamma} )^{1/\varepsilon}.
  \end{equation}
\end{lemma}
\begin{proof}
  This is \cite[Proposition~5.1]{anzeletti2023pathbypath}, whose proof only uses \eqref{est.intfome}; we repeat it for convenience. Put $\psi\defby Y-B^H(\omega)$, so that $\norm{\psi}_{1\vr;[u,v]}\le\int_u^v|b|(r,Y_r)\dd r=:I$. By \eqref{est.intfome} with $f=|b|$,
  \begin{equation*}
    I\le C\,\Xi_{T,K}(|b|)(\omega)(1+I^{1-\varepsilon})(v-u)^\gamma.
  \end{equation*}
  By Young's inequality, $C\,\Xi_{T,K}(|b|)(\omega)I^{1-\varepsilon}(v-u)^\gamma\le\frac12I+C'(\Xi_{T,K}(|b|)(\omega)(v-u)^\gamma)^{1/\varepsilon}$, and \eqref{eq:a_priori_bound} follows (after renaming constants, and since $I<\infty$).
\end{proof}
\subsection{Existence of strong solutions}
\label{sub:existence_strong_sol}
  We first show that solutions are invariant in the sense that the limit of solutions is also a solution. While this property is trivial when $b$ is continuous, it is much less so when $b$ is only in $L_qL_p$. 
\begin{proposition}\label{prop.invariance}
  Assume that \eqref{eq:parameters} holds and recall $\eta_0$ is the constant from \cref{lem.minusoneOK}.  Let $T>0$,  $x\in\R^d$, $b\in L_q([0,T];L_p(\R^d))$.
  Let $(x^n)_{n\in\Z_+}$ be a sequence of vectors in $\R^d$ converging to $x$ as $n\to\infty$.
  Let $(b^n)_{n\in\Z_+}$ be a sequence of functions in $L_q([0,T];L_p(\R^d))$ such that 
  \begin{align}
    \lim_{n\to\infty}\|b^n-b\|_{L_q([0,T];\B^{-\eta_0}_p(\Rd))}=0\quad\text{and}\quad \sup_{n}\|b^n\|_{L_q([0,T];L_p(\Rd))}<\infty. \label{con:bnb}
  \end{align} 
  Let $X^n$ be a strong solution to  \Gref{x^n;b^n} for each $n$.
  Assume that $X^n \to X$ in $C([0, T])$ in probability.
  \noindent 

  Then $X$ is a strong solution of \Gref{x;b}. Furthermore, we have
    \begin{align}
      \lim_{n\to\infty} \sup_{t \in [0, T]} \abs[\Big]{\int_0^t b^n(r,X^n_r)dr-\int_0^t b(r,X_r)dr} = 0 \quad \text{in probability}.\label{limintxn1}
    \end{align}
\end{proposition}
\begin{proof}
  Through subsequences, we can assume that $X^n \to X$ in $C([0, T])$ almost surely.
  We define $\psi^n=X^n-B^H$ so that $\psi^n \to \psi$ in $C([0, T])$ almost surely.
    Since each $X^n$ is a strong solution, $X$ is adapted to the standard augmentation of the filtration generated by $B^H$.
We show that
\begin{equation}\label{eq:psi_consistent}
  \psi_\cdot = x + \int_0^{\cdot} b(r, X_r) \dd r\quad \text{a.s.},
\end{equation}
which means that $X$ solves \Gref{x;b}. Then \eqref{limintxn1} follows, since $\int_0^\cdot b^n(r,X^n_r)\,dr=\psi^n-x^n\to\psi-x=\int_0^\cdot b(r,X_r)\,dr$ uniformly on $[0,T]$ in probability.

Fix $\gamma\in(0,1-H-\frac{Hd}p-\frac1q-\eta_0H)$ and $\varepsilon\in(0,\gamma)$, and let $\Xi_{T, K}=\Xi_{T,K,\gamma,\varepsilon,\eta_0}$ be as in \cref{lem.regularizing}.
 Let $K_m\defby\{x\in\Rd:|x|\le m\}$, $m\ge1$.

  By \cref{lem.regularizing} \ref{it.xif}, we have $\Xi_{T, K_m}(b) < \infty$ for all $m$ almost surely. 
For each $m$, \cref{lem.regularizing} and \eqref{con:bnb} imply that  $\sup_n \E\Xi_{T, K_m}(\abs{b^n}) < \infty$ 
and $\lim_{n \to \infty} \E|\Xi_{T, K_m}(b^n - b) |= 0$.
 Taking subsequences if necessary, we may suppose that 
$\liminf_{n} \Xi_{T, K_m}(\abs{b^n}) < \infty$  and $\lim_{n \to \infty} \Xi_{T, K_m}(b - b^n) = 0$ for all $m$ almost surely. 

Let $\Omega_{b,(b^n),T}$ be the intersection of the event $\Omega_{b,T}$ from \cref{cor:continuity_along_B}, the events $\Omega_{b-b^n,T,K_m}$ and $\Omega_{|b^n|,T,K_m}$, $n,m\ge1$, from \cref{lem.regularizing}, and the event on which $\lim_{n} \Xi_{T, K_m}(b - b^n) = 0$ and $\liminf_{n} \Xi_{T, K_m}(\abs{b^n}) < \infty$ for all $m$. Such event has full probability and depends only on $b$, $(b^n)$ and $T$. We fix an $\omega\in \Omega_{b,(b^n),T}$ such that, in addition, for all $n$, $X^n(\omega)$ solves \Gref{x^n;b^n} along $B^H(\omega)$ and $X^n(\omega)\to X(\omega)$ in $C([0,T])$ (this holds almost surely).
Since $\psi^n(\omega) \to \psi(\omega)$ in $C([0, T])$, we can find a $K = K_{m(\omega)}$ such that $\psi^n(\omega) \in K$ for all $n$. 
Then, the a priori estimate, \cref{lem:a_priori_bound}, implies that $\liminf_n \norm{\psi^n(\omega)}_{1\vr;[0,T]} < \infty$.
Hence there is a subsequence $(n_k)$, depending on $\omega$, such that $\sup_k\norm{\psi^{n_k}(\omega)}_{1\vr;[0,T]}<\infty$.
Applying \cref{cor:continuity_along_B} to the sequence $(\psi^{n_k}(\omega))_k$ and passing to a further subsequence, still denoted by $(n_k)$, we get
\begin{equation}\label{eq:b_X_vs_b_X_n}
  \lim_{k \to \infty} \sup_{t \in [0, T]} \abs[\Big]{\int_0^t [b(r, X^{n_k}_r(\omega)) - b(r, X_r(\omega))] \dd r}  = 0.
\end{equation}
Furthermore, by \eqref{est.intfome},
\begin{equation*}
   \sup_{t \in [0, T]} \abs[\Big]{\int_0^t [b(r, X^{n_k}_r(\omega)) - b^{n_k}(r, X^{n_k}_r(\omega))] \dd r}\les \Xi_{T, K}(b-b^{n_k})(\omega) \big(1+\sup_j\norm{\psi^{n_j}(\omega)}_{1\vr;[0,T]} \big)\to0
\end{equation*}
as $k\to\infty$. Since $\psi^{n_k}_t(\omega)-x^{n_k}=\int_0^tb^{n_k}(r,X^{n_k}_r(\omega))\,dr$, these two estimates imply that
\begin{align*}
  \lim_{k \to \infty} \sup_{t \in [0, T]} \abs[\Big]{\psi^{n_k}_t(\omega)-x^{n_k}-\int_0^t b(r, X_r(\omega)) \dd r}  = 0 ,
\end{align*}
which yields \eqref{eq:psi_consistent}, since $\psi^{n_k}(\omega)\to\psi(\omega)$ and $x^{n_k}\to x$.
\end{proof}
We are now ready to prove the existence of strong solutions to \Gref{x;b}.

\begin{proof}[Proof of \cref{l:strongsol}]
  Let $(b^n)_{n\in\Z_+}$ be the mollifying sequence of $b$ as in \cref{lem:mollifying}. Then each $b^n$ has bounded spatial derivatives and condition \eqref{con:bnb} holds.
  Let $X^n$ be the unique strong solution to  \Gref{x;b^n}, which is obtained through a fixed-point argument.
  
Fix $\ell\in(\frac{1}{1-H-\frac{Hd}{p}+\frac{Hd}{p\vee2}},2)$.
Since $b^n$ is bounded, \cref{lem:girsanov_sol} is applicable to $X^n$ (see \cref{r:boundedisok}); together with \eqref{GGbound} and \cref{lem:mollifying}(ii),(iii), it yields
\begin{equation*}
  M\defby\sup_{n}\E\exp\Big(12\Big\|\int_0^\cdot b^n(r,X^n_r)\,dr\Big\|^2_{\H}\Big)<\infty.
\end{equation*}
Hence, \eqref{con.expXi} holds for $(X^{n_1},X^{n_2})$ with $\lambda=12$, uniformly in $n_1,n_2$, which is all that is needed for \eqref{Vbound11}, \eqref{Vbound22} with $m=1$.  
Let $\delta>0$, $L\ge1$  be arbitrary. Let $n_1,n_2\in\Z_+$. Put for $t\in[0,T]$
\begin{equation*}
	V^{n_1,n_2}_t:=\int_0^t \int_0^1 \nabla b^{n_1}  (r,\theta X^{n_1}_r+ (1 - \theta) X_r^{n_2})\, d\theta d r,\quad R^{n_1,n_2}_t:=\int_0^t (b^{n_1}-b^{n_2})(r,X_r^{n_2})\,dr.
\end{equation*}

We apply estimate \eqref{eq:pathwise_stability} of \cref{l:pathwise.stabitliy} and then \eqref{Vbound11}, \eqref{Vbound22} with $m=1$ and $\eta=0$. Using the Chebyshev inequality we derive
\begin{align*}
&\P(\sup_{t \in [0, T]} |X^{n_1}_t - X^{n_2}_t| \ge \delta)\\
&\quad\le 
	\P(\|V^{n_1,n_2}\|_{\ell\var;[0,T]} \geq L)
	+\P(e^{CL^2} (|x^{n_1} - x^{n_2}| + \|R^{n_1,n_2}\|_{\ell\var;[0,T]}) \ge \delta)\\
&\quad \le C L^{-1}	+ Ce^{CL^2} (|x^{n_1} - x^{n_2}| +\|b^{n_1}-b^{n_2}\|_{L_q([0,T];\B^{-1}_p)})\delta^{-1},
\end{align*}
where $C=C(M,\ell,T,\sup_n\|b^n\|_{L_qL_p})$. Sending in the above inequality first $n_1,n_2\to\infty$ and then $L\to\infty$ we get 
$$
\lim_{n_1,n_2\to\infty}\P(\sup_{t \in [0, T]} |X^{n_1}_t - X^{n_2}_t|\ge\delta)=0,
$$
where we used \eqref{con:bnb} and that $x^n\to x $ as $n\to\infty$. 
Hence, there exists a process $X$ such that 
$X^n$ converges to $X$ uniformly on $[0,T]$ in probability. 
It follows from \cref{prop.invariance} that $X$ is a strong solution to \Gref{x;b}.

Finally, we prove \eqref{nboundsolnew} and \eqref{nboundsolquant}. Fix $\lambda>0$ and write $\|\cdot\|_{\W}\defby\|\cdot\|_{\W^{H+\frac12+\eps}_2([0,T])}$.
If $p=\infty$, then \eqref{novikovbound2} and \cref{lem:mollifying}(ii) give $\|\int_0^\cdot b^n(r,X^n_r)\,dr\|_{\W}\le\bar C\|b\|_{L_qL_p}$ almost surely, for every $n$.
If $2Hd<p<\infty$, then \eqref{novikovbound1b} and \cref{lem:mollifying}(ii) give
\begin{equation*}
  \E \exp \Big( \lambda \Big\|\int_0^\cdot b^n(r,X^n_r)\,dr\Big\|_{\W}^2 \Big)\le C'(T,\lambda,\|b^n\|_{L_qL_p})\le C'(T,\lambda,\|b\|_{L_qL_p}).
\end{equation*}
If $p=2Hd$, then \eqref{novikovbound1} with some $\eta'=\eta'(T,\lambda)$ and \cref{lem:mollifying}(iii) give
\begin{equation*}
  \E \exp \Big( \lambda \Big\|\int_0^\cdot b^n(r,X^n_r)\,dr\Big\|_{\W}^2 \Big)\le Ce^{c\Lambda_{b^n}^2(\eta')}\le Ce^{4c\Lambda_{b}^2(\eta'/2)}.
\end{equation*}
By \eqref{limintxn1}, along a subsequence, $\int_0^\cdot b^n(r,X^n_r)\,dr\to\int_0^\cdot b(r,X_r)\,dr=X-B^H-x$ uniformly on $[0,T]$ almost surely. Since $\|\cdot\|_{\W}$ is lower semicontinuous with respect to uniform convergence (by Fatou's lemma applied to the double integral defining it), Fatou's lemma yields \eqref{nboundsolquant} with $\eta\defby\eta'/2$, and in particular \eqref{nboundsolnew}.
\end{proof}






\section{Stability in moment norms}\label{sec:stability}
In \cref{l:pathwise.stabitliy}, we proved the pathwise stability estimate \eqref{eq:pathwise_stability} for solutions which satisfy exponential integrability condition \eqref{con.expXi}. In this section, we prove stability in moment norms for equations with smooth drifts (\cref{lem:moment_stability}), with constants depending on the drifts only through $\|b^i\|_{L_qL_p}$ and $\Lambda_{b^i}$; the extension to general drifts is carried out by approximation in the proof of \cref{thm:main_stability} (\cref{sec:proof_main_results}). The smooth case suffices for constructing the semiflow in \cref{sec:flow_and_path_by_path_uniqueness}.
\begin{proposition}[Stability in moment norms]\label{lem:moment_stability}
  Suppose that the condition \eqref{eq:parameters} is satisfied.
  Let $T > 0$,  $x^1, x^2 \in \mathbb{R}^d$, $b^1, b^2 \in L_q L_p$ be smooth functions with bounded derivatives, and 
    $X^i$ be the strong solution to \Gref{x^i;b^i}, and put $\gamma\defby1-H-\frac{Hd}{p}-\frac1q$.
  Then
  for every $m \in [1, \infty)$ and 
  $0 \leq s < t \leq T$,
  we have
  \begin{equation*}
    \norm{X^1_t - X^1_s - (X^2_t - X^2_s)}_{L_m(\Omega)} 
    \leq C e^{C(1\vee\norm{b^1}_{L_qL_p})^{1/{\gamma^{2}} }} (\abs{x^1 - x^2} + \norm{b^1-b^2}_{L_q \mathcal{B}_p^{-1}})
    (t-s)^{1 - H - \frac{Hd}{p} - \frac{1}{q}}.
  \end{equation*}
  In the above, $C$ is some finite constant of the form 
  \begin{align*}
    C(m, T, \norm{b^1}_{L_q L_p}, \norm{b^2}_{L_q L_p})\I_{(p>2Hd)}+\bar C(m, T, \Lambda_{b^1}(\eta ), \Lambda_{b^2}(\eta ))\I_{(p=2Hd)} 
  \end{align*} 
  for some $\eta =\eta (m,T)$.  
\end{proposition}
\begin{remark}\label{rem:dependence_of_const}
  The dependence of the constants on $\norm{b^i}_{L_q L_p}, \Lambda_{b^i}(\eta )$ is always increasing.
  We will not mention this dependence in the sequel. For $p>2Hd$, the constant $C$ already depends on $\|b^1\|_{L_qL_p}$ through the Girsanov bounds, so the explicit factor $e^{C(1\vee\|b^1\|_{L_qL_p})^{1/\gamma^2}}$ is informative only for $p=2Hd$, where $\bar C$ depends on $b^i$ only through $\Lambda_{b^i}(\eta)$.
\end{remark}
A straightforward strategy to deduce the stability estimate in $L_m(\Omega)$ is via the pathwise stability estimate \eqref{eq:pathwise_stability}. 
However, this is only possible provided that (recall $V$ from \eqref{eq:def_V_and_R})
\begin{align}\label{eq:rho_var_mom_desired}
  \mathbb{E}[e^{\lambda \norm{V}_{\ell\vr}^{\ell}}] < \infty  
\end{align}
for some  $\ell \in ((1-H-\frac{Hd}{p} + \frac{Hd}{p \vee 2})^{-1}, 2)$  and for any  $\lambda > 0$.
In view of \cref{rem:var_exp_bd} and Girsanov's argument (as in the proof of \cref{l:pathwise.stabitliy}), we have 
\begin{align*}
  \mathbb{E}[e^{\lambda \norm{V}_{\ell\vr}^{(H+\frac{Hd}{p})^{-1}}}] < \infty  \quad \text{ for some small } \lambda > 0.
\end{align*}
For \eqref{eq:rho_var_mom_desired} to hold,  it is sufficient that 
\begin{align*}
    \frac{Hd}{p} + H < 1 - H - Hd \Big( \frac{1}{p} - \frac{1}{p \vee 2} \Big),
\end{align*}
which does not necessarily hold under \eqref{eq:parameters} in general.

Hence, in the regime of \eqref{eq:parameters}, the pathwise stability estimate \eqref{eq:pathwise_stability} does not imply stability in moment norms and we need to employ a different approach.
Our strategy is to view \eqref{eq:Y_linearization} as 
a rough differential equation driven by a VMO rough path lifted from $V$ and $R$, and to employ Gr\"onwall's argument in the probabilistic setting 
(i.e., in the space of these VMO rough paths).


%

\subsection{VMO estimates}\label{subsec:vmo_estimate}
To control the moments of the rough path lifts, we employ  VMO-type estimates on $\ell$-variations of $V$ and $R$, which is the goal of this subsection.
\begin{proposition}[VMO estimate on variation]\label{lem:vmo_variation}
  Suppose that the parameters $H, p, q, d$ satisfy \eqref{eq:parameters}.
  Let $b^1, b^2 \in L_q L_p$ be smooth functions with bounded derivatives, and 
  let $\ell \in ((1-H-\frac{Hd}{p}+\frac{Hd}{p\vee 2})^{-1}, 2)$. 
  Let $X^i$ be the strong solution to \Gref{x^i;b^i},
  and we define $V$ and $R$ by \eqref{eq:def_V_and_R}.
    Let $(\mathcal{F}_t)_{t \in [0, T]}$ be the filtration generated by $B^H$.
  Then, for any $m \in [1, \infty)$ and $s < t$, we have
  \begin{align*}
    \norm{\norm{\norm{V}_{\ell\vr;[s,t]}}_{L_m(\Omega) \vert \mathcal{F}_s}}_{L_{\infty}(\Omega)} & 
    \leq C \norm{b^1}_{L_q([s,t]; L_p(\mathbb{R}^d))} (t-s)^{1-H-\frac{Hd}{p} -\frac{1}{q}}, \\
    \norm{\norm{\norm{R}_{\ell\vr;[s,t]}}_{L_m(\Omega) \vert \mathcal{F}_s}}_{L_{\infty}(\Omega)} & 
    \leq C \norm{b^1 - b^2}_{L_q([s,t]; \mathcal{B}_p^{-1})} (t-s)^{1-H-\frac{Hd}{p}-\frac{1}{q}}. 
  \end{align*}
  In the above, $C$ is some finite constant of the form 
  \begin{align*}
    C(m,\ell, T, \norm{b^1}_{L_q L_p}, \norm{b^2}_{L_q L_p})\I_{(p>2Hd)}+\bar C(m, T, \Lambda_{b^1}(\eta ), \Lambda_{b^2}(\eta ))\I_{(p=2Hd)} 
  \end{align*} 
  for some $\eta =\eta (m,T)$. 
\end{proposition}
We note that \cref{lem:vmo_variation} is an improvement of \eqref{Vbound11} and \eqref{Vbound22}. To estimate the conditional moments, one may attempt to use the argument in the proof \cref{l:pathwise.stabitliy} with conditional expectation. Indeed, for each $s$, similar to \eqref{genboundstep1} with the notations therein, we have
\begin{align*}
  \E^s_{\P} F_{s,t}=\E^s_{\wt \P_\theta} F_{s,t} \xi_s\xi_{T}^{-1}
  \le 
  \left(\E^s_{\wt \P_\theta}F_{s,t}^2 \right)^{\frac{1}{2}} \left(\E^s_{\wt \P_\theta}\xi_s^2 \xi_{T}^{-2}\right)^{\frac{1}{2}}.
\end{align*} 
However, while possible, we have not been able to obtain the desired estimates for $\E^s_{\wt \P_\theta}\xi_s^2 \xi_{T}^{-2}$. We therefore need to employ a different approach.

It is more convenient to work with the Mandelbrot--Van Ness representation of $B^H$.
The assertion of \cref{lem:vmo_variation} concerns only the conditional laws of $X^1,X^2$ given $\F_s$, and $X^1,X^2$ are strong solutions. Hence, as explained in \cref{subsec:girsanov_prelim}, we may assume that there is a two-sided Brownian motion $\tilde W$ such that \eqref{eq:mandelbrot_van_ness} holds, and we let $(\tilde{\mathcal F}_t)_{t\in\R}$ be its filtration and $\zeta^v$, $\bar B^v$ be as in \eqref{eq:MVN_decomp}.
Since $\mathcal{F}_t\subseteq \tilde{\mathcal{F}}_t$,
it suffices to prove \cref{lem:vmo_variation} with $\mathcal{F}$ replaced by $\tilde{\mathcal{F}}$ (see the end of the proof below).

We fix a $v \in [0, T)$. 
As we are assuming that $b^1$, $b^2$ are smooth, for any continuous function $\zeta$ with $\zeta_0 = 0$, the differential equation 
\begin{align*}
  d x_t = b^i(v+t, x_t) d t + d \zeta_t, \quad x_0 = y,
\end{align*}
has a unique solution. 
Let $\tilde B^H\defby\bar B^0$ be the Riemann--Liouville fractional Brownian motion \eqref{eq:RL_def} driven by the Brownian motion $(\tilde W_t)_{t\ge0}$.
We denote by $\tilde{X}^i[y, \zeta]$ the unique solution to the SDE 
\begin{align*}
  d \tilde{X}^i_t[y, \zeta] = b^i(v+t, \tilde{X}^i_t[y, \zeta]) d t 
  + d \zeta_t + d \tilde{B}_t^H, \quad 
  \tilde{X}^i_0[y, \zeta] = y.
\end{align*}
Note that $\tilde{X}^i[y, \zeta]$ is adapted to the filtration generated by $\tilde{B}^H$. 

We observe the identity 
\begin{align*}
    X^i_{v+t} = X^i_v + \int_0^{t} b^i(v+r, X^i_{v+r}) \dd r + \zeta^v_t+\bar B^v_t,
\end{align*}
where, by \eqref{eq:MVN_decomp}, $\zeta^v$ is $\tilde{\mathcal F}_v$-measurable and $\bar B^v$ is independent of $\tilde{\mathcal{F}}_v$ and equal to $\tilde{B}^H$ in law.
Since the solution $X^i$ is strong, $X^i_v$ is measurable with respect to $\tilde{\mathcal{F}}_v$. 
Therefore, for any bounded measurable $F$, $v \in [0, T)$ and $t \in [0, T-v]$, we have 
\begin{equation}\label{eq:F_X_i_conditioning}
  \E[F((X^1_{v+r}, X^2_{v+r})_{0 \leq r \leq t}) \vert \tilde{\mathcal{F}}_v] 
  = \E[F((\tilde{X}^1_{r}[y^1, \zeta], \tilde{X}^2_{r}[y^2, \zeta])_{0 \leq r \leq t})] \big\vert_{y^i = X^i_v,\ \zeta = \zeta^v} .
\end{equation}
Due to this relation, we can relate conditional moments with respect to $\tilde{\mathcal{F}}_v$ with moments of $\tilde{X}^i[y, \zeta]$ driven by $\tilde{B}^H$. 

All results of \cref{sec:int_bd} were established for $G\in\{B^H,\tilde B^H\}$ (see \cref{lem:LND_examples} and \cref{subsec:girsanov_prelim}); below we apply them with $G=\tilde B^H$, the Brownian motion $(\tilde W_t)_{t\ge0}$ and its filtration.
In particular, $K^{\tilde B^H}(t,r)=\tilde c_H(t-r)^{H-\frac12}$, and the Cameron--Martin space $\H^{\tilde B^H}$ satisfies \eqref{GGbound}.

In particular, by setting 
\begin{align}\label{eq:def_tilde_psi}
    \tilde{\psi}^i_t \defby \tilde{X}^i_t[y^i, \zeta] - \zeta_t - \tilde{B}^H_t - y^i = \int_0^t b^i(v+r, \tilde{X}^i_r[y^i, \zeta]) \dd r,
\end{align}
we observe that $\tilde{X}^i[y^i,\zeta]-\zeta$ solves \eqref{eqn.XfG} with $G=\tilde{B}^H$, initial value $y^i$ and drift $\tilde{b}^i(r,x)\defby b^i(v+r,x+\zeta_r)$, which is bounded and satisfies $\|\tilde b^i\|_{L_q([0,T-v];L_p)}\le\|b^i\|_{L_qL_p}$ and $\Lambda_{\tilde b^i}\le\Lambda_{b^i}$.
Hence \cref{lem:girsanov_sol} (applicable by \cref{r:boundedisok}) and \eqref{GGbound} imply that, for some $\eps>0$ and all $\lambda \geq 0$,
\begin{align}\label{eq:tilde_psi_novikov}
  \E[e^{\lambda \norm{\tilde{\psi}^i}_{\H^{\tilde{B}^H}}^2}]\le\E[e^{\lambda C\norm{\tilde{\psi}^i}_{\W^{H+\frac{1}{2}+\eps}_2([0, T-v])}^2}] \leq C_1(\lambda, \norm{b^i}_{L_q L_p}) \I_{\{p > 2Hd\}}
  + C_2(\lambda, \Lambda_{b^i}(\eta)) \I_{\{p = 2Hd\}}
\end{align}
with $\eta = \eta(\lambda, T)$, uniformly in $y^i$ and $\zeta$.
Hence, Novikov's condition is satisfied and Girsanov's theorem (\cref{lem:girsanov}) is applicable to $\tilde{B}^H + \tilde{\psi}^i$ and to $\tilde{B}^H+\theta\tilde\psi^1+(1-\theta)\tilde\psi^2$, $\theta\in[0,1]$.  

We set
\begin{align}
  \tilde{V}_t[y^1, y^2, \zeta] &\defby 
  \int_0^t \int_0^1 \nabla b^1(v+r, \theta \tilde{X}^1_r[y^1, \zeta]
  + (1-\theta) \tilde{X}^2_r[y^2, \zeta]) \dd \theta \dd r. \label{eq:def_tilde_V} 
\end{align}

\begin{lemma}\label{lem:tilde_V_var}
  In the setting of \cref{lem:vmo_variation}, we define $\tilde{V} = \tilde{V}[y^1, y^2, \zeta]$ by \eqref{eq:def_tilde_V}. 
  Then, for any $t \leq T-v$ we have
  \begin{multline}\label{eq:tilde_V_var}
    \norm{\norm{\tilde{V}}_{\ell\vr;[0, t]}}_{L_m(\Omega)}
    \leq C_0(\ell, m, T)
    \norm{b^1}_{L_q([v,v+t]; L_p(\mathbb{R}^d))} t^{1-H-\frac{Hd}{p} - \frac{1}{q}} \\
    \times [ C_1(\norm{b^1}_{L_q L_p}, \norm{b^2}_{L_q L_p}) \I_{\{\frac{p}{Hd} > 2\}} 
    + C_2(\Lambda_{b^1}(\eta ), \Lambda_{b^2}(\eta )) \I_{\{\frac{p}{Hd} = 2\}}]
  \end{multline}
  with $\eta  = \eta (m, T)$. The estimate is uniform over $y^1$, $y^2$ and $\zeta$.
\end{lemma}
\begin{proof}
  We proceed as in \cref{l:pathwise.stabitliy}.
  We observe
  \begin{align*}
    \tilde{V}_t[y^1, y^2, \zeta]  = 
    \int_0^1  \tilde{V}^{\theta}_t \dd \theta  
  \end{align*}
  where $\tilde{b}^{1, \theta}(r, x) \defby b^1(v+r, x + \theta y^1 + (1-\theta) y^2 + \zeta_r)$ and (recall that $\tilde{\psi}^i$ is defined by \eqref{eq:def_tilde_psi})
  \begin{align*}
    \tilde{V}^{\theta}_t \defby \int_0^t
    \nabla \tilde{b}^{1, \theta}(r, \tilde{B}^H_r + \theta \tilde{\psi}^1_r
    + (1-\theta) \tilde{\psi}^2_r) \dd r.
  \end{align*}
  By the triangle inequality, it suffices to prove 
  \begin{align}\label{eq:tilde_V_theta_var}
    \norm{\norm{\tilde{V}^{\theta}}_{\ell\vr;[0, t]}}_{L_m(\Omega)} 
    \leq C \norm{b^1}_{L_q([v, v+t]; L_p(\mathbb{R}^d))} t^{1-H-\frac{Hd}{p} - \frac{1}{q}},
  \end{align}
  where the constant $C$ is chosen as in \eqref{eq:tilde_V_var}.

  Setting 
  \begin{align*}
    \tilde{\xi}_t(\tilde{\phi}) \defby 
    \exp \Big( -\int_0^t (K^{\tilde B^H})^{-1} \tilde{\phi}(s) \dd \tilde{W}_s - \frac{1}{2} \int_0^t \abs{(K^{\tilde B^H})^{-1} \tilde{\phi} (s)}^2 \dd s  \Big),
  \end{align*}
  Girsanov's theorem (see the estimate \eqref{genboundstep1}) yields 
  \begin{multline}\label{eq:tilde_V_theta_var_gir}
    \norm{\norm{\tilde{V}^{\theta}}_{\ell\vr;[0, t]}}_{L_m(\Omega)} \\
    \leq \norm[\Big]{\norm[\Big]{\int_0^{\cdot} \nabla \tilde{b}^{1, \theta}(r, \tilde{B}^H_r) \dd r}_{\ell\vr;[0, t]}}_{L_{2m}(\Omega)}
    \E[\tilde{\xi}_t(\theta \tilde{\psi}^1 + (1-\theta) \tilde{\psi}^2)^{-1}]^{\frac{1}{2m}}.
  \end{multline}
    Applying \cref{cor:integral_estimates}(i) with $G=\tilde{B}^H$ and $\eta=0$ (similarly to \eqref{eq:nabla_b_1_B_theta}), we have  
  \begin{equation}\label{eq:tilde_b_1_theta_var}
    \norm[\Big]{\norm[\Big]{\int_0^{\cdot} \nabla \tilde{b}^{1, \theta}(r, \tilde{B}^H_r) \dd r}_{\ell\vr;[0, t]}}_{L_{2m}(\Omega)}
    \leq C(\ell, m) 
    \norm{\tilde{b}^{1, \theta}}_{L_q([0, t]; L_p(\mathbb{R}^d))} t^{1 - H - \frac{Hd}{p } - \frac{1}{q}}.
  \end{equation} 
  As in \eqref{factor2new}, the estimates \eqref{eq:gir_density_mom} and \eqref{eq:tilde_psi_novikov} yield 
  \begin{multline}\label{eq:tilde_xi_minus_1}
    \E[\tilde{\xi}_t(\theta \tilde{\psi}^1 + (1-\theta) \tilde{\psi}^2)^{-1}] \\
    \leq C_1(\norm{b^1}_{L_q L_p}, \norm{b^2}_{L_q L_p}) \I_{\{\frac{p}{dH} > 2\}} 
    + C_2(\Lambda_{b^1}(\eta ), \Lambda_{b^2}(\eta )) \I_{\{\frac{p}{dH} = 2\}}
  \end{multline}
  with $\eta  = \eta (T)$. 
  Noting that 
  $\norm{\tilde{b}^{1, \theta} }_{L_q([0, t]; L_p(\mathbb{R}^d))}
  = \norm{b^1 }_{L_q([v, v+t]; L_p(\mathbb{R}^d))}$,
  the combination of \cref{eq:tilde_V_theta_var_gir,eq:tilde_b_1_theta_var,eq:tilde_xi_minus_1} yields \eqref{eq:tilde_V_var}.
\end{proof}

\begin{proof}[Proof of \cref{lem:vmo_variation}]
    Let $v<t$. In view of \eqref{eq:F_X_i_conditioning},
  \begin{align*}
    \norm{\norm{V}_{\ell\vr;[v, t]}}_{L_m(\Omega) \vert \tilde{\mathcal{F}}_v}
    = \norm{\norm{\tilde{V}[y^1, y^2, \zeta]}_{\ell\vr;[0, t-v]}}_{L_m(\Omega)} \vert_{y^i = X^i_v, \zeta = \zeta^v}.
  \end{align*}
  Since the bound of \cref{lem:tilde_V_var} is uniform in $(y^1,y^2,\zeta)$, the estimate for $V$ follows with $\tilde{\mathcal{F}}_v$ in place of $\mathcal{F}_v$, and hence for $\mathcal{F}_v\subseteq\tilde{\mathcal{F}}_v$ by the tower property, since $\E[\,\cdot\,|\mathcal{F}_v]=\E[\E[\,\cdot\,|\tilde{\mathcal{F}}_v]|\mathcal{F}_v]$ and the bound is deterministic.
  The estimate for $R$ is proved analogously: one replaces $\nabla\tilde b^{1,\theta}$ by $\tilde b^1-\tilde b^2$ (with $\theta=0$), and \cref{cor:integral_estimates}(i) by \cref{lem:integral_along_fbm} with $(\alpha,\rho)=(-1,p)$, which is admissible by \cref{lem.minusoneOK}.
\end{proof}

\subsection{Lifting paths}\label{subsec:lift}
As in the previous section, we assume that $b^1, b^2$ are smooth and 
let $(\mathcal{F}_t)_{t \in [0, T]}$ be the filtration generated by $B^H$.
We define $V$ and $R$ by  \eqref{eq:def_V_and_R}.
\cref{lem:vmo_variation} implies that  
\begin{align*}
  \norm{\norm{V_t - V_s }_{L_m(\Omega) \vert \mathcal{F}_s}}_{L_{\infty}(\Omega)} 
  \leq C(m, T) \norm{b^1}_{L_q([s,t]; L_p(\mathbb{R}^d))} (t-s)^{1-H-\frac{Hd}{p} - \frac{1}{q}},
\end{align*}
and a similar bound for $R$. Thus $V$ and $R$ have vanishing mean oscillation with exponent $\gamma=1 - H - \frac{Hd}{p}-\frac1q$, which can be arbitrarily close to $0$ under \eqref{eq:parameters}.
Although $V$ has finite $\ell$-variation for some $\ell<2$, this variation norm is not exponentially integrable in general (see the beginning of \cref{sec:stability}), so that we work with moments instead.
To estimate $\int X^\Delta\dd V$ in $L_m(\Omega)$ by a sewing argument, the first-order germ $X^\Delta_s(V_t-V_s)$ produces a remainder of order only $2\gamma$, which may be less than $1$.
Therefore, we use germs built from iterated integrals of $V$ up to order $n=\lfloor\gamma^{-1}\rfloor$, whose remainder has order $(n+1)\gamma>1$ (\cref{lem:Z-V-remainder}). That is,  
we have to lift $V$ and $R$ as a rough path of sufficiently large order.
Here we recall notation and very basic properties from the theory of rough paths.

We fix an $\ell \in ((1-H-\frac{Hd}{p}+\frac{Hd}{p\vee 2})^{-1}, 2)$. 
By \cref{lem:vmo_variation}, we know the VMO-type estimates on $\norm{V}_{\ell\vr;[s,t]}$ and $\norm{R}_{\ell\vr;[s,t]}$.
Let $I$ and $\Theta$ be the index sets defined by 
\begin{align}\label{eq:index_set}
  I \defby \{1, 2, \ldots, d\}^2 \tand \Theta \defby \bigcup_{n=1}^{\infty} I^{\set{1, 2, \ldots, n}}.
\end{align}
For $w \in \Theta$, 
we write $\abs{w} = n$ if $w \in I^{\set{1, 2, \ldots, n}}$.
For $w \in \Theta$ and $k < \abs{w}$ we set
\begin{align*}
  w^+_k \defby w \vert_{\{1,2, \ldots, k\}}, \quad w^-_k \defby w \vert_{\{\abs{w}-k+1, \abs{w}-k+2, \ldots, \abs{w}\}}.
\end{align*}
We define the lift $\mathbb{V} = (\mathbb{V}^w)_{w \in \Theta}$ from $V$ as follows. 
If $\abs{w} = 1$ then we set $\mathbb{V}^w_{s,t} \defby V^{w(1)}_t-V^{w(1)}_s$. For $\abs{w}=n$, we define inductively
\begin{align*}
  \mathbb{V}^{w}_{s, t} \defby \int_s^t \mathbb{V}^{w_{n-1}^+}_{s, r} \dd V_r^{w(\abs{w})} 
\end{align*}
as Young's integral \cite{young}.
The lift $\mathbb{V}$ satisfies the algebraic condition called Chen's relation:
\begin{align}\label{eq:chen}
  \mathbb{V}^w_{s, t} - \mathbb{V}^w_{s, u} - \mathbb{V}^w_{u, t}
  = \sum_{k=1}^{n - 1} \mathbb{V}^{w^+_{n-k}}_{s, u} \mathbb{V}^{w^-_k}_{u, t}, \quad \forall s\le u\le t.
\end{align}
It also follows from \cite{lyons1998} that for each $w\in \Theta$, there exists a (deterministic) constant $C(w, \ell)$ such that
\begin{align}\label{est.Vwpathwise}
   |\mathbb{V}^w_{s,t}|\le C(w, \ell)\|V\|_{\ell\vr;[s,t]}^{\abs{w}}, \quad \forall s\le t.
 \end{align}
In the same way, we can construct a joint lift  
\begin{align*}
  (\mathbb{U}^{j w})_{j \in \{1, \ldots, d\}; w \in \Theta}
\end{align*}
of $R$ and $V$ as follows. For
$\abs{w} = 1$ we set 
\begin{align*}
  \mathbb{U}^{j w}_{s, t} \defby \int_s^t R^j_{s, r} \dd V^{w(1)}_r
\end{align*}
and inductively for $\abs{w}>1$,
\begin{align*}
  \mathbb{U}^{j w}_{s, t} \defby \int_s^t \mathbb{U}^{j w^+_{\abs{w} - 1}}_{s, r} \dd V_r^{w(\abs{w})},
\end{align*}
where integrals are defined in Young sense.
We similarly have Chen's identity
\begin{align}\label{eq:chen_U}
  \mathbb{U}^{j w}_{s, t} - \mathbb{U}^{j w}_{s, u} - \mathbb{U}^{j w}_{u, t}
    = (R^j_u - R^j_s) \mathbb{V}^{w}_{u, t} + \sum_{k=1}^{\abs{w} -1} \mathbb{U}_{s, u}^{j w^+_{\abs{w} - k}}
    \mathbb{V}^{w^-_k}_{u, t}\quad \forall s\le u\le t.
\end{align}
 and 
\begin{equation}\label{est.Uwpathwise}
  \abs{\mathbb{U}^{j w}_{s,t}} \leq C(w, \ell) \norm{R}_{\ell\vr;[s,t]} \norm{V}_{\ell\vr;[s,t]}^{\abs{w}},\quad \forall s\le t.
\end{equation}
\begin{lemma}\label{lem:sig_vmo}
  Suppose that the parameters $H, p, q, d$ satisfy \eqref{eq:parameters}.
  Let $b^1, b^2 \in L_q L_p$ be smooth functions with bounded derivatives, and 
  let $\ell \in ((1-H-\frac{Hd}{p}+\frac{Hd}{p\vee 2})^{-1}, 2)$. 
  Let $X^i$ be the strong solution to \Gref{x^i;b^i},
  and we define $V$ and $R$ by \eqref{eq:def_V_and_R}.
    Let $(\mathcal{F}_t)_{t \in [0, T]}$ be the filtration generated by $B^H$.
    If $\frac{p}{Hd} > 2$, for any $m \in (0, \infty)$, $w \in \Theta$ and $s < t$,  
  with $C = C(m, w, \ell, T, \norm{b^1}_{L_q L_p}, \norm{b^2}_{L_q L_p})$ we have
  \begin{align*}
    \norm{\norm{\mathbb{V}^w_{s,t}}_{L_m(\Omega) \vert \mathcal{F}_s}}_{L_{\infty}(\Omega)} 
    \leq C \norm{b^1}_{L_q([s,t]; L_p(\mathbb{R}^d))}^{|w|} (t-s)^{(1-H-\frac{Hd}{p} - \frac{1}{q})\abs{w}}, 
  \end{align*}
  \begin{align*}
    \norm{\norm{\mathbb{U}^{jw}_{s,t} }_{L_m(\Omega) \vert \mathcal{F}_s}}_{L_{\infty}(\Omega)} 
    \leq C \norm{b^1}_{L_q([s,t]; L_p(\mathbb{R}^d))}^{\abs{w}} 
    \norm{b^1 - b^2}_{L_q([s,t]; \mathcal{B}_p^{-1})} (t-s)^{(1-H-\frac{Hd}{p} - \frac{1}{q})(1+\abs{w})}. 
  \end{align*}
  If $\frac{p}{Hd} = 2$, similar estimates hold with some $\eta  = \eta (m, T)$ and 
  \begin{align*}
    C = C(m, w, \ell, T, \Lambda_{b^1}(\eta ), \Lambda_{b^2}(\eta )).
  \end{align*}
\end{lemma}
\begin{proof}
    It easily follows from \cref{lem:vmo_variation}, \eqref{est.Vwpathwise}, \eqref{est.Uwpathwise} and the conditional H\"older inequality.
\end{proof}


\subsection{Moment stability via linear rough differential equations}
In this subsection, we prove \cref{lem:moment_stability}.
Throughout, we assume that $b^1$ and $b^2$ are smooth, and
we fix $m \in [2, \infty)$ and 
\begin{align*}
  \gamma \defby 1 - H - \frac{Hd}{p} - \frac{1}{q}.
\end{align*}
Let $X^i$ be the strong solution to \Gref{x^i;b^i}.
As in the previous subsections, we define $V$ and $R$ by \eqref{eq:def_V_and_R}; we write $V^{ji}$ for the $(j,i)$-entry of the matrix-valued process $V$, that is, $V^{ji}_t=\int_0^t\int_0^1\partial_jb^{1,i}(r,\theta X^1_r+(1-\theta)X^2_r)\,d\theta\,dr$.
(The case $m\in[1,2)$ of \cref{lem:moment_stability} follows from the case $m=2$ by Jensen's inequality.) 
The signatures $(\mathbb{V}^w)_{w}$ and $(\mathbb{U}^{jw})_{j, w}$ are discussed in \cref{subsec:lift}.

Recalling \eqref{eq:Y_linearization}, we set $Y \defby X^{\Delta}$ and $y \defby x^1 - x^2$ so that
\begin{align}\label{eq:Z-equation}
  Y_t^i = y^i + \sum_{j=1}^d \int_0^t Y_r^{j} \dd V_r^{j i} + R^i_t, \quad i=1, 2, \ldots, d.
\end{align}
Our strategy is to employ Gr\"onwall's arguments directly in $L_m(\Omega)$, using the VMO estimates from \cref{lem:sig_vmo}.

For $w \in \Theta$ with $w(i) = (a_1, a_2)$ we write $w(i)(j) \defby a_j$ ($j = 1, 2$).
We denote by $\chain^{j i}_n$ the subset of $\Theta$ such that $w \in \chain^{j i}_n$ if and only if
\begin{align*}
  \begin{gathered}
    \abs{w} \leq n, \quad w(1)(1)=j, \quad w(\abs{w})(2) = i, \\
    w(i)(2) = w(i+1)(1) \quad \forall i \leq \abs{w}-1.
  \end{gathered}
\end{align*}
Iterating equation \eqref{eq:Z-equation} $n$ times shows that, up to a remainder which is of order $n+1$ in $V$, the increment $Y_t^i-Y_s^i$ equals
$R^i_t-R^i_s+\sum_{j=1}^d\sum_{w\in\chain^{ji}_n} Y_s^{j}  \mathbb{V}_{s,t}^{w} +  \sum_{j=1}^d\sum_{w\in\chain^{ji}_{n-1}} \mathbb{U}^{jw}_{s,t}$.
This motivates the following quantity
\begin{multline*}
  \llbracket Y \rrbracket_{s, t}
  \defby
  \sup_{s \leq u < v \leq t} \frac{\norm{Y_v - Y_u}_{L_m(\Omega)}}{(v-u)^{\gamma}} 
  + \sum_{i = 1}^d \sum_{k=1}^{\lfloor \gamma^{-1} \rfloor -1} \sup_{s\leq u < v \leq t} (v-u)^{-(k+1)\gamma} \\
  \times \norm[\Big]{Y^i_v - Y^i_u - (R^i_v - R^i_u) - \sum_{j=1}^d \sum_{w \in \chain_{k}^{j i}} Y_u^j \mathbb{V}^{w}_{u, v}
  - \indic_{\{k \geq 2\}} \sum_{j=1}^d \sum_{w \in \chain_{k-1}^{j i}} \mathbb{U}_{u, v}^{j w}}_{L_m(\Omega)}.
\end{multline*}
Since $b^1$, $b^2$ are smooth with bounded derivatives, $|\mathbb{V}^w_{u,v}|\le C(v-u)^{\abs{w}}$, $|\mathbb{U}^{jw}_{u,v}|\le C(v-u)^{1+\abs{w}}$ and $|Y_v-Y_u|\le C(v-u)$, so that $\llbracket Y \rrbracket_{0, T} < \infty$.
\begin{lemma}\label{lem:Z-V-remainder}
  We set
  \begin{align*}
    A_{s, t}^i \defby \sum_{j=1}^d \Big( \sum_{w \in \chain_n^{j i}} Y_s^j \mathbb{V}_{s, t}^w
    + \sum_{w \in \chain_{n-1}^{j i}} \mathbb{U}_{s, t}^{j w} \Big ),
  \end{align*}
  with $n = \lfloor \gamma^{-1} \rfloor$. 
  Then there is a finite constant $C$ of the form
  \begin{align*}
    C=C(m, T, \norm{b^1}_{L_q L_p}, \norm{b^2}_{L_q L_p})\1_{(p>2Hd)}+ \bar C(m, T, \Lambda_{b^1}(\eta ), \Lambda_{b^2}(\eta ))\1_{(p=2Hd)},
  \end{align*}
  with $\eta  = \eta (m, T)$, such that for every $i \in \{1, \ldots, d\}$ and for every $0 \leq s < t \leq T$,
  we have 
  \begin{equation*}
    \norm[\Big]{\sum_{j=1}^d \int_s^t Y_r^j \dd V^{j i}_r - A_{s, t}^i}_{L_m(\Omega)} 
        \leq C \big(\norm{b^1}_{L_q([s, t]; L_p)} \vee \norm{b^1}_{L_q([s, t]; L_p)}^n\big)
    \llbracket Y \rrbracket_{s, t} (t-s)^{(n+1)\gamma}.
  \end{equation*}
\end{lemma}
\begin{proof}
  Since 
  \begin{align*}
    \abs{\mathbb{V}_{s,t}^w} \leq C(\norm{\nabla b^1}_{L_\infty(\Rd)})(t-s)^{\abs{w}}, \quad 
    \abs{\mathbb{U}_{s,t}^{jw}} \leq C(\norm{b^1}_{L_\infty(\Rd)}, \norm{b^2}_{L_\infty(\Rd)}) (t-s)^{1+\abs{w}},    
  \end{align*}
  we have the following pathwise convergence
  \begin{align*}
    \lim_{\substack{\pi \text{ is a partition of } [s,t], \\ \abs{\pi} \to 0}} \sum_{[u, v] \in \pi} A_{u, v}^i
    &= 
    \lim_{\substack{\pi \text{ is a partition of } [s,t], \\ \abs{\pi} \to 0}} \sum_{[u, v] \in \pi} \sum_{j=1}^d Y_u^j (V^{j i}_{v} - V^{j i}_{u}) \\
    &= \sum_{j=1}^d \int_s^t Y_r^j \dd V_r^{j i}.
  \end{align*}
  We have
  \begin{multline*}
    \delta A_{s, u, t}^i = \sum_{j=1}^d \Big\{- \sum_{w \in \chain_{n}^{j i}} (Y^{j}_{u}-Y^j_s) \mathbb{V}_{u, t}^w \\
    + \sum_{w \in \chain_n^{j i}} Y^j_s (\mathbb{V}^w_{s, t} - \mathbb{V}^w_{s, u} - \mathbb{V}^w_{u, t})
    + \sum_{w \in \chain_{n-1}^{j i}} (\mathbb{U}^{j w}_{s, t} - \mathbb{U}^{j w}_{s, u} - \mathbb{U}^{j w}_{u, t} ) \Big\}.
  \end{multline*}
  Using the Chen's identities  \eqref{eq:chen} and \eqref{eq:chen_U}, we have
  \begin{align*}
    \delta A_{s, u, t}^i
    =& - \sum_{j=1}^d \sum_{w \in \chain_{n-1}^{j i}} \Big( Y_u^j - Y_s^j - (R^j_{u} - R^j_{s}) \\
    &\hspace{2cm} - \sum_{l=1}^d
    \sum_{w_1 \in \chain_{n - \abs{w}}^{l j}} Y^l_s \mathbb{V}_{s, u}^{w_1}
    - \sum_{l=1}^d \sum_{w_2 \in \chain_{n - \abs{w} - 1}^{l j}} \mathbb{U}^{l w_2}_{s, u}
    \Big) \mathbb{V}_{u, t}^w \\
    &- \sum_{j=1}^d \sum_{w \in \chain_{n}^{j i} \setminus \chain_{n-1}^{j i}} (Y^j_{u} - Y^j_{s}) \mathbb{V}_{u, t}^w.
  \end{align*}
    Each summand is of the form $Z\,\mathbb{V}^w_{u,t}$, where $Z$ is $\mathcal{F}_u$-measurable and, by definition of $\llbracket Y \rrbracket$, $\|Z\|_{L_m(\Omega)}\le\llbracket Y \rrbracket_{s,t}(u-s)^{(n-\abs{w}+1)\gamma}$.
  Using $\E|Z\mathbb{V}^w_{u,t}|^m=\E\big[|Z|^m\E[|\mathbb{V}^w_{u,t}|^m|\mathcal{F}_u]\big]$ and the estimate $\norm{\norm{\mathbb{V}_{u, t}^w }_{L_m(\Omega) \vert \mathcal{F}_u}}_{L_{\infty}(\Omega)}
  \leq C \norm{b^1}_{L_q([u,t];L_p)}^{\abs{w}} (t-u)^{\abs{w}\gamma}$
  from \cref{lem:sig_vmo}, we obtain
  \begin{align*}
    \norm{\delta A_{s, u, t}^i}_{L_m(\Omega)} \leq C (\norm{b^1}_{L_q([s,t];L_p)} \vee \norm{b^1}_{L_q([s,t];L_p)}^n)
    \llbracket Y \rrbracket_{s,t} (t-s)^{(n+1)\gamma}.
  \end{align*}
  By our choice of $n$, we have $(n+1)\gamma>1$. The claimed estimate follows from the sewing lemma in the Banach space $L_m(\Omega)$ (see e.g. \cite[Lemma~4.2]{friz14}), applied on $[s,t]$ with the constant $\llbracket Y \rrbracket_{s,t}$, which dominates $\llbracket Y \rrbracket_{u,v}$ for $[u,v]\subseteq[s,t]$.
\end{proof}
\begin{lemma}\label{lem:Z_Gr\"onwall}
  In the above setting,  
  we have for $0 \leq s < t \leq T$
  \begin{align*}
    \norm{Y_t - Y_s}_{L_m(\Omega)} \leq C 
    e^{C (1 \vee \Gamma )^{\frac{n}{\gamma}}}
    \big\{\abs{y} + \norm{b^1 - b^2}_{L_q \mathcal{B}_p^{-1}} \big\} (t-s)^{\gamma},
  \end{align*}
  where $n = \lfloor \gamma^{-1} \rfloor$, $\Gamma \defby \norm{b^1}_{L_qL_p}$ and $y = Y_0 = x^1 - x^2$.
  The constant is chosen as in \cref{lem:Z-V-remainder}.
\end{lemma}
\begin{proof}
  By \eqref{eq:Z-equation} we have
  \begin{multline*}
    Y^i_{s, t} - R^i_{s, t} - \sum_{j=1}^d \Big( \sum_{w \in \chain_n^{j i}} Y_s^j \mathbb{V}_{s, t}^w
    + \sum_{w \in \chain_{n-1}^{j i}} \mathbb{U}_{s, t}^{j w} \Big ) \\
    = \sum_{j=1}^d \int_s^t Y^j_r \dd V^{j i}_r - \sum_{j=1}^d \Big( \sum_{w \in \chain_n^{j i}} Y_s^j \mathbb{V}_{s, t}^w
    + \sum_{w \in \chain_{n-1}^{j i}} \mathbb{U}_{s, t}^{j w} \Big ).
  \end{multline*}
  Note by \cref{lem:sig_vmo} that 
  \begin{align*}
    \norm{Y_s^j \mathbb{V}_{s, t}^w}_{L_m(\Omega)} &\leq C \norm{b^1}_{L_qL_p}^{\abs{w}} \norm{Y_s^j}_{L_m(\Omega)} (t-s)^{\abs{w}\gamma}, \\
      \norm{\mathbb{U}^{j w}_{s, t}}_{L_m(\Omega)} &\leq C \norm{b^1}_{L_qL_p}^{\abs{w}} \norm{b^1 - b^2}_{L_q \mathcal{B}_p^{-1}} (t-s)^{(\abs{w} + 1) \gamma}.
  \end{align*}
    Therefore, if $t-s \leq 1$,  \cref{lem:Z-V-remainder} yields (with $\Gamma = \norm{b^1}_{L_qL_p}$), for every $k\in\{1,\ldots,n-1\}$,
  \begin{multline*}
    \norm[\Big]{Y^i_{s, t} - R^i_{s, t} - \sum_{j=1}^d \Big( \sum_{w \in \chain_k^{j i}} Y_s^j \mathbb{V}_{s, t}^w
    + \sum_{w \in \chain_{k-1}^{j i}} \mathbb{U}_{s, t}^{j w} \Big )}_{L_m(\Omega)} \\
    \leq C (\Gamma \vee \Gamma^n) \Big\{ \llbracket Y \rrbracket_{s, t} (t-s)^{(n+1)\gamma} 
    + \norm{Y_s}_{L_m(\Omega)}(t-s)^{(k+1)\gamma} \\
    + \norm{b^1 - b^2}_{L_q \mathcal{B}_p^{-1}} (t-s)^{(k+1)\gamma} \Big\}.
  \end{multline*}
    Moreover, $\|Y_{s,t}\|_{L_m(\Omega)}\le\|\int_s^tY_r\,dV_r\|_{L_m(\Omega)}+\|R_{s,t}\|_{L_m(\Omega)}$, which is bounded in the same way (the case $k=0$), using \cref{lem:vmo_variation} for $R$.
  Applying these bounds on every subinterval $[u,v]\subseteq[s,t]$, this implies that, with some constant $K \geq 1$ (whose dependence is as in \cref{lem:Z-V-remainder})
  \begin{equation*}
     \llbracket Y \rrbracket_{s,t} 
    \leq K (1 \vee \Gamma^n) \big\{ \llbracket Y \rrbracket_{s, t} (t-s)^{\gamma} + \sup_{r \in [s, t]} \norm{Y_r}_{L_m(\Omega)} 
    + \norm{b^1 - b^2}_{L_q \mathcal{B}_p^{-1}}\big\}.
  \end{equation*}
  Hence, if 
  \begin{align}\label{eq:smallts}
    (t-s)^{\gamma} \leq \frac{1}{2} (K (1 \vee \Gamma^n))^{-1},
  \end{align}
  we have
  \begin{align*}
    \llbracket Y \rrbracket_{s,t}
    \leq 
    2 K (1 \vee \Gamma^n)
    \big\{  \sup_{r \in [s, t]} \norm{Y_r}_{L_m(\Omega)} 
    + \norm{b^1 - b^2}_{L_q \mathcal{B}_p^{-1}}\big\}.
  \end{align*}
  In particular,
  \begin{align}\label{eq:Gr\"onwall}
    \norm{Y_t - Y_s}_{L_m(\Omega)} \leq 2K (1 \vee \Gamma^n) \big\{\sup_{r \in [s, t]}\norm{Y_r}_{L_m(\Omega)} 
    + \norm{b^1 - b^2}_{L_q \mathcal{B}_p^{-1}}\big\} (t-s)^{\gamma}.
  \end{align}
    By the Gr\"onwall lemma (\cref{lem:rough-Gr\"onwall}) applied to $t\mapsto\|Y_t\|_{L_m(\Omega)}$ (in place of $G$) with $w(s,t)\defby t-s$, $\kappa\defby\gamma^{-1}$, $L\defby(2K(1\vee\Gamma^n))^{-1/\gamma}$ and $\Upsilon(s,t)\defby2K(1\vee\Gamma^n)\norm{b^1 - b^2}_{L_q \mathcal{B}_p^{-1}}(t-s)^\gamma$, we obtain
  \begin{align*}
    \sup_{t \in [0, T]} \norm{Y_t}_{L_m(\Omega)}
    \le Ce^{C T (1\vee \Gamma)^{\frac{n}{\gamma}}}\big\{ \abs{y} + \norm{b^1 - b^2}_{L_q \mathcal{B}_p^{-1}}\big\}.
  \end{align*}
  Plugging this estimate into \eqref{eq:Gr\"onwall}, we obtain the claimed estimate, first under the condition \eqref{eq:smallts}, then extended straightforwardly to all $s\le t\le T$. 
\end{proof}
\begin{proof}[Proof of \cref{lem:moment_stability}]
  It follows from \cref{lem:Z_Gr\"onwall}.
\end{proof}

To construct semiflow in the next section, we need the following slight modification of \cref{lem:moment_stability}. 
\begin{corollary}\label{cor:stability_from_s}
  We write $(\mathcal{F}_t)$ for the filtration generated by $B^H$. 
  Let $s \in [0, T)$, $x^{i, s}$ $(i=1,2)$ be $\mathcal{F}_s$-measurable random variables in $\mathbb{R}^d$ and $b^i$ $(i = 1, 2)$ be 
  smooth functions in $[0, T] \times \mathbb{R}^d$ with bounded derivatives.
    Assume \eqref{eq:parameters} and $b^i\in L_qL_p$. Suppose that $X^{i, s}$ solves the SDE  
  \begin{align*}
    X^{i, s}_t = x^{i, s} + \int_s^t b^i(r, X^{i, s}_r) \dd r + B_t^H - B_s^H, \quad \forall t \geq s.
  \end{align*} 
  If $\frac{p}{Hd} > 2$, for any $m \in [1, \infty)$ and $s \leq u < t \leq T$, we have 
  \begin{multline*}
    \norm{X^{1, s}_t - X^{1, s}_u - (X^{2, s}_t - X^{2, s}_u)}_{L_m(\Omega)} \\
    \leq C (\norm{x^{1, s} - x^{2, s}}_{L_m(\Omega)} + \norm{b^1 - b^2}_{L_q \mathcal{B}_p^{-1}}) (t-u)^{1 - H - \frac{Hd}{p } - \frac{1}{q}}.
  \end{multline*}
  with $C = C(m, T, \norm{b^1}_{L_q L_p}, \norm{b^2}_{L_q L_p})$. If $\frac{p}{Hd} = 2$, a similar estimate holds
  with $\eta  = \eta (m, T)$ and $C$ replaced by $Ce^{C(1\vee\|b^1\|_{L_qL_p})^{1/\gamma^2}}$, where $C = C(m, T, \Lambda_{b^1}(\eta ), \Lambda_{b^2}(\eta ))$ and $\gamma=1-H-\frac{Hd}p-\frac1q$, as in \cref{lem:moment_stability}.
\end{corollary}
\begin{proof}
        Since $x^{i,s}$ is $\F_s$-measurable and $b^i$ is smooth, $X^{i,s}$ is measurable with respect to $B^H$, and the claim concerns only the law of $(X^{1,s},X^{2,s})$. Hence, as explained in \cref{subsec:girsanov_prelim}, we may assume that there is a two-sided Brownian motion $\tilde W$ such that \eqref{eq:mandelbrot_van_ness} holds. With $\zeta^s,\bar B^s$ from \eqref{eq:MVN_decomp}, we set $\zeta_t\defby\zeta^s_{t-s}$ and $\bar B_t\defby\bar B^s_{t-s}$ for $t\in[s,T]$,
  so that $B^H_t-B^H_s=\zeta_t+\bar B_t$, the path $\zeta$ is $\tilde{\mathcal{F}}_s$-measurable, and $\bar B_{s+\cdot}$ is a Riemann--Liouville fractional Brownian motion independent of $\tilde{\mathcal{F}}_s$; moreover, $x^{i,s}$ is $\tilde{\mathcal F}_s$-measurable since $\F_s\subseteq\tilde{\mathcal F}_s$.
  Put $\bar{X}^{i, s} \defby X^{i, s} - \zeta$ and $\bar{b}^i(t, x) \defby b^i(t, x+\zeta_t)$.
  We then have
  \begin{align*}
    \bar{X}^{i, s}_t = x^{i, s} + \int_s^t \bar{b}^i(r,\bar{X}_r^{i,s}) \dd r + \bar{B}_{t}, \quad t\in[s,T].
  \end{align*}
  Conditionally on $\tilde{\mathcal{F}}_s$, this is an equation of the form \Gref{x;b} with $B^H$ replaced by a Riemann--Liouville fractional Brownian motion, deterministic initial values and smooth drifts with bounded derivatives satisfying $\|\bar b^i\|_{L_q([s,T];L_p)}\le\|b^i\|_{L_qL_p}$, $\Lambda_{\bar b^i}\le\Lambda_{b^i}$ and $\|\bar b^1-\bar b^2\|_{L_q([s,T];\mathcal{B}^{-1}_p)}\le\|b^1-b^2\|_{L_q\mathcal{B}^{-1}_p}$.
  The proof of \cref{lem:moment_stability} applies verbatim in this setting: for $G=\tilde B^H$ and $v\ge0$, the process $G_{v+\cdot}-\E[G_{v+\cdot}|\mathcal{F}^W_v]$ is again a Riemann--Liouville fractional Brownian motion independent of $\mathcal{F}^W_v$, so that the conditioning argument of \cref{subsec:vmo_estimate} is available, and \cref{sec:int_bd} applies to $G=\tilde B^H$.
  Therefore, almost surely,
  \begin{align*}
    \norm{\bar{X}^{1, s}_t-\bar{X}^{1, s}_u - (\bar{X}^{2, s}_t-\bar{X}^{2, s}_u)}_{L_m(\Omega) \vert \tilde{\mathcal{F}}_s}
    \leq C (\abs{x^{1, s} - x^{2, s}} + \norm{b^1 - b^2}_{L_q \mathcal{B}_p^{-1}})(t-u)^{1 - H - \frac{Hd}{p } - \frac{1}{q}}
  \end{align*}
  with a deterministic constant $C$ as in the statement.
  Since $\bar{X}^{1, s} - \bar{X}^{2, s} = X^{1,s} - X^{2, s}$, taking the $L_m(\Omega)$-norm yields the claimed estimate.  
\end{proof}
\section{Semiflow and path-by-path uniqueness} 
\label{sec:flow_and_path_by_path_uniqueness}
    In this section, we construct the semiflow of \cref{thm:main_semiflow} (\cref{subsec:semiflow}) and deduce path-by-path uniqueness from it (\cref{subsec:pbp}).
  In \cite{Shaposhnikov:2016aa}, Shaposhnikov showed that path-by-path uniqueness can be obtained from the existence of an almost Lipschitz semiflow. This line of arguments is later polished and streamlined in \cite{galeatigerencser22,anzeletti2023pathbypath}, and we follow \cite[Section~5 and Appendix~A]{anzeletti2023pathbypath}.
  The two ingredients are the regularizing property of fractional Brownian paths (\cref{lem:fbm_regularizing}), which says that almost every path of $B^H$ is $(1-\varepsilon,\gamma)$-regularizing for $b$ and $|b|$ in the sense of \cite[Definition~A.1]{anzeletti2023pathbypath}, and the stability estimates in moment norms (\cref{cor:stability_from_s}).
  In contrast to \cite{anzeletti2023pathbypath}, we construct the semiflow as a limit of the flows of equations with smooth drifts, which avoids the use of strong solutions with random initial data.
\subsection{Construction of the semiflow}\label{subsec:semiflow}

  \begin{remark}\label{rem:flow}
    A semiflow is called a \emph{flow} if the map $x \mapsto \Phi^{s, x}_t$ is a homeomorphism of $\mathbb{R}^d$.
    It is possible to construct a flow for \eqref{eq:main-sde} under \eqref{eq:parameters}, following for instance the arguments in \cite{galeatigerencser22}. 
    However, we will  not  address flow property in our work for conciseness.
  \end{remark}
  %


  The following pathwise version of \cref{prop.invariance} identifies limits of solutions of equations with smooth drifts.
  \begin{lemma}\label{lem:pathwise_invariance}
    Assume \eqref{eq:parameters}, let $b\in L_qL_p$ and let $(b^n)$ be a subsequence of the mollifying sequence of $b$ from \cref{lem:mollifying}. Then, possibly after passing to a further subsequence of $(b^n)$, there exists an event $\Omega_{b,(b^n),T}$ of full probability with the following property.
    Let $\omega\in\Omega_{b,(b^n),T}$, $s\in[0,T)$ and let $y^n,y\colon[s,T]\to\Rd$ be continuous paths such that $y^n\to y$ uniformly on $[s,T]$ and
    \begin{equation*}
      y^n_t=y^n_s+\int_s^tb^n(r,y^n_r)\dd r+B^H_t(\omega)-B^H_s(\omega),\quad t\in[s,T],\ n\ge1.
    \end{equation*}
        Then $\int_s^T|b(r,y_r)|\dd r<\infty$ and $y_t=y_s+\int_s^tb(r,y_r)\dd r+B^H_t(\omega)-B^H_s(\omega)$ for $t\in[s,T]$.
  \end{lemma}
  \begin{proof}
    Let $\Omega_{b,(b^n),T}$ be the event constructed in the proof of \cref{prop.invariance}, for the subsequence chosen there; it depends only on $b$, $(b^n)$ and $T$.
    The proof of \cref{prop.invariance} applies verbatim, with $[0,T]$ replaced by $[s,T]$, $\psi^n\defby y^n-B^H(\omega)$ and $\psi\defby y-B^H(\omega)$: it uses only $\omega\in\Omega_{b,(b^n),T}$, the uniform convergence $\psi^n\to\psi$, the a priori bound \eqref{eq:a_priori_bound} on $[s,T]$ (which gives $\liminf_n\norm{\psi^n}_{1\vr;[s,T]}<\infty$ and hence an $\omega$-dependent subsequence with bounded $1$-variation), \cref{cor:continuity_along_B} and \eqref{est.intfome}.
  \end{proof}

\begin{proof}[Proof of \cref{thm:main_semiflow}]
    Put $\gamma\defby1-H-\frac{Hd}{p}-\frac1q>0$.
  Let $(b^n)_{n=1}^{\infty}$ be the mollifying sequence of $b$ from \cref{lem:mollifying}. By \cref{lem:mollifying} and \cref{rem:smooth_drift}, $b^n$ is a smooth function with bounded derivatives, $b^n \to b$ in $L_q \mathcal{B}_p^{-1}$, and $\|b^n\|_{L_qL_p}$ and $\Lambda_{b^n}(\eta)$ are bounded uniformly in $n$ for each $\eta>0$, so that the constants below can be chosen independently of $n$.
  For $(s,x)\in[0,T)\times\Rd$, let $\Psi^{s, x, n}$ be the solution to
  \begin{equation*}
    \Psi^{s, x, n}_t = x + \int_s^t b^n(r, \Psi^{s, x, n}_r) \dd r + B_t^H - B_s^H,\quad t\in[s,T],
  \end{equation*}
  and set $\Psi^{s,x,n}_t\defby x$ for $t<s$. Since $b^n$ is Lipschitz in $x$ uniformly in $t$, this equation is solved path by path: for every $\omega$, $(s,t,x)\mapsto\Psi^{s,x,n}_t(\omega)$ is continuous, satisfies the semiflow property $\Psi^{s,x,n}_t=\Psi^{u,\Psi^{s,x,n}_u,n}_t$ for all $s\le u\le t$ and $x$, and $\Psi^{s,x,n}_t$ is a measurable function of $(B^H_r-B^H_s)_{r\in[s,t]}$.
  Fix $m\ge2$, to be taken large in Step~2. Constants $C$ below depend on $m,T$ and on $b$ through $\|b\|_{L_qL_p}$ and $\Lambda_b(\eta)$, but not on $n,s,t,x$.

  \textbf{Step 1: moment estimates.} By \cref{rem:int_f_sol} (which bounds the drift increments by $C|t-t'|^{1-\frac{Hd}p-\frac1q}=C|t-t'|^{\gamma+H}$)\footnote{Here \cref{rem:int_f_sol} is applied to the process $u\mapsto\Psi^{s,x,n}_{s+u}$, which is a strong solution driven by the fractional Brownian motion $B^H_{s+\cdot}-B^H_s$, with drift $b^n(s+\cdot,\cdot)$.} and since $\|B^H_t-B^H_{t'}\|_{L_m(\Omega)}\le C|t-t'|^H$, we have for $s\le t'\le t$
  \begin{equation}\label{eq:Psi_t_t_prime}
    \norm{\Psi_t^{s, x, n} - \Psi_{t'}^{s, x, n}}_{L_m(\Omega)}
        \leq C (\abs{t-t'}^{\gamma+H}+\abs{t-t'}^{H})\le C\abs{t-t'}^{H},
  \end{equation}
  By \cref{cor:stability_from_s} (applied with $u=s$ and $b^1=b^2=b^n$; note that $\Psi^{s,x,n}_t-\Psi^{s,x',n}_t=(x-x')+\big[(\Psi^{s,x,n}_t-x)-(\Psi^{s,x',n}_t-x')\big]$), we have
  \begin{equation}\label{eq:Psi_x_x_prime}
    \norm{\Psi^{s, x, n}_t - \Psi^{s, x', n}_t}_{L_m(\Omega)}
    \leq C \abs{x - x'},
  \end{equation}
  and for $s' < s\le t$, using the semiflow property and \cref{cor:stability_from_s} with the $\mathcal{F}_s$-measurable initial value $\Psi^{s',x,n}_s$, followed by \eqref{eq:Psi_t_t_prime},
  \begin{equation}\label{eq:Psi_s_s_prime}
    \norm{\Psi^{s, x, n}_t - \Psi^{s', x, n}_t}_{L_m(\Omega)}
    = \norm{\Psi^{s, x, n}_t - \Psi^{s, \Psi_s^{s', x, n}, n}_t}_{L_m(\Omega)}
    \leq C \norm{x - \Psi_s^{s', x, n}}_{L_m(\Omega)}
    \leq C \abs{s - s'}^{H}.
  \end{equation}
  The cases where $t<s$ or $t'<s$ reduce to these by the convention $\Psi^{s,x,n}_t=x$ for $t<s$.
  Finally, by \cref{cor:stability_from_s} with $b^1=b^n$, $b^2=b^{n'}$,
  \begin{equation}\label{eq:Psi_m_n}
    \norm{\Psi^{s, x, n}_t - \Psi^{s, x, n'}_t}_{L_m(\Omega)}
    \leq C \norm{b^n - b^{n'}}_{L_q \mathcal{B}_p^{-1}}.
  \end{equation}
  Write $z=(s,t,x)$, $z'=(s',t',x')$, $\rho(z,z')\defby\abs{s-s'}^H + \abs{t-t'}^{H} + \abs{x-x'}$ and $D^{n,n'}(z)\defby\Psi^{s,x,n}_t-\Psi^{s,x,n'}_t$. By \cref{eq:Psi_t_t_prime,eq:Psi_x_x_prime,eq:Psi_s_s_prime}, $\|\Psi^{n}(z)-\Psi^{n}(z')\|_{L_m(\Omega)}\le C\rho(z,z')$, and by \eqref{eq:Psi_m_n} and the inequality $\min(a,b)\le a^{\varepsilon}b^{1-\varepsilon}$, for every $\varepsilon\in(0,1)$,
  \begin{equation*}
    \norm{D^{n,n'}(z)-D^{n,n'}(z')}_{L_m(\Omega)}
    \leq C \norm{b^n - b^{n'}}_{L_q \mathcal{B}_p^{-1}}^{\varepsilon} \rho(z,z')^{1 - \varepsilon}.
  \end{equation*}

  \textbf{Step 2: construction of $\Phi$.} Let $L\in\N$ and $Z_L\defby[0,T]^2\times\{|x|\le L\}$. Since $m$ can be chosen arbitrarily large, the Kolmogorov continuity theorem (in the form of the Garsia--Rodemich--Rumsey inequality) and the previous estimates yield, for every $\varepsilon\in(0,1)$,
  \begin{equation*}
    \Big\|\sup_{z\in Z_L}|D^{n,n'}(z)|\Big\|_{L_m(\Omega)}\le C_L\norm{b^n - b^{n'}}_{L_q \mathcal{B}_p^{-1}}^{\varepsilon},\qquad
    \Big\|\sup_{z\neq z'\in Z_L}\frac{|\Psi^n(z)-\Psi^n(z')|}{\rho(z,z')^{1-\varepsilon}}\Big\|_{L_m(\Omega)}\le C_{L,\varepsilon},
  \end{equation*}
  with constants independent of $n,n'$. Since $(b^n)$ converges in $L_q\B^{-1}_p$, passing to a further subsequence, we find an event $\Omega_1$ of full probability on which $\Psi^n(\omega)$ converges uniformly on every $Z_L$ to a continuous limit $\Phi(\omega)$, and on which, by Fatou's lemma, $\Phi(\omega)$ satisfies the displayed H\"older estimate for every rational $\varepsilon\in(0,1)$ and every $K\in\N$; this suffices, since smaller $\varepsilon$ gives a stronger estimate on bounded sets, and since $\rho^{1-\varepsilon}\le C(|s-s'|^{H(1-\varepsilon)}+|t-t'|^{H(1-\varepsilon)}+|x-x'|^{1-\varepsilon})$ it yields \eqref{eq:semiflow_holder} for every $\varepsilon\in(0,H)$ (take $\varepsilon'=\varepsilon/H$). On $\Omega\setminus\Omega_1$ we set $\Phi\defby0$.
  Then $\Phi$ is measurable, and $\Phi^{s,x}_t$ is measurable with respect to the completion of $\sigma(B^H_r-B^H_s:r\in[s,t])$, as an almost sure limit of such random variables.

  \textbf{Step 3: semiflow property.} Let $\omega\in\Omega_1$ and $s\le u\le t$, $x\in\Rd$. Since $\Psi^{s,x,n}_u(\omega)\to\Phi^{s,x}_u(\omega)$, and $\Psi^n(\omega)\to\Phi(\omega)$ uniformly on compact sets with $\Phi(\omega)$ continuous, we may pass to the limit in $\Psi^{s,x,n}_t(\omega)=\Psi^{u,\Psi^{s,x,n}_u(\omega),n}_t(\omega)$ to get $\Phi^{s,x}_t(\omega)=\Phi^{u,\Phi^{s,x}_u(\omega)}_t(\omega)$. This is \cref{def:semiflow}(ii).

  \textbf{Step 4: $\Phi$ solves the equation.} We apply \cref{lem:pathwise_invariance} to the subsequence of $(b^n)$ obtained in Step~2, and pass to the further subsequence provided there; along it, $\Psi^n(\omega)\to\Phi(\omega)$ uniformly on compact sets for every $\omega\in\Omega_1$. Let $\omega\in\Omega_1\cap\Omega_{b,(b^n),T}$, with $\Omega_{b,(b^n),T}$ from \cref{lem:pathwise_invariance}, and let $(s,x)\in[0,T)\times\Rd$. The paths $y^n\defby\Psi^{s,x,n}(\omega)|_{[s,T]}$ solve the equations in \cref{lem:pathwise_invariance} and converge uniformly to $y\defby\Phi^{s,x}(\omega)|_{[s,T]}$. Hence, by \cref{lem:pathwise_invariance}, $\int_s^T|b(r,\Phi^{s,x}_r(\omega))|\dd r<\infty$ and
  \begin{equation*}
    \Phi^{s, x}_t(\omega) = x + \int_s^t b(r, \Phi^{s, x}_r(\omega)) \dd r + B_t^H(\omega) - B_s^H(\omega),\quad t\in[s,T].
  \end{equation*}
  Since the event $\Omega_1\cap\Omega_{b,(b^n),T}$ does not depend on $(s,x)$, this proves \cref{def:semiflow}(iii).
\end{proof}

\subsection{Path-by-path uniqueness}\label{subsec:pbp}
We now identify any solution along a typical path of $B^H$ with the semiflow, following \cite[Proposition~5.2 and the proof of Theorem~2.2]{anzeletti2023pathbypath}; see also \cite[Theorem~A.2]{anzeletti2023pathbypath} for the abstract formulation.
Throughout this subsection, $\Phi$ is the semiflow from \cref{thm:main_semiflow}, $\Omega_\Phi$ denotes the event of full probability on which $\Phi$ is a semiflow satisfying \eqref{eq:semiflow_holder}, and we fix
\begin{equation}\label{eq:pbp_exponents}
  \gamma\in\Big(0,1-H-\frac{Hd}{p}-\frac1q\Big),\qquad \varepsilon\in(0,\gamma),\qquad \beta\defby\min\Big(1-\varepsilon+\gamma,\frac\gamma\varepsilon\Big)>1.
\end{equation}
For a closed ball $K\subset\Rd$ centred at $0$, we let $\Xi_{T,K}$ and $\Omega_{f,T,K}$, $f\in\{b,|b|\}$, be as in \cref{lem:fbm_regularizing}, with these $\gamma$ and $\varepsilon$.

\begin{proposition}[A priori comparison with the semiflow]\label{prop:comparison_semiflow}
  Assume \eqref{eq:parameters} and let $b\in L_qL_p$. Let $K\subset\Rd$ be a closed ball centred at $0$, and let $\omega\in\Omega_\Phi\cap\Omega_{b,T,K}\cap\Omega_{|b|,T,K}$.
  Let $Y$ be a solution to \eqref{fSDE} along $B^H(\omega)$, that is, $\int_0^T|b(r,Y_r)|\dd r<\infty$ and $Y_t=Y_0+\int_0^tb(r,Y_r)\dd r+B^H_t(\omega)$ for $t\in[0,T]$.
  Suppose that $Y_t-B^H_t(\omega)\in K$ and $\Phi^{s,Y_s}_t(\omega)-B^H_t(\omega)\in K$ for every $0\le s\le t\le T$.
  Then there exists a finite constant $C=C(\gamma,\varepsilon,\Xi_{T,K}(b)(\omega),\Xi_{T,K}(|b|)(\omega))$ such that, with the control $w(s,t)\defby\int_s^t|b(r,Y_r)|\dd r+(t-s)$,
  \begin{equation}\label{est.apriY}
    |Y_t-\Phi^{s,Y_s}_t(\omega)|\le C\big(w(s,t)^{1-\varepsilon+\gamma}+w(s,t)^{\gamma/\varepsilon}\big)\quad\forall\,0\le s\le t\le T.
  \end{equation}
\end{proposition}
\begin{proof}
  We omit the dependence on $\omega$. Fix $0\le s\le t\le T$ and put $\psi\defby Y-B^H$, $\xi\defby\Phi^{s,Y_s}-B^H$ on $[s,t]$. Then $\psi_s=\xi_s$, both paths take values in $K$, and
  \begin{equation*}
    \norm{\psi}_{1\vr;[s,t]}\le\int_s^t|b|(r,Y_r)\dd r,\qquad \norm{\xi}_{1\vr;[s,t]}\le\int_s^t|b|(r,\Phi^{s,Y_s}_r)\dd r<\infty
  \end{equation*}
  by \cref{def:semiflow}(iii). We write $I_Y\defby\int_s^t|b|(r,Y_r)\dd r$, $I_\Phi\defby\int_s^t|b|(r,\Phi^{s,Y_s}_r)\dd r$ and $\Delta\defby|I_\Phi-I_Y|$.

  \textit{Step 1.} Writing $I_\Phi-I_Y=\int_s^t[|b|(r,B^H_r+\xi_r)-|b|(r,B^H_r+\xi_s)]\dd r-\int_s^t[|b|(r,B^H_r+\psi_r)-|b|(r,B^H_r+\psi_s)]\dd r$ and applying \eqref{est.young} with $f=|b|$ to both terms, we get
  \begin{equation*}
    \Delta\le C\,\Xi_{T,K}(|b|)\big(I_\Phi^{1-\varepsilon}+I_Y^{1-\varepsilon}\big)(t-s)^\gamma\le C\,\Xi_{T,K}(|b|)\big(\Delta^{1-\varepsilon}+2I_Y^{1-\varepsilon}\big)(t-s)^\gamma,
  \end{equation*}
  where we used $I_\Phi\le\Delta+I_Y$ and the subadditivity of $a\mapsto a^{1-\varepsilon}$. By Young's inequality, $C\,\Xi_{T,K}(|b|)\Delta^{1-\varepsilon}(t-s)^\gamma\le\frac12\Delta+C'(\Xi_{T,K}(|b|)(t-s)^\gamma)^{1/\varepsilon}$, and hence
  \begin{equation}\label{est.absb}
    \Delta\le C\big(I_Y^{1-\varepsilon}(t-s)^\gamma+(t-s)^{\gamma/\varepsilon}\big),
  \end{equation}
  with $C$ depending on $\Xi_{T,K}(|b|)$, $\gamma$ and $\varepsilon$.

  \textit{Step 2.} Using the equations for $Y$ and $\Phi^{s,Y_s}$ and $\psi_s=\xi_s$, we have
  \begin{equation*}
    |\Phi^{s,Y_s}_t-Y_t|\le\Big|\int_s^t[b(r,B^H_r+\xi_r)-b(r,B^H_r+\xi_s)]\dd r\Big|+\Big|\int_s^t[b(r,B^H_r+\psi_r)-b(r,B^H_r+\psi_s)]\dd r\Big|.
  \end{equation*}
  By \eqref{est.young} with $f=b$, $I_\Phi\le\Delta+I_Y$ and \eqref{est.absb},
  \begin{equation*}
    |\Phi^{s,Y_s}_t-Y_t|\le C\,\Xi_{T,K}(b)\big(I_\Phi^{1-\varepsilon}+I_Y^{1-\varepsilon}\big)(t-s)^\gamma\le C\big(I_Y^{1-\varepsilon}(t-s)^\gamma+(t-s)^{\gamma/\varepsilon}\big),
  \end{equation*}
  where, in the last step, we used $I_\Phi^{1-\varepsilon}\le I_Y^{1-\varepsilon}+\Delta^{1-\varepsilon}$, the bound $\Delta^{1-\varepsilon}(t-s)^\gamma\le C(I_Y^{(1-\varepsilon)^2}(t-s)^{\gamma(2-\varepsilon)}+(t-s)^{\gamma/\varepsilon})$, which follows from \eqref{est.absb}, and $I_Y^{(1-\varepsilon)^2}(t-s)^{\gamma(2-\varepsilon)}\le I_Y^{1-\varepsilon}(t-s)^\gamma+(t-s)^{\gamma/\varepsilon}$, which is Young's inequality with exponents $\frac1{1-\varepsilon}$ and $\frac1\varepsilon$.
  Since $I_Y\le w(s,t)$ and $t-s\le w(s,t)$, this yields \eqref{est.apriY}.
\end{proof}

\begin{proposition}[Path-by-path uniqueness]\label{prop:pbp_uniqueness}
  Assume \eqref{eq:parameters} and let $b\in L_qL_p$. Let $K_N\defby\{x\in\Rd:|x|\le N\}$ and
  \begin{equation*}
    \Omega_{b,T}^*\defby\Omega_\Phi\cap\bigcap_{N\in\N}\big(\Omega_{b,T,K_N}\cap\Omega_{|b|,T,K_N}\big),
  \end{equation*}
  which is an event of full probability depending only on $b$ and $T$. Then for every $\omega\in\Omega_{b,T}^*$, every solution $Y$ to \eqref{fSDE} along $B^H(\omega)$ coincides with $\Phi^{0,Y_0}(\omega)$ on $[0,T]$.
  In particular, path-by-path uniqueness holds for \eqref{eq:main-sde}, and for every $\omega\in\Omega_{b,T}^*$ and $x\in\Rd$, $\Phi^{0,x}(\omega)$ is the unique solution to \eqref{fSDE} along $B^H(\omega)$ with initial value $x$.
\end{proposition}
\begin{proof}
  Let $\omega\in\Omega_{b,T}^*$ and let $Y$ be a solution along $B^H(\omega)$; we omit the dependence on $\omega$.
  By continuity of $Y$, $B^H$ and of $(s,t,x)\mapsto\Phi^{s,x}_t$, the sets $\{Y_t-B^H_t,\ \Phi^{s,Y_s}_t-B^H_t: 0\le s\le t\le T\}$ and $\{Y_t,\ \Phi^{s,Y_s}_t: 0\le s\le t\le T\}$ are bounded. Hence there is $N$ such that all these points belong to $K_N$, and \cref{prop:comparison_semiflow} applies with $K=K_N$.
  Fix $\kappa\in(\beta^{-1},1)$, which is possible since $\beta>1$, and let $C_\kappa=C(\omega,\kappa,N)$ be the H\"older constant from \eqref{eq:semiflow_holder} on $K_N$.
  Fix $\tau\in(0,T]$ and put $F(t)\defby\Phi^{t,Y_t}_\tau$, $t\in[0,\tau]$. For $0\le s\le t\le\tau$, the semiflow property gives $\Phi^{s,Y_s}_\tau=\Phi^{t,\Phi^{s,Y_s}_t}_\tau$, and hence, by \eqref{eq:semiflow_holder} and \eqref{est.apriY},
  \begin{equation*}
    |F(t)-F(s)|=|\Phi^{t,Y_t}_\tau-\Phi^{t,\Phi^{s,Y_s}_t}_\tau|\le C_\kappa|Y_t-\Phi^{s,Y_s}_t|^\kappa\le C\big(w(s,t)^{\kappa(1-\varepsilon+\gamma)}+w(s,t)^{\kappa\gamma/\varepsilon}\big).
  \end{equation*}
  Both exponents are at least $\kappa\beta>1$. For a partition $\{t_i\}$ of $[0,\tau]$ with mesh tending to $0$, we have $\max_iw(t_{i-1},t_i)\to0$ since $w$ is a control, and for $\theta>1$, $\sum_iw(t_{i-1},t_i)^\theta\le\max_iw(t_{i-1},t_i)^{\theta-1}w(0,\tau)\to0$. Hence $F(\tau)=F(0)$, that is, $Y_\tau=\Phi^{\tau,Y_\tau}_\tau=\Phi^{0,Y_0}_\tau$. Since $\tau$ is arbitrary, $Y=\Phi^{0,Y_0}$ on $[0,T]$.
  The last assertion follows since $\Phi^{0,x}(\omega)$ is a solution by \cref{def:semiflow}(iii).
\end{proof}

\section{Proofs of the main results}\label{sec:proof_main_results}

\begin{proof}[Proof of \cref{thm:main_uniqueness}]
  (i) Strong existence is proven in \cref{l:strongsol}. Alternatively, $\Phi^{0,x}$ from \cref{thm:main_semiflow} is a strong solution, since $\Phi^{0,x}_t$ is measurable with respect to the completion of $\sigma(B^H_r:r\le t)$.

  (ii) This is \cref{prop:pbp_uniqueness}: for every $\omega$ in the event $\Omega^*_{b,T}$ of full probability, which does not depend on the initial value, and every $x\in\Rd$, the path $\Phi^{0,x}(\omega)$ is the unique solution to \eqref{fSDE} along $B^H(\omega)$.
\end{proof}

By \cref{thm:main_uniqueness}, for each $x\in\Rd$ and $b\in L_qL_p$ the solution to \Gref{x;b} is unique, and it coincides almost surely with the strong solution constructed in \cref{l:strongsol}. This justifies the phrase ``the solution'' in \cref{thm:main_stability}. Note that the proof of \cref{thm:main_uniqueness} does not rely on \cref{thm:main_stability}.

\begin{proof}[Proof of \cref{thm:main_stability}]
  Let $X^i$ be the solution to \Gref{x^i;b^i}, $i=1,2$, which is the strong solution from \cref{l:strongsol}. Put $\gamma\defby1-H-\frac{Hd}p-\frac1q$ and fix $\ell\in\big((1-H-\frac{Hd}{p}+\frac{Hd}{p\vee2})^{-1},2\big)$, which is possible by \cref{lem.minusoneOK}.
  In this proof, $C$ denotes a finite constant whose dependence is as in \eqref{eq:C_dep} if $p>2Hd$ and as in \eqref{eq:C_dep_2} if $p=2Hd$, for a suitable $\eta=\eta(d,H,p,q,T,m)$, and which may change from line to line.
  Since $\|Z\|_{L_m(\Omega)}\le\|Z\|_{L_{m\vee 2}(\Omega)}$, it suffices to consider $m\ge2$.

  \textbf{Step 1: exponential moments of the drifts.} By \eqref{nboundsolquant} and \eqref{GGbound}, for every $\lambda>0$,
  \begin{equation*}
    M(\lambda)\defby\max_{i=1,2}\E\exp\Big(\lambda\Big\|\int_0^\cdot b^i(r,X^i_r)\,dr\Big\|^2_{\H}\Big)\le C(\lambda),
  \end{equation*}
  where $C(\lambda)$ depends on $\lambda$ and otherwise as described above (and is finite for $p=\infty$ by \eqref{novikovbound2}). In particular, $X^1$ and $X^2$ satisfy \eqref{con.expXi}, and the constant $C_0$ in \eqref{Vbound11}--\eqref{Vbound22} (which depends on $M(12)$) is of the form $C$.

  \textbf{Step 2: pathwise stability.} By \cref{l:pathwise.stabitliy}, almost surely
  \begin{equation*}
    \sup_{t\in[0,T]}|X^1_t-X^2_t|\le 4e^{c\|V\|^\ell_{\ell\var;[0,T]}}\big(|x^1-x^2|+\|R\|_{\ell\var;[0,T]}\big),
  \end{equation*}
  with $V$, $R$ given by \eqref{eq:def_V_and_R} and $c=c(\ell)$. We set
  \begin{equation*}
    A\defby\log4+c\|V\|^\ell_{\ell\var;[0,T]},\qquad D\defby\|R\|_{\ell\var;[0,T]},
  \end{equation*}
  so that \eqref{eq:X_1_X_2_pw_bd} holds. By \eqref{Vbound11} with Besov parameter $0$ (in place of $\eta$) and $m\ell$ in place of $m$, and since $\|f\|_{L_q\B^0_p}\le C\|f\|_{L_qL_p}$,
  \begin{equation*}
    \|A\|_{L_m(\Omega)}\le\log4+c\,\big\|\|V\|_{\ell\var;[0,T]}\big\|^\ell_{L_{m\ell}(\Omega)}\le \log 4+C\|b^1\|_{L_qL_p}^\ell T^{\ell\gamma}\le C(1+\|b^1\|_{L_qL_p})^2,
  \end{equation*}
  where we used $\ell<2$. By \eqref{Vbound22} with Besov parameter $0$, $\|D\|_{L_m(\Omega)}\le C\|b^1-b^2\|_{L_q\B^{-1}_p}T^{\gamma}$. This proves part (i).

  \textbf{Step 3: stability in moment norms.} Let $(b^{i,n})_{n\ge1}$ be the mollifying sequence of $b^i$ from \cref{lem:mollifying}, and let $X^{i,n}$ be the strong solution to \Gref{x^i;b^{i,n}}, which exists and is unique by \cref{rem:smooth_drift}.
  By \cref{lem:moment_stability} applied to $(b^{1,n},b^{2,n})$, for $0\le s<t\le T$,
  \begin{equation*}
    \norm{X^{1,n}_t - X^{1,n}_s - (X^{2,n}_t - X^{2,n}_s)}_{L_m(\Omega)}
    \leq C_n e^{C_n(1\vee\norm{b^{1,n}}_{L_qL_p})^{1/{\gamma^{2}} }} (\abs{x^1 - x^2} + \norm{b^{1,n}-b^{2,n}}_{L_q \mathcal{B}_p^{-1}})
    (t-s)^{\gamma},
  \end{equation*}
  where $C_n$ depends increasingly on $\|b^{i,n}\|_{L_qL_p}$ if $p>2Hd$, and on $\Lambda_{b^{i,n}}(\eta)$ if $p=2Hd$, $i=1,2$ (see \cref{rem:dependence_of_const}). By \cref{lem:mollifying}(ii),(iii), $\|b^{i,n}\|_{L_qL_p}\le C\|b^i\|_{L_qL_p}$ and $\Lambda_{b^{i,n}}(\eta)\le2\Lambda_{b^i}(\eta/2)$; hence $C_n$ and the exponential factor are bounded uniformly in $n$ by constants of the form $C$ and $e^{C(1\vee\|b^1\|_{L_qL_p})^{1/\gamma^2}}$, respectively (with $\eta/2$ in place of $\eta$). If $p>2Hd$, the exponential factor is absorbed into $C$.
  Furthermore, by \cref{lem:mollifying}(i), $b^{i,n}\to b^i$ in $L_q\B^{-1}_p$, so that $\norm{b^{1,n}-b^{2,n}}_{L_q \mathcal{B}_p^{-1}}\to\norm{b^{1}-b^{2}}_{L_q \mathcal{B}_p^{-1}}$.
  Finally, the proof of \cref{l:strongsol} shows that $X^{i,n}\to X^i$ uniformly on $[0,T]$ in probability, where we used that the limit there is a strong solution to \Gref{x^i;b^i} and hence coincides with $X^i$ by \cref{thm:main_uniqueness}.
  Passing to an almost surely convergent subsequence and applying Fatou's lemma, we obtain \eqref{eq:moment_stab} for $p>2Hd$ and its modification stated in part (ii) for $p=2Hd$.
\end{proof}

\appendix
\section{Approximation of the drift}\label{app:aux}

\begin{lemma}[Mollifying sequence]\label{lem:mollifying}
  Let $f$ be a function in $L_q([0,T],L_p(\Rd))$ with $p,q\in[1,\infty]$. For each integer $n\ge1$, define
  \begin{equation*}
    f^n(t,x) \defby \I_{\{\|f(t,\cdot)\|_{L_p(\Rd)}\le n\}}\,P_{1/n}f(t,\cdot)(x),\quad (t,x)\in[0,T]\times\Rd.
  \end{equation*}
  Then the following holds.
  \begin{enumerate}[label={(\roman*)}, font=\normalfont, leftmargin=*]
        \item[(0)] Each $f^n$ is bounded, and smooth in $x$ with all spatial derivatives bounded uniformly in $(t,x)\in[0,T]\times\Rd$.
    \item $\lim_n f^n=f$ in $L_q([0,T],\B^{-\varepsilon}_p(\Rd))$ for every $\varepsilon>0$;
    \item $\sup_n\|f^n\|_{L_q([0,T],\B^{-\varepsilon}_p(\Rd))}\le C\sup_n\|f^n\|_{L_q([0,T],L_p(\Rd))}\le C\|f\|_{L_q([0,T],L_p(\Rd))}$;
    \item $\Lambda_{f^n}(2\eta)\le2\Lambda_f(\eta)$ for every $n\ge1$ and $\eta>0$ (with the convention of \cref{rem:Lambda_q_infty} if $q=\infty$).
  \end{enumerate}
\end{lemma}
\begin{proof}
  (0) For any multi-index $k$, Young's inequality gives
  \begin{equation*}
    |\partial^k_xP_{1/n}g(x)|\le\|\partial^kp_{n^{-1}\mathbb{I}}\|_{L_{p'}(\Rd)}\|g\|_{L_p(\Rd)},
  \end{equation*}
  where $\frac1p+\frac1{p'}=1$. Applying this with $g=f(t,\cdot)$ on the set $\{\|f(t,\cdot)\|_{L_p}\le n\}$ yields the claim.

  (i) Recall the standard estimate $\|P_{\delta}g-g\|_{\B^{-\varepsilon}_p}\le C\delta^{\varepsilon/2}\|g\|_{\B^0_p}\le C\delta^{\varepsilon/2}\|g\|_{L_p}$ for $\delta\in(0,1]$, $\varepsilon\in(0,2]$ (see \cite[Chapter~2]{BCD11}). Hence, for a.e. $t$,
  \begin{equation*}
    \|f^n(t)-f(t)\|_{\B^{-\varepsilon}_p}\le Cn^{-\varepsilon/2}\|f(t)\|_{L_p}+C\I_{\{\|f(t)\|_{L_p}>n\}}\|f(t)\|_{L_p}.
  \end{equation*}
  The $L_q([0,T])$-norm of the first term tends to $0$. The second term vanishes for $n\ge\|f\|_{L_\infty L_p}$ if $q=\infty$, and its $L_q([0,T])$-norm tends to $0$ by dominated convergence if $q<\infty$.

  (ii) This follows from $\|P_{1/n}g\|_{L_p}\le\|g\|_{L_p}$ and $L_p\hookrightarrow\B^0_p\hookrightarrow\B^{-\varepsilon}_p$.

  (iii) Fix $\eta>0$, put $\Lambda\defby\Lambda_f(\eta)$ and $g\defby|f|\I_{\{|f|>\Lambda\}}$, so that $\|g\|_{L_qL_p}\le\eta$. Since $P_{1/n}$ is a Markov operator, $|f^n|\le P_{1/n}|f|\le P_{1/n}g+\Lambda$ pointwise. On the set $\{|f^n|>2\Lambda\}$, this gives $P_{1/n}g\ge|f^n|-\Lambda>\frac12|f^n|$. Therefore $|f^n|\I_{\{|f^n|>2\Lambda\}}\le2P_{1/n}g$, and
  \begin{equation*}
    \|f^n\I_{\{|f^n|>2\Lambda\}}\|_{L_qL_p}\le2\|P_{1/n}g\|_{L_qL_p}\le2\|g\|_{L_qL_p}\le2\eta,
  \end{equation*}
  that is, $\Lambda_{f^n}(2\eta)\le2\Lambda$.
\end{proof}
\begin{remark}\label{rem:smooth_drift}
  In \cref{sec:strong.existence,sec:stability,sec:flow_and_path_by_path_uniqueness}, a ``smooth function with bounded derivatives'' $b\colon[0,T]\times\Rd\to\Rd$ means a bounded measurable function which is smooth in $x$, with all spatial derivatives bounded uniformly in $(t,x)$. This is all that is used about such drifts (existence and uniqueness of solutions by the Cauchy--Lipschitz theorem, and the mean value theorem in $x$), and by \cref{lem:mollifying}(0) the functions $f^n$ are of this form.
\end{remark}
  
\bibliographystyle{alpha}
\bibliography{ref}
\end{document}